\documentclass[11pt]{article}
\usepackage{amsmath,amsthm,amsfonts,amssymb,amscd}
\usepackage{fancyhdr}
\usepackage{setspace}
\newtheorem{problem}{PROBLEM}
\newtheorem{theorem}{THEOREM}[section]
\theoremstyle{definition}
\newtheorem{corollary}{Corollary}[theorem]
\newtheorem{definition}[theorem]{Definition}
\newtheorem{remark}{Remark}[section]

\newtheorem{lemma}[theorem]{Lemma}
\newtheorem{proposition}[theorem]{Proposition}
\numberwithin{equation}{subsection}
\newcommand{\C}{\mathbf{C}}
\newcommand{\N}{\mathbb{N}}
\newcommand{\Z}{\mathbb{Z}}
\newcommand{\Q}{\mathbb{Q}}
\newcommand{\f}{\mathbf{F}}

\newcommand{\pv}{\mathbf{E}}
\newcommand{\ex}{\mathfrak{e}}

\newcommand{\li}{\mathfrak{l}}
\newcommand{\Pp}{\mathfrak{M}}

\newcommand{\pp}{\mathfrak{m}}
\newcommand{\pq}{\mathfrak{n}}

\newcommand{\mpp}{\mathfrak{p}}

\newcommand{\Ge}{\mathfrak{E}}

\newcommand{\U}{\mathcal{U}}
\newcommand{\pr}{\mathcal{P}}
\newcommand{\R}{\mathcal{R}}
\newcommand{\Dl}{\mathfrak{d}}
\newcommand{\G}{\mathbb{G}}
\newcommand{\pt}{\textbf{P}}
\newcommand{\sg}{\mathbb{H}}
\newcommand{\s}{\mathfrak{s}}

\newcommand{\ls}{\mathcal{S}}
\newcommand{\Cd}{\mathfrak{L}}

\newcommand{\K}{\mathbf{K}}
\newcommand{\bh}{\mathbf{H}}
\newcommand{\M}{\mathbf{M}}
\newcommand{\IL}{\mathbf{L}}
\newcommand{\x}{\mathfrak{x}}
\newcommand{\ft}{\mathfrak{t}}
\newcommand{\z}{\mathfrak{z}}
\newcommand{\y}{\mathfrak{y}}
\newcommand{\rf}{\mathfrak{f}}
\newcommand{\rs}{\mathrm{S}}
\newcommand{\rt}{\mathrm{T}}
\newcommand{\bt}{\begin{theorem}}
\newcommand{\et}{\end{theorem}}
\newcommand{\bco}{\begin{corollary}}
\newcommand{\eco}{\end{corollary}}
\newcommand{\bd}{\begin{definition}}
\newcommand{\ed}{\end{definition}}
\newcommand{\bp}{\begin{problem}}
\newcommand{\ep}{\end{problem}}
\newcommand{\bl}{\begin{lemma}}
\newcommand{\el}{\end{lemma}}
\newcommand{\bprop}{\begin{proposition}}
\newcommand{\eprop}{\end{proposition}}
\newcommand{\br}{\begin{remark}}
\newcommand{\er}{\end{remark}}
\newcommand{\bpf}{\begin{proof}}
\newcommand{\epf}{\end{proof}}
\begin{document}

\title{Extensions by Antiderivatives, Exponentials of Integrals and by Iterated Logarithms}
\author{V. Ravi Srinivasan\footnote{This work is a part of the author's PhD thesis; email: varadhu\_ravi@ou.edu} \\Graduate Student\\ University of Oklahoma}
\maketitle
\begin{abstract}
Let $\f$ be a characteristic zero differential field with an
algebraically closed field of constants $\C$, $\pv\supset\f$ be a
no new constant extension by antiderivatives of $\f$ and let
$\y_1,\cdots\y_n$ be antiderivatives of $\pv$. The antiderivatives
$\y_1,\cdots,\y_n$ of $\pv$ are called J-I-E antiderivatives if
$\y'_i\in\pv$ satisfies certain conditions.  We will discuss a new
proof for the Kolchin-Ostrowski theorem and generalize this
theorem for a tower of extensions by J-I-E antiderivatives and use
this generalized version of the theorem to classify the finitely
differentially generated subfields of this tower. In the process,
we will show that the J-I-E antiderivatives are algebraically
independent over the ground differential field. An example of a
J-I-E tower is iterated antiderivative extensions of the field of
rational functions $\C(x)$ generated by iterated logarithms,
closed at each stage by all (translation) automorphisms. We
analyze the algebraic and differential structure of these
extensions. In particular, we show that the nth iterated
logarithms and their translates are algebraically independent over
the field generated by all lower lever iterated logarithms. Our
analysis provides an algorithm for determining the differential
field generated by any rational expression in iterated logarithms.
These results ultimately rest on the Kolchin--Ostrowski theorem
and it applies to antiderivatives and exponentials of integrals.
Regarding the latter, in Part II of this paper we will present
similar results for iterated exponentials closed under all
(scaling) automorphisms.

\end{abstract}

\tableofcontents \newpage
\section{Introduction}\label{introd}

All the fields considered in this paper are of characteristic
zero. If $\f$ is a field and $':\f\to\f$ a linear map satisfying
the condition $(uv)'=u'v+uv'$ for all $u,v\in\f$ then we will call
the map $'$, a \textsl{derivation} of $\f$. A \textsl{differential
field} is a field $\f$ with a derivation. If $\f$ is a
differential field then one can easily see that
$\C:=\{c\in\f|c'=0\}$ is also a differential field. We will call
$\C$, the field of \textsl{constants} of $\f$. Let $\pv$ and $\f$
be differential fields and let $\pv\supseteq\f$. We say that $\pv$
is a \textsl{differential field extension} of $\f$ if the
derivation of $\pv$ restricted to $\f$ is the derivation of $\f$.
A differential field extension $\pv$ of $\f$  will be called a
\textsl{No New Constants} (NNC) extension of $\f$ if the field of
constants of $\pv$ and $\f$ are the same.

Let $\pv\supset\f$ be a NNC extension. If $\x\in\pv$ and
$\x'\in\f$ then we call $\x$ an \textsl{antiderivative} of an
element (namely, $\x'$) of $\f$, and if $\pv=\f(\x_1\cdots,\x_n)$
for some antiderivatives $\x_1,\cdots,\x_n\in\pv$ of $\f$ then we
will call $\pv$ an \textsl{extension of $\f$ by antiderivatives}.
If $\ex\in\pv$ and $\frac{\ex'}{\ex}\in\f$ then we call $\ex$ an
\textsl{exponential of an integral} of an element (namely,
$\frac{\ex'}{\ex}$) of $\f$, and if $\pv=\f(\ex_1\cdots,\ex_m)$
for some exponentials of integrals $\ex_1\cdots,\ex_m\in\pv$ of
$\f$ then we will call $\pv$\textsl{ an extension of $\f$ by
exponentials of integrals}.

In section \ref{sec kol-ost} we will give a new proof for the
following well known theorem: Let $\f$ be a differential field
with an algebraically closed field of constants $\C$ and let
$\pv\supset\f$ be a NNC extension. Let $\x_1,\cdots,\x_n\in\pv$,
$\ex_1,\cdots,\ex_m\in\pv$ where $\x_i$'s are
antiderivatives($\x'_i\in\f$) and $\ex_i$'s are exponentials of
integrals($\frac{\ex'_i}{\ex_i}\in\f$). Then $\x_1,\cdots,\x_n$,
$\ex_1,\cdots,\ex_m$ are algebraically dependent over $\f$ only if
there are $c_i\in\C$, not all zero, such that $\sum^n_{i=1}
c_i\x_i\in\f$ or there are $n_i\in\Z$, not all zero, such that
$\prod^m_{i=1}\ex^{n_i}_i\in\f$. Thus the algebraic dependence of
$\x_1,\cdots,\x_n$, $\ex_1,\cdots,\ex_m$ over $\f$ becomes a non
trivial linear dependence of $\x_1,\cdots,\x_n$ over $\f$ or there
is a non trivial power product relation among $\ex_1,\cdots,\ex_m$
over $\f$. This theorem is known as the Kolchin-Ostrowski theorem
and it appears as theorem \ref{Kol-Ost} in this paper. A short
note about the history of this theorem is also provided in the
beginning of section \ref{sec kol-ost}.

In section \ref{algo for anti and expo} we will give an algorithm
to compute the differential subfields of an extension $\pv$ of
$\f$ by antiderivatives or by exponentials of integrals. The
extension $\pv$ of $\f$ is assumed to be purely transcendental
over $\f$. Moreover, when $\f$ can be realized as the field of
fractions of a polynomial ring over $\C$ that lives inside $\f$
then for any given intermediate differential subfield
$\pv\supset{\bf K}\supset\f$, our algorithm also computes the
subgroup of differential automorphisms of $\pv$ over $\f$ fixing
${\bf K}$. We complete section \ref{algo for anti and expo} by
proving that if $\pv\supset\f$ is a NNC extension and
$\pv\setminus\f$ contains an antiderivative of $\f$ then there is
an infinite tower of extensions by antiderivatives with the ground
field $\f$ and not imbeddable in any finite tower of
Picard-Vessiot extensions with the ground field $\f$.

We will investigate in section \ref{the jie tower} a special tower
of extensions by antiderivatives, namely the J-I-E tower. We
classify the finitely differentially generated subfields of this
tower. A J-I-E tower exist for any differential field $\f$ that
has a proper antiderivative extension and it may contain
non-elementary functions.

A tower of extensions by iterated logarithms is an example of
J-I-E tower. For a vector $\vec{c}:=(c_1,\cdots,c_n)\in\C^n$,
where $\C$ is an algebraically closed-characteristic zero
differential field with a trivial derivation, we call
$\x[\vec{c},n]:=\log(\log(\cdots\log(x+c_1)\cdots+c_{n-1})+c_n)$
an iterated logarithm of level $n$. In section \ref{the iter
logs}, we give meanings for these iterated logarithms and produce
an algorithm to compute the differential subfields of differential
field extensions by iterated logarithms. In the process, we will
also show that the iterated logarithms are algebraically
independent over $\C(x)$, where $x$ is an element whose derivative
equals 1. In Section \ref{examples} we will provide some examples
of extensions by iterated logarithms and show how our algorithm
works.

{\bf Picard-Vessiot Theory:} Here we will recall some definitions
and state several results from differential Galois theory. One may
find proofs for these results in \cite{Magid 1}. Let $(\f, ')$ be
a differential field with an algebraically closed field of
constants $\C$ and let $\pv$ be any differential field extension
of $\f$. The differential Galois group $\G(\pv|\f)$ is the group
of all differential automorphisms of $\pv$ fixing every element of
$\f$, that is, $\G(\pv|\f):=\{\sigma\in
Aut(\pv|\f)|\sigma(u')=\sigma(u)'\ \forall u\in\pv\}$. Sometimes
we denote $\G(\pv|\f)$ by $\G$ without referring to ground
differential field $\f$ and its extension $\pv$. Let $L(y)$ be a
monic homogeneous linear differential operator of order $n$ over a
differential field $\f$. A differential field extension
$\pv\supseteq \f$ is called a Picard-Vessiot(P-V) extension of
$\f$ for $L(y)$ if the following conditions holds:
\begin{enumerate}
    \item $\pv$ is generated over $\f$ as a differential field by the
    set $V$ of solutions of $L(y)=0$ in $\pv \ (\pv=\f<V>)$
    \item $\pv$ contains a full set of solutions of $L(y)=0$\ (there are
    $y_i\in V, 1\leq i\leq n$, with the wronskian $w(y_1,\cdots,y_n)\ \neq 0$)
    \item Every constant of $\pv$ lies in $\f$.
\end{enumerate}

A Picard-Vessiot extension exists for a given monic homogeneous
linear differential operator $L(y)$ in the case that the field of
constants $\C$ of $\f$ is algebraically closed and it is unique up
to differential automorphisms fixing $\f$. If $\pv$ is a P-V
extension of $\f$ then the set of all elements fixed by the
differential Galois group $\G(\pv|\f)$ is $\f$, that is,
$\pv^{\G}=\{a\in \pv\mid \sigma(a)=a \ \text{for all}\ \sigma\in
\G\}=\f$. The differential Galois group of a P-V extension is an
algebraic matrix group over the field of constants.

If $\pv_i$ is a Picard-Vessiot extension of $\f$ for $1\leq i\leq
n$ then there is a Picard-Vessiot extension $\pv$ of $\f$ such
that $\pv\supseteq\pv_i\supseteq\f$ and $\pv$ is the compositum of
its subfields $\pv_i$.

There is a Fundamental theorem in this context. Let $\f$ be a
differential field with algebraically closed field of constants
$\C$, and let $\pv\supseteq \f$ be a P-V extension. Then the
differential Galois group of $\pv$ over $\f$ is naturally an
algebraic group over $\C$ and there is a lattice inverting
bijective correspondence between
$$ \{\pv\supseteq {\bf K}\supseteq \f\mid {\bf K} \ \text{is an
intermediate differential field}\} $$ and
$$\{ \sg\leq \G\mid \sg \ \text{is a Zariski closed subgroup of} \ \G\} $$
given by $$ {\bf K}\mapsto \G(\pv|{\bf K}) \ \text{and}\
\sg\mapsto \pv^\sg.$$ The intermediate field ${\bf K}$ is a P-V
extension of $\f$ if and only if the subgroup $\sg=\G(\pv|{\bf
K})$ is normal in $\G$; if it is, then $$
\G(\pv^\sg|\f)=\frac{\G}{\sg}.$$ Let $\G^0$ be the connected
component of the identity in $\G(\pv|\f)$, and let $\pv^0$ be the
corresponding intermediate field. Then $\pv^0$ is the algebraic
closure of $\f$ in $\pv$, $\pv^0$ is a finite Galois extension of
$\f$ with Galois group $\frac{\G(\pv|\f)}{\G^0(\pv|\f)}$, and the
transcendence degree of $\pv$ over $\pv^0$ is dim$(\G^0(\pv|\f))$.

Analogous to the algebraic closure of a given field, we may define
a Picard-Vessiot closure of a given differential field $\f$. The
Picard-Vessiot closure $\f_1$ of $\f_0:=\f$ is a differential
field extension of $\f_0$ such that
\begin{itemize}
    \item $\f_1$ is a union of Picard-Vessiot extensions of $\f_0$
    \item Every Picard-Vessiot extension of $\f_0$ has an isomorphic
    copy in $\f_1$.
\end{itemize}
The Picard-Vessiot closure $\f_1$ of $\f_0$ need not be
``closed''. That is, there are linear homogeneous differential
equations over $\f_1$ whose solutions may not be in $\f_1$ (see
theorem \ref{simple version}). This leads us to consider a chain
of Picard-Vessiot closures of $\f_0$. A finite tower of
Picard-Vessiot closures of $\f_0$ is a chain
$$\f_0\subseteq\f_1\subseteq\f_2\subseteq\cdots\subseteq\f_n,$$
where $\f_0:=\f, n\in\N$ and $\f_i$ is the Picard-Vessiot closure
of $\f_{i-1}$, for all $1\leq i\leq n$. Finally we define the
complete Picard-Vessiot closure $\f_\infty$ of $\f$ as the union
$\cup^\infty_{i=0}\f_i$. The differential field $\f_\infty$ is
``closed''. If $\pv$ is a normal differential subfield of
$\f_\infty$ then every automorphism of $\phi\in\G(\pv|\f)$ extends
to an automorphism $\Phi\in\G(\f_\infty|\f)$ and every
automorphism $\Phi\in\G(\f_\infty|\f)$ also restricts to a
$\phi\in\G(\pv|\f)$. We also note that the fixed field of
$\G(\f_\infty|\f)$ is $\f$. For details see \cite{Magid 3}.

\section{The Kolchin-Ostrowski Theorem}\label{sec kol-ost}

Throughout this paper, $\f$ denotes a characteristic zero
differential field with an algebraically closed field of constants
$\C$. Sometimes we will denote the field of constants $\C$ of $\f$
by $\C_\f$. Let us recall some definitions from section
\ref{introd}

\bd Let $\pv\supset\f$ be a differential field extension of $\f$.
An element $\x\in\pv$ is called an \textsl{antiderivative} of an
element of $\f$ if $\x'\in\f$. A No New Constant(NNC) extension
$\pv\supset\f$ is called an \textsl{extension by antiderivatives}
of $\f$ if for $i=1,2,\cdots,n$ there exists $\x_i\in\pv$ such
that $\x'_i\in\f$ and $\pv=\f(\x_1,\x_2,\cdots,\x_n)$. \ed

\bd

Let $\pv\subset\f$ be a differential field extension of $\f$. An
element $\ex\in\pv$ is called an \textsl{exponential of an
integral} of an element of $\f$ if $\frac{\ex'}{\ex}\in\f$. A NNC
extension $\pv\supset\f$ is an extension by \textsl{exponential of
integrals} of $\f$ if for $i=1,2,\cdots,n$ there exists
$\ex_i\in\pv$ such that $\frac{\ex'_i}{\ex_i}\in\f$ and
$\pv=\f(\ex_1,\ex_2,\cdots,\ex_n)$.

\ed

In this section we will prove the Kolchin-Ostrowski theorem, which
states

\bt\label{Kol-Ost}(\textsf{Kolchin-Ostrowski}) Let $\pv\supset\f$
be a NNC differential field extension and let
$\x_1,\cdots,\x_n\in\pv$, $\ex_1,\cdots,\ex_m\in\pv\setminus\{0\}$
be such that $\x_i$ is an antiderivative of an element $\f$ for
each $i$ ($\x'_i\in\f$) and $\ex_i$ is an exponential of an
integral of an element of $\f$ for each $i$
($\frac{\ex'_j}{\ex_j}\in\f$). Then either
$\x_1,\cdots,\x_n$,$\ex_1,\cdots,\ex_m$ are algebraically
independent over $\f$ or there exist
$(c_1,\cdots,c_n)\in\C^n\setminus\{0\}$ such that
$\sum^n_{i=1}c_i\x_i\in\f$ or there exist
$(r_1,\cdots,r_m)\in\Z^n\setminus\{0\}$ such that
$\prod^m_{j=1}\ex^{r_j}_j\in\f$. \et

\subsection{Algebraic Dependence of Antiderivatives.}
In his paper \cite{Ostrowski}, A. Ostrowski proves that a set of
antiderivatives $\{\x_1,\cdots$, $\x_n\}$ of $\f$ is either
algebraically independent over $\f$ or there are constants
$c_i\in\C$ not all zero such that $\sum^n_{i=1}c_i\x_i\in\f$. In
his setting, $\f$ is a differential field of meromorphic functions
and $\C=\mathbb{C}$, the field of complex numbers. Later,
Ostrowski's result was generalized by Kolchin \cite{Kolchin} to
theorem \ref{Kol-Ost}. In their papers \cite{j.ax} and
\cite{M.Ros}, J. Ax and M. Rosenlicht also presented proofs of
theorem \ref{Kol-Ost}. The proof we are going to present is
elementary and differ from the proofs listed above.

\bt\label{anti elem} Let $\pv\supset\f$ be a differential field
extension and let $\x\in\pv$ be an antiderivative. Then either
$\x$ is transcendental over $\f$ or $\x\in\f$. \et

\bpf Let $\C_\f$ denote the field of constants of $\f$ and suppose
that $\x$ is algebraic over $\f$. Then there is a monic
irreducible polynomial $P(x)=\sum^n_{i=0}a_ix^i\in\f[x]$ such that
$P(\x)=0$. Note that $(P(\x))'=0$, that is $\x$ is a solution of
the polynomial
$$\sum^{n}_{i=1}(ia_i\x'-a'_{i-1})x^{i-1}\in\f[x].$$
Since the degree of the above polynomial $< n$, it has to be the
zero polynomial. In particular $n\x'=a'_{n-1}$, that is
$(\x-b)'=0$, where $b:=\frac{a_{n-1}}{n}\in\f$. Observe that
$\x-b$ is algebraic over $\f$ (since $\x$ and $b$ are algebraic)
and therefore there is a monic irreducible polynomial
$Q(x)=\sum^m_{i=0}b_ix^i\in\f[x]$ such that $Q(\x-b)=0$. Again
taking the derivative of the equation $Q(\x-b)=0$, we note that
$\x-b$ is a solution of the polynomial
$$\sum^{m}_{i=1}b'_{i-1}x^{i-1}\in\f[x].$$ Since the degree of the
above polynomial is $<m$, it has to be the zero polynomial. Thus
$b_i\in\C_\f$ and therefore the polynomial $Q(x)$ has coefficients
in $\C$. Since $\C_\f$ is algebraically closed and $\x-b$ is a
zero of $Q(x)$ we obtain $\x-b\in\C_\f$. Now $b\in\f$ will imply
that $\x=b+c\in\f$, where $c:=\x-b$.

Note that we do not require the constants of $\f$ and $\pv$ to be
the same to prove this theorem. The above theorem is also proved
in \cite{Kaplansky}, page 23 and \cite{Magid 1}, page 7.\epf

Let $\pv\supseteq\f$ be an extension by antiderivatives
$\x_1,\cdots,\x_n\in\pv\setminus\C$ of $\f$ . That is,
$\pv=\f(\x_1,\cdots,\x_n)$, $\x'_i\in\f$ and $\x_i\notin\C$ for
all $1\leq i\leq n$. Since $\pv$ is a NNC extension of $\f$, the
differential subfield $\pv_i=\f(\x_i)$ of $\pv$ is also a NNC
extension of $\f$. Let $f_i:=\x'_i\in\f$  and observe that
$$\x''_i=\frac{f'_i}{f_i}\x'_i.$$
Thus $\x_i$ is a solution of a second order linear homogeneous
differential equation over $\f$. Moreover, if $V_i$ is the vector
space spanned by the unity $1\in\C$ and $\x_i$ over $\C$ then
$\pv_i=\f\langle V_i\rangle$--the differential field generated by
$\f$ and $V_i$. The full set of solutions of the differential
equation $Y''=\frac{f'}{f}Y'$ is the vector space $V_i$. Thus we
see that $\pv_i$ is a Picard-Vessiot extension of $\f$. Since a
compositum of Picard-Vessiot extensions is again a Picard-Vessiot
extension(see \cite{Magid 1}, page 28-29),
$\pv:=\pv_1\cdot\pv_2\cdots\pv_n$ is also a Picard-Vessiot
extension of $\f$.

Assume that $\x_i\notin\f$ for each $i$. If
$\sigma\in\G(\pv_i|\f)$ then
\begin{equation}\sigma(\x_i)'=\sigma(\x'_i)=\sigma(f_i)=f_i=\x'_i.\end{equation}
Thus $\sigma(\x_i)'=\x'_i$, which implies
$\big(\sigma(\x_i)-\x_i\big)'=0$. Since $\pv$ is a NNC extension
of $\f$, there is a $c_{i\sigma}\in\C$ such that
$\sigma(\x_i)-\x_i=c_{i\sigma}$, that is,
$\sigma(\x_i)=\x_i+c_{i\sigma}$. On the other hand, for any
$c\in\C$, the automorphism $\sigma_{ic}:\pv_i\to\pv_i$ defined as
$\sigma_{ic}(\x_i)=\x_i+c$ and $\sigma(f)=f$ for all $f\in\f$ can
be readily seen as a differential automorphism. Thus
$\G(\pv_i|\f)$ injects into $(\C,+)$ as an algebraic subgroup for
each $i$. Note that $(\C,+)$ has no non trivial algebraic
subgroups and since $\x_i\notin\f$, from the fundamental theorem,
we see that $\G(\pv_i|\f)\simeq(\C,+)$ and that the extension
$\pv_i$ of $\f$ has no intermediate differential subfields.  Any
automorphism of $\pv$ fixing $\f$ is completely determined by its
action on $\x_1,\cdots,\x_n$ and thus we have a map
$\sigma\mapsto(c_{1\sigma},\cdots,c_{n\sigma})$, an algebraic
group homomorphism from $\G$ to $(\C,+)^n$. This map is clearly
injective. From this observation, we see that the differential
Galois group $\G(\pv|\f)$ is isomorphic to an algebraic subgroup
of $(\C,+)^n$. Note that $\G(\pv|\f)$ could be a proper algebraic
subgroup of $(\C,+)^n$; depending on whether all the
antiderivatives are algebraically independent over $\f$ or not. We
will discuss about the nature of the algebraic dependence of
antiderivatives in the next theorem.

We will do a similar analysis for the extensions by exponentials
of integrals of $\f$ in subsection \ref{expo iteg}.

\bt \label{transcen}Let $\pv\supset\f$ be a NNC differential field
extension and for $1=1,2,\cdots,n$ let $\x_i\in\pv$ be
antiderivatives of $\f$. Then either $\x_i$'s are algebraically
independent over $\f$ or there is a tuple
$(c_1,\cdots,c_n)\in\C^n\setminus\{0\}$ such that
$\sum^n_{i=1}c_i\x_i\in\f$.\et

\bpf[Proof 1.] First we will present Kolchin's proof. Observe that
$\pv=\f(\x_1,\x_2,\cdots,\x_n)$ is a Picard-Vessiot extension of
$\f$ and for every $\sigma\in\G(\pv|\f)$ we see that
$\sigma(\x_i)=\x_i+c_{i\sigma}$. Thus, as noted earlier,
$\G(\pv|\f)$ imbeds into $(\C^n, +)$ as an algebraic subgroup.
Suppose that the $\x_i's$ are algebraically dependent and say
$\x_1$ is algebraic over $\f(\x_2,\x_3,\cdots,\x_n)$. We may also
assume that $\x_i$'s$\notin\f$ for any $i$ (otherwise there is
nothing to prove).

Since $\x_1$ is an antiderivative of an element of $\f$ and $\x_1$
is algebraic over $\f(\x_2,\x_3,\cdots,\x_n)$ from theorem
\ref{anti elem} we obtain $\x_1\in\f(\x_2,\x_3,\cdots,\x_n)$ and
thus $\G(\pv|\f)\hookrightarrow (\C^n, +)$ is not a surjection. In
particular, if $\sigma\in\G(\pv|\f)$ fixes $\x_2,\cdots,\x_n$ then
$\sigma$ fixes $\x_1$ too. Therefore
$$\G(\pv|\f)=\{(d_1,d_2,\cdots,d_n)\in\C^n|L_i(d_1,d_2,\cdots,d_n)=0, 1\leq i\leq t\},$$
where $L_i$ is a linear form over $\C$ for each $i$. Now for any
$\sigma\in\G(\pv|\f)$ and $L\in\{L_i|1\leq i\leq t\}$,
\begin{align*}\sigma(L(\x_1,\x_2,\cdots,\x_n))&=L(\sigma(\x_1),\sigma(\x_2),\cdots,\sigma(\x_n))\\
&=L(\x_1+d_1,\x_2+d_2,\cdots,\x_n+d_n)\\
&=L(\x_1,\x_2,\cdots,\x_n)+L(d_1,d_2,\cdots,d_n)\\
&=L(\x_1,\x_2,\cdots,\x_n)\quad\text{since}\ L(d_1,\cdots,d_n)=0
\end{align*}
and thus $L(\x_1,\x_2,\cdots,\x_n)\in\pv^{\G(\pv|\f)}$. From
Galois theory we know that $\pv^{\G(\pv|\f)}=\f$. Since $L$ is a
linear form over $\C$, we obtain $L(\x_1,\x_2,\cdots,\x_n)$
$=\sum^n_{i=1}c_i\x_i\in\f$. \epf

\bpf[Proof 2.]This proof does not require Galois theory. For every
tuple $(c_1,\cdots,c_n)$ $\in\C^n\setminus\{0\}$ let us assume
that $\sum^n_{i=1}c_i\x_i\notin\f$ . Theorem \ref{anti elem} and
our assumption that $\sum^n_{i=1}c_i\x_i\notin\f$ guarantees us a
nonempty algebraically independent subset $S$ of $\{\x_i|1\leq
i\leq n\}$ over $\f$. We may assume that
$S=\{\x_2,\x_3,\cdots,\x_n\}$. Again from theorem \ref{anti elem},
we see that $\x_1$ is transcendental over $\f(S)$ or
$\x_1\in\f(S)$. We will show that the latter case is not possible
and this will prove the theorem.

Suppose that $\x_1\in\f(S)$ and let $t$ be the largest positive
integer such that
$$\sum^t_{i=1} c_i\x_i\in\f(S_t),$$
where $c_i\in\C, c_1=1$ and $S_t:=S\setminus \{\x_i\big|2\leq
i\leq t\}$.

Since $|S|<\infty$ and $t\geq 1$, such a $t$ exist and since
$\sum^n_{i=1}c_i\x_i\notin\f$, $S_t\neq\emptyset$. In particular,
$t<n$ and thus $\x_{t+1}\in S_t$. For notational convenience let
$\x:=\x_{t+1}$. Then

$$\sum^t_{i=1} c_i\x_i=\frac{P}{Q}$$ where $P:=\sum^r_{i=0} a_i
\x^i$, $Q:=\sum^s_{i=0} b_i \x^i$, $b_s=1$, $a_r\neq 0, a_i,b_i\in
\mathbf{K}:=\f(S_t\setminus\{\x\})$ and $(P,Q)=1$. Differentiating
the above equation, we get $\sum^t_{i=1}
c_i\x'_i=\frac{P'Q-PQ'}{Q^2}$ and thus
\begin{equation}\label{std deriv eqn} fQ^2=P'Q-PQ',\end{equation}
where $f:=\sum^t_{i=1} c_i\x'_i$. If $f=0$ then $(\sum^t_{i=1}
c_i\x_i)'=0$ and since $\pv$ is a NNC extension of $\f$,
$\sum^t_{i=1} c_i\x_i\in\C\subset\f$, a contradiction to our
assumption that $\sum^t_{i=1} c_i\x_i\notin\f$. Thus $f\neq 0$.
Now suppose that deg$Q\geq 1$. From the above equation we see that
$Q$ divides $P'Q-PQ'$, which implies $Q$ divides $PQ'$ and since
$(P,Q)=1$, $Q$ divides $Q'$. Thus $s=$deg$Q\leq$ deg$Q'$. But then
deg$Q'$=deg$((s\x'+b'_{s-1})\x^{s-1}+\cdots+b_1\x'+b'_0)$ $\leq
s-1$, a contradiction. Thus deg$Q=0$, that is $Q\in\mathbf{K}$.

Hence we may assume that $\sum^t_{i=1} c_i\x_i=P$ and note that
\begin{equation}\label{kol eqn}f=P'.\end{equation}

Case 1: deg$(P)=0$, that is
$P\in\mathbf{K}=\f(S_t\setminus\{\x\})$.

Then  $\sum^t_{i=1} c_i\x_i=P\in\f(S_t\setminus\{\x\})$. Since
$\x=\x_{t+1}$, we obtain $\sum^{t+1}_{i=1}
c_i\x_i\in\f(S_t\setminus\{\x_{t+1}\})$, where $c_{t+1}:=0$. This
contradicts the maximality of $t$.

Case 2: deg$(P)>1$

From equation \ref{kol eqn} we see that
\begin{equation}\label{kol
main eqn}f=a'_r\x^r+(ra_r\x'+
a'_{r-1})\x^{r-1}+\cdots+a_1\x'+a'_0.
\end{equation}
Thus comparing the coefficients of $\x^r$ we get $a'_r=0$, that is
$a_r\in\C$. Since $r-1\geq 1$  comparing the coefficients of
$\x^{r-1}$, we get
\begin{align*}&ra_r\x'+a'_{r-1}=0\\
\implies&\x'=\big(\frac{-a_{r-1}}{ra_r}\big)'\\
\implies&\x=\frac{-a_{r-1}}{ra_r}+c_1 \end{align*} for some
$c_1\in\C$ and thus $\x=\frac{-a_{r-1}}{ra_r}+c_1\in\mathbf{K}$, a
contradiction to the assumption that $\x$ is transcendental over
$\mathbf{K}$.

Case 3: deg$P=1$

Finally if deg$P=1$ then $P=a_1\x+a_0=\sum^{t}_{i=1} c_i\x_i$ and
therefore taking the derivative we have
$$
a'_1\x+a_1\x'+a'_0=f.$$ Thus comparing the coefficients, we obtain
$a'_1=0$ that is $a_1\in\C$ and $a_1\x'+a'_0=f$. Now letting
$c_{t+1}:=-a_1$ and substituting $\x_{t+1}$ for $\x$, we get
$\sum^{t+1}_{i=1}
c_i\x_i=a_0\in\mathbf{K}=\f(S_t\setminus\{\x_{t+1}\})$ and this
again contradicts the maximality of $t$. Hence the theorem. \epf

\subsection{Exponentials of an Integrals.}\label{expo iteg}

Here we will prove theorems analogous to theorems \ref{anti elem}
and \ref{transcen} for the exponential of an integral setting.

\bt\label{expo elem}

Let $\pv\supset\f$ be a differential field extension. If there is
a $\ex\in\pv$ such that $\frac{\ex'}{\ex}\in\f$ then either $\ex$
is transcendental over $\f$ or there is an $n\in\N$ such that
$\ex^n\in\f$.\et

\bpf  Suppose that $\ex$ is algebraic over $\f$,
$\frac{\ex'}{\ex}=f\in\f$ and let
$P(x)=\sum^n_{i=0}a_ix^i\in\f[x]$ be the monic irreducible
polynomial of $\ex$. Then $P(\ex)=0$ and therefore $(P(\ex))'=0$,
which implies $\ex$ is a solution of the polynomial
$$P_1:=nfx^n+\sum^{n-1}_{i=0}(a'_{i}-ia_{i}f)x^i\in\f[x].$$
Since $P$ is the monic irreducible polynomial of $\ex$, we have
$nfP=P_1$. Thus comparing the coefficients of $nfP$ and $P_1$ we
obtain $nfa_0=a'_0$ and since $nf\ex^n=(\ex^n)'$, we obtain
$\big(\frac{\ex^n}{a_0}\big)'=0$ (P is irreducible so $a_0\neq
0$). Note that $\ex$ and $a_0$ are algebraic over $\f$ so
$\frac{\ex^n}{a_0}$ is also algebraic over $\f$. Since
$(\frac{\ex^n}{a_0})'=0$, as in the proof of theorem \ref{anti
elem}, we obtain $\frac{\ex^n}{a_0}=c\in\C_\f$ and thus
$\ex^n=ca_0\in\f$.

This theorem is also proved in \cite{Kaplansky}, page 24 and
\cite{Magid 1}, page 8. \epf

\bt\label{alg indp}

Let $\pv\supset\f$ be a NNC differential field extension and for
$i=1,2,\cdots,n$ let $\ex_i\in\pv\setminus\{0\}$ be such that
$\frac{\ex'}{\ex}\in\f$. Then either $\ex_1,\cdots,\ex_n$ are
algebraically independent or there exist
$(k_1,\cdots,k_n)\in\Z^n\setminus\{0\}$ such that the power
product $\prod^n_{i=1}\ex^{k_i}_i\in\f$.

\et

\bpf The proof of this theorem very much mimics the proof of
theorem \ref{transcen}. Let us assume that
$\prod^n_{i=1}\ex^{k_i}_i\notin\f$ for any
$(k_1,\cdots,k_n)\in\Z^n\setminus\{0\}$. Then from theorem
\ref{expo elem} we see that there is a nonempty algebraically
independent set $S\subset\{\ex_i|1\leq i\leq n\}$ and we may
assume that $S=\{\ex_2,\cdots,\ex_n\}$. From theorem \ref{expo
elem} we see that either $\ex_1$ is transcendental over $\f(S)$ or
there is a $k_1\in\N$ such that $\ex^{k_1}_1\in\f(S)$. We will
show that the latter is not possible and this will prove the
theorem.

Suppose that there is a $k_1\in\N$ such that
$\ex^{k_1}_1\in\f(S)$. Let $t$ be the largest positive integer
such that the power product
$$\prod^t_{i=1}\ex^{k_i}_i\in\f(S_t),$$

where $k_i\in\Z$ for $2\leq i\leq t$ and
$S_t=S\setminus\{\ex_i|2\leq i\leq t\}$.  Since
$\prod^n_{i=1}\ex^{k_i}_i\notin\f$ we obtain $S_t\neq\emptyset$.
Indeed $\ex_{t+1}\in S_t$. Let $\ex:=\ex_{t+1}$ and write

$$\prod^t_{i=1}\ex^{k_i}_i=\frac{P}{Q}$$ where $P:=\sum^l_{i=0}
a_i \ex^i $, $Q:=\sum^m_{i=0}b_i \ex^i$, $(P,Q)=1$ $b_m=1$,
$a_l\neq 0, a_i,b_i\in\f(S_t\setminus\{\ex\})$. Differentiating
the above equation, we get
$$\big(\prod^t_{i=1}\ex^{k_i}_i\big)'=\frac{P'Q-PQ'}{Q^2},$$

Let $f_i:=\frac{\ex'_i}{\ex_i}$, $g:=\frac{\ex'}{\ex}$,
$P=\sum^l_{i=0} a_i \ex^i$ and $Q=\sum^m_{i=0} b_i \ex^i$. Note
that $g,f_i\in\f$ and
\begin{align*}
(\prod^t_{i=1}\ex^{k_i}_i)'&=\sum^t_{j=1}(\ex^{k_j}_j)'
\prod^t_{i=1, i\neq j}\ex^{k_i}_i\\
&=\sum^t_{j=1}k_j\ex'_j\ex^{k_j-1}_j \prod^t_{i=1, i\neq
j}\ex^{k_i}_i,
\end{align*}
which implies
\begin{equation}\label{prod deriv}(\prod^t_{i=1}\ex^{k_i}_i)'=\big(\sum^t_{i=1}
k_jf_j\big)\prod^t_{i=1} \ex^{k_i}_i\end{equation} and thus
$\Big(\frac{P}{Q}\Big)'=\big(\sum^t_{i=1}k_jf_j\big)\frac{P}{Q}.$
Hence \begin{equation}\label{poly expo 2}
QP'-PQ'=\big(\sum^t_{i=1}k_jf_j\big)PQ.
\end{equation}
Since
\begin{align*}QP'-PQ'&=((a'_l+la_lg)\ex^{l+m}+\cdots+a'_0b_0)\\
&-(ma_lg\ex^{l+m}+\cdots+a_0b'_0)\\
&=(a'_l+(l-m)a_lg)\ex^{l+m}+\cdots+a'_0b_0-a_0b'_0,
\end{align*} and
$$PQ=a_l\ex^{l+m}+(a_lb_{m-1}+a_{l-1})\ex^{l+m-1}+\cdots+a_0b_0,$$
substituting in equation \ref{poly expo 2} we get
\begin{align*}
&(a'_l+(l-m)a_lg)\ex^{l+m}+\cdots+a'_0b_0-a_0b'_0=\big(\sum^t_{i=1}k_jf_j\big)
(a_l\ex^{l+m}\\&+(a_lb_{m-1}+a_{l-1})\ex^{l+m-1}+\cdots+a_0b_0).
\end{align*}
The LHS and RHS are polynomial in $\ex$ with coefficients in
$\f(S_t\setminus\{\ex\})$. Since $\pv\supset\f$ is a NNC extension
and $\prod^t_{i=1}\ex^{k_i}_i\notin\f$ we have
$\sum^t_{i=1}k_jf_j\neq 0$ and therefore both the LHS and RHS are
of degree $l+m$. Thus comparing the coefficients of $\ex^{l+m}$ we
get
\begin{align*}
a'_l&+(l-m)a_lg=(\sum^t_{i=1}k_jf_j) a_l\\
\implies a'_l&=[(\sum^t_{i=1}k_jf_j)+(m-l)g] a_l.
\end{align*}

We observe that

\begin{equation}\label{expo ode 1}a'_l=(\sum^{t+1}_{i=1}k_jf_j) a_l,\end{equation}
where $k_{t+1}:=m-l$ and $f_{t+1}:=g$.

We also know that $\prod^{t+1}_{i=1}\ex^{k_i}_i$ is also a nonzero
solution of the equation \ref{expo ode 1} and therefore
$\Big(\frac{\prod^{t+1}_{i=1}\ex^{k_i}_i}{a_l}\Big)'=0$. Since
$\pv$ and $\f$ have the same field of constants, there is an
$\alpha\in\C\setminus\{0\}$ such that
$\prod^{t+1}_{i=1}\ex^{k_i}_i=\alpha a_l$. Now
$a_l\in\f(S_t\setminus\{\ex\})$ will imply
$\prod^{t+1}_{i=1}\ex^{k_i}_i\in\f(S_t\setminus\{\ex_{t+1}\})$, a
contradiction to the maximality of $t$. Hence the theorem. \epf

{\bf The Kolchin-Ostrowski Theorem}

\bpf[Proof of theorem \ref{Kol-Ost}]

Let us assume that $\x_1,\cdots,\x_n,\ex_1,\cdots,\ex_m$ are
algebraically dependent over $\f$ and also that
$\ex_1,\cdots,\ex_m$ are algebraically independent over $\f$.
(Note that if $\ex_1,\cdots,\ex_m$ are algebraically dependent
over $\f$ we may apply theorem \ref{alg indp} to prove this
theorem.) Let us prove that there are constants $c_i\in\C$ not all
zero such that $\sum^n_{i=1}c_i\x_i\in\f$.

It is clear from our assumption that $\x_1,\cdots,\x_n$ is
algebraically dependent over $\mathbf{K}:=\f(\ex_1,\cdots,\ex_m)$.
Since $\x_1,\cdots,\x_n$ are antiderivatives of $\f$ they are also
antiderivatives of $\mathbf{K}$ and thus theorem \ref{transcen} is
applicable with $\mathbf{K}$ as the ground field. Thus there are
constants $c_i\in\C$ not all zero such that
$\sum^n_{i=1}c_i\x_i\in\mathbf{K}$. Choose a subset
$S\subset\{\ex_1,\cdots,\ex_m\}$ so that
$\sum^n_{i=1}c_i\x_i\in\f(S)$ but not in any of the subfields
$\f(S_1)$, where $S_1$ is a proper subset of $S$.

We claim that $S=\emptyset$ and this will prove that
$\sum^n_{i=1}c_i\x_i\in\f$. Suppose not. Then there is a $\ex\in
S$ and we may write
\begin{equation}\label{kol ost eqn}\sum^n_{i=1}c_i\x_i=\frac{P}{Q},\end{equation} where
$P,Q\in\f(S\setminus\{\ex\})[\ex]$, $(P,Q)=1$ and $Q$ a monic
polynomial. Let $f=(\sum^n_{i=1}c_i\x_i)'$. Note that $f\in\f$ and
if $f=0$ then $(\sum^n_{i=1}c_i\x_i)'=0$ and since the extensions
are NNC, we see that $\sum^n_{i=1}c_i\x_i=\alpha\in\C\subset\f$
and we are done. So we assume $f\neq 0$ and note that this
condition also says that $P\neq 0$. Now Differentiating the
equation \ref{kol ost eqn} we obtain
\begin{equation}\label{kol ost prim eq}fQ^2=P'Q-Q'P.\end{equation}

Hereafter one can complete the proof by precisely following the
part of the proof of theorem \ref{transcen} that follows after
equation \ref{std deriv eqn}. Here I will give an alternate
argument which is also applicable for the part of the proof of
theorem \ref{transcen} that follows after equation \ref{std deriv
eqn}.

Note that deg$(P'Q-Q'P)\leq r+s$ and deg$(fQ^2)=$deg$(Q^2)=2s$.

Case 1: deg$Q>$deg$P$.

In this case we see that $r+s<deg(Q^2)=2s$. Since the leading
coefficient $f$ of the LHS of \ref{kol ost prim eq} is nonzero, we
obtain that $\ex$ is algebraic over $\f(S\setminus\{\ex\})$, a
contradiction.

Case 2: deg$(Q)<$deg$(P)$

Let $\frac{\ex'}{\ex}=g\in\f$, $P=\sum^r_{i=0}a_i\ex^i$, $a_r\neq
0$, $Q=\sum^s_{i=0}b_i\ex^i$ and $b_s=1$. Note that
$P'Q-Q'P=(a'_r-(r-s)a_rg)\ex^{r+s}+\cdots$. If
$(a'_r-(r-s)a_rg)\neq 0$ then $r+s=deg(P'Q-Q'P)$ and since $s<r$,
deg$(Q^2)=2s<r+s$, which implies $\ex$ is algebraic over
$\f(S\setminus\{\ex\})$, a contradiction to our assumption that
$\ex_i$'s are algebraically independent over $\f$. Thus
$a'_r-(r-s)a_rg=0$, that is $a'_r=(r-s)ga_r$. Note that $a_r\neq
0$ and since $(\ex^{r-s})'=(r-s)g\ex^{r-s}$ and $r\neq s$, there
is a constant $\alpha\in\C\setminus\{0\}$ such that
$\ex^{r-s}=\alpha a_r\in\f(S\setminus\{\ex\})$ again contradicting
the algebraic independency of $\ex_i$'s over $\f$.

Case 3: deg$(P)=$deg$(Q)$

Since deg$(Q^2)=2s$, deg$(P'Q-Q'P)\leq 2s$ and $f\neq 0$, we have
$f=a'_r-(r-s)ga_r$ and this equation further reduces to $f=a'_r$
since $r=s$. Now the facts $(\sum^n_{i=1}c_i\x_i)'=f$ and
$\mathbf{K}$ is a NNC extension together will imply that
$\sum^n_{i=1}c_i\x_i=a_r+\alpha$ for some $\alpha\in\C$. Thus
$\sum^n_{i=1}c_i\x_i=a_r+\alpha\in\f(S\setminus\{\ex\})$, a
contradiction to the minimality of $S$.

Thus $S$ has to be the empty set and hence the theorem. \epf

\section{Extensions by antiderivatives and by exponentials of
integrals}\label{algo for anti and expo}

Let $\pv\supset\f$ be an extension by antiderivatives
$\x_1,\cdots,\x_n$ of $\f$. We know from theorem \ref{transcen}
that the set of antiderivatives $\{\x_i|1\leq i\leq n\}$ is either
algebraically independent or there are constants $c_i\in\C$ not
all zero such that $\sum^n_{i=1}c_i\x_i\in\f$. Also note that if
$\x_1,\cdots, \x_n$ is algebraically dependent over $\f$ then we
may chose a transcendence base $S\subset\{\x_1,\cdots,\x_n\}$ of
$\pv$ over $\f$ and this makes $\pv$ algebraic over $\f(S)$. But
then each $\x\in\{\x_1,\cdots,\x_n\}\setminus S$ becomes algebraic
over $\f(S)$ and therefore from theorem \ref{anti elem} we obtain
$\x\in\f(S)$ which implies $\pv=\f(S)$. In other words extensions
by antiderivatives are purely transcendental. Thus, to study an
extension by antiderivatives $\x_1,\cdots,\x_n$ of $\f$, we may
very well assume that $\x_1\cdots,\x_n$ are algebraically
independent over $\f$.

In this section we will prove the following theorem

\begin{theorem}\label{struct theorem for antideriv}

Let $\pv=\f(\x_1,\cdots,\x_n)$ be an extension by antiderivatives
$\x_1,\cdots, \x_n$ of $\f$ and let $\x_1,\cdots, \x_n$ be
algebraically independent over $\f$. Let $u\in\pv$ and
$u=\frac{P}{Q}$, $P,Q\in\f[\x_1,\cdots,\x_n]$ and $(P,Q)=1$. Then
there is a $t\in\N$ and $\f-linear$ forms $D_i\in
Span_\f\{\x_1,\cdots,\x_n\}$ for $1\leq i\leq t$ such that
$$\f\langle u\rangle=\f(D_i|1\leq i\leq t).$$ Moreover these
linear forms $D_i$ can be explicitly computed for $P$ and $Q$.

\end{theorem}

A much stronger result can be obtained using Galois theory and
that is, if ${\bf K}$ is an intermediate differential subfield of
$\pv|\f$ then
\begin{equation} \label{Galois stru of anti}{\bf K}=\f(L_i|1\leq i\leq t),
\end{equation} where the linear forms are over $\C$. That is $L_i\in Span_
\C\{\x_1\cdots,\x_n\}$. This follows immediately from the
following three facts 1. The extension $\pv\supset\f$ is a P-V
extension with a differential Galois group $(\C,+)^n$. 2. There is
a bijective correspondence between the algebraic subgroups of
$(\C,+)^n$ and the intermediate differential subfields of
$\pv|\f$; see the fundamental theorem stated in section
\ref{introd}. 3. The algebraic subgroups of $(\C,+)^n$ are
solution sets of linear forms over $\C$.

Though we know the structure of intermediate differential
subfields of $\pv|\f$ it is not clear how to obtain those linear
forms for a given intermediate differential subfield. The theorem
\ref{struct theorem for antideriv} shows that there is a way to
figure out linear forms(not over $\C$ but over $\f$) for singly
differentially generated subfields of $\pv$ containing $\f$ and
since a finitely differentially generated subfield is a compositum
of singly differentially generated subfields of $\pv$ containing
$\f$, we may generalize the theorem \ref{struct theorem for
antideriv} for any finitely generated differential subfield of
$\pv$ containing $\f$. We will prove a similar result for
extensions by exponentials of integrals and will also prove a
similar structure theorem for NNC extensions of the form
$\f(\x_1,\cdots,\x_n,\ex_1,\cdots,\ex_m)$, where $\x'_i\in\f$ and
$\frac{\ex'_i}{\ex_i}\in\f$ and
$\x_1,\cdots,\x_n,\ex_1,\cdots,\ex_m$ are algebraically
independent over $\f$.

To prove theorem \ref{struct theorem for antideriv} we need some
results about several variable polynomials over a commutative ring
with unity, which will be dealt in the following section.

\subsection{Multivariable Taylor formula}

Let $\R$ be an integral domain with $\Q\subseteq\R$  and let
$\R[y_1,\cdots,y_n]$ be the polynomial ring over
$n-$indeterminates $y_1,\cdots,y_n$. Let
$P:=P(y_1,\cdots,y_n)\in\R[y_1,\cdots,y_n]$,
$(r_1,\cdots,r_n)$$\in\R^n$ and denote $P(y_1+r_1,\cdots,y_n+r_n)$
by $\tilde{P}$. Let $\frac{\partial}{\partial y_i}$ denote the
standard partial derivation on the ring $\R[y_1,\cdots,y_n]$. From
the Taylor series expansion of $\tilde{P}$, we have

\begin{equation}\label{taylor for multi}
\tilde{P}=P+\sum^n_{i=1}r_i\frac{\partial P}{\partial y_i}
+\frac{1}{2!}\sum^n_{j=1}\sum^n_{i=1}r_ir_j\frac{\partial^2
P}{\partial y_j\partial y_i}+\cdots
\end{equation}

\bprop\label{trans poly prop}

Let $P\in\R[y_1,\cdots,y_n]$ and for $1\leq i\leq n$ let
$r_i\in\R$. Suppose that $P$ divides
$\tilde{P}:=P(y_1+r_1,\cdots,y_n+r_n)$. Then $P=\tilde{P}$ and
$\sum^n_{i=1}r_i\frac{\partial H_j}{\partial y_i}=0$ for every
homogeneous component $H_j$ of total degree $j$ of $P$. In
particular $H_j=\tilde{H}_j$ for every $j$ and
$\sum^n_{i=1}r_i\frac{\partial P}{\partial y_i}=0$. \eprop

\bpf Rewrite the equation \ref{taylor for multi} as

\begin{equation}\label{taylor for multi rewr}\tilde{P}-P=\sum^n_{i=1}r_i\frac{\partial P}{\partial y_i}
+\frac{1}{2!}\sum^n_{j=1}\sum^n_{i=1}r_ir_j\frac{\partial^2
P}{\partial y_j\partial y_i}+\cdots\end{equation}If $P\in\R$ then
the proposition follows immediately. Assume that $P$ has a
monomial whose total degree is $\geq 1$. We observe that the
operator $\sum^n_{i=1}r_i\frac{\partial }{\partial y_i}$ is
applied to a monomial of $P$ reduces the total degree of that
monomial by one and the operator
$\sum^n_{j=1}\sum^n_{i=1}r_ir_j\frac{\partial^2 }{\partial
y_j\partial y_i}$ applied to a monomial reduces its total degree
by two and so on... Thus the total degree of the RHS of equation
\ref{taylor for multi rewr} is less than the total degree of $P$.
Clearly, $P$ divides $\tilde{P}$ implies $P$ divides the LHS of
equation \ref{taylor for multi rewr} and therefore $P$ divides the
RHS whose total degree is less than that of $P$. Thus RHS of
\ref{taylor for multi rewr} equals 0, that is
\begin{equation}\label{taylor tail multi}\sum^n_{i=1}r_i\frac{\partial P}{\partial y_i}
+\frac{1}{2!}\sum^n_{j=1}\sum^n_{i=1}r_ir_j\frac{\partial^2
P}{\partial y_j\partial y_i}+\cdots=0\end{equation}and hence
$\tilde{P}=P$.

Let $P=\sum^k_{j=0}H_j$, where $H_j$ is the homogenous component
of total degree $j$ of $P$. Again we observe that when the
operator $\sum^n_{i=1}r_i\frac{\partial }{\partial y_i}$ is
applied to a monomial of $H_j$, the degree of that monomial goes
down by one. Therefore, if this operator is applied to a
homogenous component $H_j$, either $\sum^n_{i=1}r_i\frac{\partial
H_j}{\partial y_i}=0$, as cancellation of monomials may occur, or
the total degree of $\sum^n_{i=1}r_i\frac{\partial H_j}{\partial
y_i}$ has to be one lesser than that of $H_j$.

Now consider the homogeneous component $H_k$. We know that the
$\sum^n_{i=1}r_i\frac{\partial H_k}{\partial y_i}=0$ or the total
degree of $\sum^n_{i=1}r_i\frac{\partial H_k}{\partial y_i}$ is
$k-1$. The latter cannot happen since from  equation \ref{taylor
tail multi} we have
\begin{equation}\label{deriv of taylor}-\sum^n_{i=1}r_i\frac{\partial H_k}{\partial
y_i}=\sum^{k-1}_{l=0}\sum^n_{i=1}r_i\frac{\partial H_l}{\partial
y_i} +\frac{1}{2!}\sum^n_{j=1}\sum^n_{i=1}r_ir_j\frac{\partial^2
P}{\partial y_j\partial y_i}+\cdots\end{equation} and the RHS of
the above equation is of total degree $\leq k-2$. Thus
$\sum^n_{i=1}r_i\frac{\partial H_k}{\partial y_i}=0$. Note that
$\sum^n_{i=1}r_i\frac{\partial H_k}{\partial y_i}=0$ implies
$\sum^n_{j=1}\sum^n_{i=1}r_ir_j\frac{\partial^2 H_k}{\partial
y_j\partial y_i}=0$ and so on... Therefore from equation
\ref{taylor for multi} we get $H_k=\tilde{H}_k$.

Now substituting
$$\sum^n_{j=1}\sum^n_{i=1}r_ir_j\frac{\partial^2 P}{\partial
y_j\partial
y_i}=\sum^{k-1}_{l=0}\sum^n_{j=1}\sum^n_{i=1}r_ir_j\frac{\partial^2
H_l}{\partial y_j\partial y_i}$$

and $\sum^n_{i=1}r_i\frac{\partial H_k}{\partial y_i}=0$ in
equation \ref{deriv of taylor}, we get
$$-\sum^n_{i=1}r_i\frac{\partial H_{k-1}}{\partial
y_i}=\sum^{k-2}_{l=0}\sum^n_{i=1}r_i\frac{\partial H_l}{\partial
y_i}+\frac{1}{2!}\sum^{k-1}_{l=0}\sum^n_{j=1}\sum^n_{i=1}r_ir_j\frac{\partial^2
H_l}{\partial y_j\partial y_i}+\cdots.$$ By comparing the total
degrees of the LHS and RHS, we conclude that
$\sum^n_{i=1}r_i\frac{\partial H_{k-1}}{\partial y_i}=0$ and thus
$H_{k-1}=\tilde{H}_{k-1}$. Similarly we can show that
$\sum^n_{i=1}r_i\frac{\partial H_j}{\partial y_i}=0$ for every
$j$. From this equation it is easy to see that $H_i=\tilde{H}_i$
and $\sum^n_{i=1}r_i\frac{\partial P}{\partial y_i}=0$. \epf

\bprop \label{variety desc of a homo poly}

For every homogeneous polynomial $P\in\R[y_1,\cdots,y_n]$ there is
a system $\{D_j\}$ of linear forms over $\R$ such that
$P=\tilde{P}$ for some $(r_1,\cdots,r_n)\in\R^n$ if and only if
$(r_1,\cdots,r_n)\in\R^n$ is a solution of the system $\{D_j\}$.

\eprop

\bpf

Suppose that $P=\tilde{P}$ for some $(r_1,\cdots,r_n)\in\R^n$ then
from proposition \ref{trans poly prop} we see that
\begin{equation}\label{homo poly deriv eqn}\sum^n_{i=1}r_i\frac{\partial P}{\partial
y_i}=0.\end{equation}

By grouping all the monomials, we could rewrite
$\sum^n_{i=1}r_i\frac{\partial P}{\partial y_i}$ as
$$\sum^t_{j=1}D_j(r_1,\cdots,r_n)X_{\omega_j},$$ where $\{D_j\}$
is a system of linear forms over $\R$ and $X^{\omega_j}$
represents a primitive monomial that appears in
$\sum^n_{i=1}r_i\frac{\partial P}{\partial y_i}$. Thus equation
\ref{homo poly deriv eqn} becomes
$$\sum^t_{j=1}D_j(r_1,\cdots,r_n)X_{\omega_j}=0$$ and clearly $P$
satisfies equation \ref{homo poly deriv eqn} if and only if the
the tuple $(r_1,\cdots,$ $r_n)\in\R^n$ satisfies the system
$\{D_j\}$. \epf

\bprop\label{reduction of poly to c poly}

Let $\R:=\C[x_1,\cdots,x_m]$ be a polynomial ring and let
$D(y_1,\cdots,$ $y_n)$ be a linear form over the ring $\R$ with
variables $y_1,\cdots,y_n$. Then there is a system $\{L_j\}$ of
linear forms over $\C$ such that $D(c_1,\cdots,c_n)=0$ for
$(c_1\cdots,c_n)\in\C^n$ if and only if $(c_1\cdots,c_n)\in\C^n$
is a solution of the system $\{L_j\}$ \eprop

\bpf By viewing the polynomial
$D(y_1,\cdots,y_n)\in\R[y_1,\cdots,y_n]$ as a polynomial over the
ring $\C[y_1,\cdots,y_n]$ with variables $x_1,\cdots,x_m$, we
obtain vectors $\omega_j:=(\omega_{j1},\cdots\omega_{jm})$
$\in\mathbb{W}^m$, where $\mathbb{W}:=\N\cup\{0\}$ and linear
forms $L_j(y_1,\cdots,y_n)\in span_{\C}\{x_1,$ $\cdots,x_m\}$ such
that
$$D(y_1,\cdots,y_n)=\sum^t_{j=1}L_j(y_1,\cdots,y_n)X_{\omega_j},$$
where $X_{\omega_j}$ is the primitive monomial
$x^{\omega_{j1}}_1\cdots x^{\omega_{jm}}_m$. Since primitive
monomials are linearly independent over constants, we see that
$D(c_1,$ $\cdots,c_n)=0$ if and only if $(c_1,\cdots,c_n)$ is a
solution of the system $\{L_j|1\leq j\leq t\}$ of linear forms
over $\C$. \epf

\bpf[Proof of theorem \ref{struct theorem for antideriv}.] Let
$\G:=\G(\pv|\f)$ and let $\sg\leq\G$ be the group of all
automorphisms that fixes $\f\langle u\rangle$. Then for any
$\sigma=(c_{1\sigma},\cdots,c_{n\sigma})\in\sg$ we have
$\sigma(u)=u$, that is
$$\frac{\sigma(P)}{\sigma(Q)}=\frac{P}{Q}$$ and thus \begin{equation}
\label{eqn stru of anti}\sigma(P)Q=\sigma(Q)P.\end{equation}

Since $(P,Q)=1$, from equation \ref{eqn stru of anti} we see that
$P$ divides $\sigma(P)$ and $Q$ divides $\sigma(Q)$. Note that
$\sigma(P)=P(\x_1+c_{1\sigma},\cdots,\x_n+c_{n\sigma})$ and
$\sigma(Q)=Q(\x_1+c_{1\sigma},\cdots,\x_n+c_{n\sigma})$ and
therefore from proposition \ref{trans poly prop} we obtain
\begin{equation}\label{fixing P and Q}\sigma(P)=P\qquad \text{and}\qquad\sigma(Q)=Q.\end{equation}
If both $P,Q\in\f$ then $\G$ fixes $u$ and thus $\f\langle
u\rangle=\f$. Let us assume $P\notin\f$ and denote
$P(\x_1+c_{1\sigma},\cdots,\x_n+c_{n\sigma})$ by $\tilde{P}$. Now
apply propositions \ref{trans poly prop} and \ref{variety desc of
a homo poly}  with $R:=\f$  to get linear forms $\{A_i|1\leq i\leq
s\}$ $\subset Span_\f\{\x_1\cdots,\x_n\}$ such that
$A_i(c_{1\sigma},\cdots,c_{n\sigma})=0$ iff $\sigma(P)=P$. We also
see that $A_i$ is fixed by all $\sigma\in\sg$. Therefore from the
fundamental theorem we conclude that $P\in\f(A_i|1\leq i\leq
s)\subset\f\langle u\rangle$. Similarly if $Q\notin\f$ then one
can find these linear forms for $Q$ say $\{B_i|1\leq i\leq t\}$
$\subset Span_\f\{\x_1\cdots,\x_n\}$ such that $Q\in\f(B_i|1\leq
i\leq t)$$\subseteq\f\langle u\rangle$. Now
$u=\frac{P}{Q}\in\f(D_i|1\leq i\leq r)$, where $\{D_i|1\leq i\leq
r\}$ $=\{A_i|1\leq i\leq s\}$ $\cup\{B_i|1\leq i\leq t\}$. On the
other hand, both the fields $\f(A_i|1\leq i\leq s)$ and
$\f(B_i|1\leq i\leq t)$ are subfields of $\f\langle u\rangle$.
Thus we see that
$$\f\langle u\rangle=\f(D_i|1\leq i\leq t).$$\epf

\begin{remark}({\bf Algorithm})

Let $\f(\x_1,\cdots,\x_n)$ be an extension by antiderivatives
$\x_1,\cdots,\x_n$ of $\f$ and assume that $\x_1,\cdots,\x_n$ are
algebraically independent over $\f$. Let $u=\frac{P}{Q}$,
$P,Q\in\f[\x_1,\cdots,\x_n]$ and $(P,Q)=1$. To compute the
differential field $\f\langle u\rangle$ we do the following:

1. Observe from equation \ref{fixing P and Q} that $\sigma(u)=u$
if and only if $\sigma(P)=P$ and $\sigma(Q)=Q$.

2. Find all tuples $(c_1,\cdots,c_n)\in\C^n$ such that
$P=P(\x+c_1,\cdots,\x_n+c_n)$ and $Q=Q(\x+c_1,\cdots,\x_n+c_n)$.
Steps 2a, 2b and 2c computes the same.

2a. From proposition \ref{trans poly prop}, we see that
$P=P(\x+c_1,\cdots,\x_n+c_n)$ if and only if
$$\sum^n_{i=1}c_i\frac{\partial P}{\partial y_i}=0$$ and similarly $Q=Q(\x+c_1,\cdots,\x_n+c_n)$
if and only if $$\sum^n_{i=1}c_i\frac{\partial Q}{\partial
y_i}=0.$$

2b. We rewrite the above equations as
$$\sum^t_{j=1}A_j(c_1,\cdots,c_n)X_{\omega_j}=0$$ and
$$\sum^s_{j=1}B_j(c_1,\cdots,c_n)Y_{\omega_j}=0,$$  where $\{A_j|1\leq j\leq t\}\subset$
$Span_\f\{\x_1,\cdots,\x_n\}$ is a system of linear forms over
$\f$ and $X^{\omega_j}$ represents a primitive monomial that
appears in $\sum^n_{i=1}r_i\frac{\partial P}{\partial y_i}$ and
$\{B_j|1\leq j\leq s\}\subset$ $Span_\f\{\x_1,\cdots,\x_n\}$ is a
system of linear forms over $\f$ and $Y^{\omega_j}$ represents a
primitive monomial that appears in $\sum^n_{i=1}c_i\frac{\partial
Q}{\partial y_i}$

2c. Observe that the displayed equations from 2b holds if and only
if $A_j(c_1,$ $\cdots,c_n)=0$ for all $1\leq j\leq t$ and
$B_j(c_1,\cdots,$ $c_n)=0$ for all $1\leq j\leq s$. Thus
$\sigma(u)=u$ if and only if
$\sigma:=(c_{1\sigma},\cdots,c_{n\sigma})$ is a solution of the
system $\{D_i|1\leq i\leq r\}$ $=\{A_i|1\leq i\leq s\}$
$\cup\{B_i|1\leq i\leq t\}$.

3. Thus the algebraic subgroup of all automorphisms of $\G$ that
fixes $u$ also fixes $\{D_j|1\leq j\leq r\}$ and vice versa.
Therefore from the fundamental theorem we conclude that $\f\langle
u\rangle$ equals the differential field $\f(D_j|1\leq j\leq r)$.

4. Finally, if $\f$ is a fraction field of a polynomial ring
$\R:=\C[x_1,\cdots,x_s]\subset\f$ then from proposition
\ref{reduction of poly to c poly} we see that each of the $D_j$'s
can be reduced to a finite set of linear forms $L_{ji}$, $1\leq
i\leq m_j$ over $\C$ and thus $\f\langle u\rangle=$ $\f(D_j|1\leq
j\leq r)=$ $\f(L_i|1\leq i\leq m)$, where $\{L_i|1\leq i\leq m\}$
$=\cup^r_{j=1}\{L_{ji}|1\leq i\leq m_j\}$.
\end{remark}

\bprop \label{Q- relatively prime}Let
$\f(\x_1,\cdots,\x_l)\supset\f$ be an extension by antiderivatives
$\x_1,\cdots,\x_l$ of $\f$ and suppose that $\x_1,\cdots,\x_l$ are
algebraically independent over $\f$. If $R\in\f[\x_1,\cdots,\x_l]$
is an irreducible polynomial then the polynomials $R$ and $R'$ are
relatively prime. \eprop

\bpf

Let $R\in\f[\x_1,\cdots,\x_l]$ be an irreducible polynomial.
Suppose that $R'$ and $R$ are not relatively prime. Then $R$,
being irreducible, has to divide $R'$. Observe that the total
degree of $R'$ is $\leq$ the total degree of $R$ and since $R$
divides $R'$, the total degree of $R$ equals the total degree of
$R'$. Thus
$$R'=\rf R$$ for some $\rf\in\f$. Let $\G$ be the differential
Galois group of $\f(\x_1,\cdots,\x_l)$ over $\f$  and let
$\sigma\in\G$. We observe that $\sigma(\x_i)=\s_i+c_{i\sigma}$,
$c_{i\sigma}\in\C$ and therefore
$\sigma(R)=R(\x_1+c_{1\sigma},\cdots,\x_l+c_{l\sigma})$. We also
observe that $R'=\rf R$ implies $\sigma(R)=c_\sigma R$ for some
$c_\sigma\in\C^\times$. Then $R$ divides $\sigma(R)$ and thus from
proposition \ref{trans poly prop} we obtain $\sigma(R)=R$. Thus
every automorphism of $G$ has to fix $R$ and since
$\f(\x_1,\cdots,\x_l)$ is a Picard-Vessiot extension of $\f$, we
obtain $R\in\f$, a contradiction. \epf

\bt \label{anti deriv in field cri}Let
$\f(\x_1,\cdots,\x_l)\supset\f$ be an extension by antiderivatives
$\x_1,\cdots,\x_l$ of $\f$. Let $S,T\in\f[\x_1,\cdots,\x_n]$ be
relatively prime polynomials and assume that $T$ has an
irreducible factor $R\in\f[x_1,\cdots,\x_l]$ such that $R^2$ does
not divide $T$. Then there is no $\y\in\f(\x_1,\cdots,\x_l)$ such
that $\y'=\frac{S}{T}$.

\et

\bpf

Suppose that there is a $\y\in\f(\x_1,\cdots,\x_l)$ such that
$\y'=\frac{S}{T}$. There are relatively prime polynomials
$P,Q\in\f[\x_1,\cdots,\x_l]$ such that $\y=\frac{P}{Q}$. Thus
taking the derivative we arrive at
\begin{equation}Q^2 S=T(P'Q-Q'P).\end{equation}

Note that $R$ is an irreducible factor of $T$ and therefore from
the above equation $R$ divides $Q^2S$. Since $S$ and $T$ are
relatively prime, $R$ has to divides $Q^2$, which implies $R$
divides $Q$. Let $n$ be the largest integer so that $R^n$ divides
$Q$. Then $R^{n+1}$ divides $Q^2$ and again from the above
displayed equation, $R^{n+1}$ divides $T(P'Q-Q'P)$. Note that $R$
divides $T$ but $R^2$ does not and thus $R^n$ divides $P'Q-Q'P$.
Since $R^n$ divides $Q$, and $P$ and $Q$ are relatively prime, we
obtain $R^n$ divides $Q'$. Let $H\in\f[\x_1\cdots,\x_l]$ be a
polynomial such that $Q=R^nH$. Note that $R$ and $H$ are
relatively prime polynomials.  Then $R^n$ divides
$Q'=nR^{n-1}R'H+R^nH'$ implies $R$ divides $R'$, which
 contradicts proposition \ref{Q- relatively prime}. \epf

\subsection{Extensions by exponentials of integrals}

Let $\f$ be a differential field with an algebraically closed
field of constants $\C$. Let $\pv\supset\f$ be an extension by
exponentials of integrals $\ex_1,\cdots,\ex_n$ of $\f$ and $\G$
the group of all differential automorphisms of $\pv$ over $\f$.
Since $f_i:=\frac{\ex'_i}{\ex_i}\in\f$, $\ex_i$ satisfies the
first order linear homogeneous differential equation
$\ex'_i=f_i\ex_i$. For any $\sigma\in\G$,
$\sigma(\ex_i)'=f_i\sigma(\ex_i)$ and thus
$(\frac{\sigma(\ex_i)}{\ex_i})'=0$. Since $\pv$ is a NNC extension
of $\f$, there is a $c_{i\sigma}\in\C\setminus\{0\}$ such that
$\frac{\sigma(\ex_i)}{\ex_i}=c_{i\sigma}$. Thus
$\sigma(\ex_i)=c_{i\sigma}\ex_i$. Also note that the action of
$\sigma$ on the elements $\ex_i$ completely determines the
automorphism $\sigma$. For any $\phi,\sigma\in\G$,
\begin{equation}\phi(\sigma(\ex_i))=\phi(c_{i\sigma}\ex_i)=c_{i\phi}c_{i\sigma}\ex_i
=c_{i\sigma}c_{i\phi}\ex_i=\sigma(\phi(\ex_i)).\end{equation} Thus
$\G$ is a commutative group and also the map
$\sigma\mapsto(c_{1\sigma},\cdots,c_{n\sigma})$ is an injective
algebraic group homomorphism from $\G$ to
$(\C\setminus\{0\},\times)^n$.

If $\pv=\f(\ex)$, $\frac{\ex'}{\ex}\in\f$ then $\G$ is an
algebraic subgroup of $(\C\setminus\{0\},\times)$. Thus if $\G$ is
non trivial then it has to be a finite subgroup of
$(\C\setminus\{0\}, \times)$. Note that $\G$ could be a finite
subgroup of $(\C\setminus\{0\}, \times)$; for example, let
$\f=\mathbb{C}(x)$ and let $\pv=\f(\sqrt[n]{x})$, $n\geq 2$. Then
we have the equation
$$(\sqrt[n]{x})'=\frac{1}{nx}\sqrt[n]{x}.$$ Thus $\pv$ is an
extension by an exponential of an integral $\sqrt[n]{x}$ of $\f$.
Clearly $\sqrt[n]{x}\notin\f$(therefore $\G$ is not the trivial
group) and for any automorphism $\sigma\in\G$
\begin{align*}&\sigma(\sqrt[n]{x})=c_{\sigma}\sqrt[n]{x}\\
\iff&(\sigma(\sqrt[n]{x}))^n=c_{\sigma}^n(\sqrt[n]{x})^n\\
\iff&\sigma(x)=c_{\sigma}^nx\\
\iff&1=c_{\sigma}^n
\end{align*}

In fact one can also show that $\G$ is the group of nth roots of
unity (follows from the fact that the ordinary Galois group and
the differential Galois group are the same if the extension $\pv$
of $\f$ is finite).

Let $\Pp:=\{\prod^k_{i=1}\ex^{m_i}_i|m_i\in\Z^*\}$, the set of all
power products of $\{\ex_i|1\leq i\leq n\}$.  We will now prove
the following theorem

\bt \label{stru theorem for expos} Let
$\pv=\f(\ex_1,\cdots,\ex_n)$ be an extension of $\f$ by
exponentials of integrals $\ex_1,\cdots,\ex_n$ of $\f$ and let
$\ex_1,\cdots,\ex_n$ are algebraically independent over $\f$. Let
$u=\frac{P}{Q}$, $P,Q\in\f[\ex_1,\cdots,\ex_n]$ and $(P,Q)=1$.
Then there are power products $\mpp_j\in\Pp$, $1\leq j\leq t$ such
that
$$\f\langle u\rangle=\f(\mpp_1,\cdots,\mpp_t).$$ Moreover, we may
explicitly compute the power products $\mpp_i$ from $P$ and $Q$.
\et

\bpf Let $\G:=\G(\pv|\f)$ and let $\sg\leq\G$ be the group of all
automorphisms of $\G$ that fixes $u$. So, for $\sigma\in\sg$ we
have $\sigma(u)=u$ and therefore $\sigma(P)Q=\sigma(Q)P$. Thus $P$
divides $\sigma(P)Q$ and since $(P,Q)=1$, $P$ divides $\sigma(P)$
and similarly $Q$ divides $\sigma(Q)$.

We may assume either $P$ or $Q$ is not in $\f$; otherwise the
differential field $\f\langle u\rangle=\f$. Assume that
$P\notin\f$ and write
\begin{equation}\label{power prod rep}P=\sum^r_{i=1}f_{\pp_i}
\pp_i,\end{equation} where $\pp_i$ are primitive monomials and
$f_{\pp_i}\in\f$. Note that
$$\sigma(P)=\sum^r_{i=1}f_{\pp_i}
\pp_i(c_{1\sigma},\cdots,c_{n\sigma})\pp_i$$ and since $c_{i\sigma
}\in\C\setminus\{0\}$, $\pp_i(c_{1\sigma},\cdots,c_{n\sigma})\neq
0$. Thus $P$ and $\sigma(P)$ have the same number terms and every
monomial that appears in $P$ also appears in $\sigma(P)$ and vice
versa. But $P$ divides $\sigma(P)$ and therefore there is a
$d_\sigma\in\f$ such that $\sigma(P)=d_\sigma P$. In fact
$\pp_i(c_{\sigma 1},\cdots, c_{\sigma n})=d_\sigma$ for all $i$
and thus $d_\sigma\in\C\setminus\{0\}$.

This shows that $\frac{\pp_i}{\pp_1}$ is fixed by every
$\sigma\in\sg$. Thus, from fundamental theorem, we obtain
$\C\langle u\rangle\supset \f(\frac{\pp_i}{\pp_1}|1\leq i\leq r)$.
Since $Q$ also divides $\sigma(Q)$, writing
$Q=\sum^s_{j=1}g_{\pq_j} \pq_j$ similar to equation \ref{power
prod rep}, we conclude that there is a
$e_\sigma\in\C\setminus\{0\}$ such that $\sigma(Q)=e_\sigma Q$.
Since $\sigma(\frac{P}{Q})=\frac{P}{Q}$, we have
$d_\sigma=e_\sigma$ and thus $\sigma$ fixes $\frac{\pq_j}{\pp_1}$.
Thus $\C\langle u\rangle\supset \f(\frac{\pq_j}{\pp_1}|1\leq j\leq
t)$.
 Now we have
$$\C\langle u\rangle\supset
\f\Big(\frac{\pp_i}{\pp_1},\frac{\pq_j}{\pp_1}|1\leq i\leq r,
1\leq j\leq s\Big).$$ On the other hand we could write
$$u=\frac{P}{Q}=\frac{\sum^r_{i=1}f_{\pp_i}
\frac{\pp_i}{\pp_1}}{\sum^r_{j=1}g_{\pq_j} \frac{\pq_j}{\pp_1}}$$
for all $1\leq i,j\leq r$. Hence from fundamental theorem it
follows that
$$\f\langle u\rangle=\f(\mpp_1,\cdots,\mpp_t),$$ where
$\{\mpp_1,\cdots,\mpp_t\}$
$=\{\frac{\pp_i}{\pp_1},\frac{\pq_j}{\pp_1}|1\leq i\leq r, 1\leq
j\leq s\}$. \epf


Now we will prove a theorem which is a combination of theorems
\ref{struct theorem for antideriv} and \ref{stru theorem for
expos}.

\bt\label{main struct thereom} Let $\pv\subset\f$ be a NNC
extension and let $\pv=\f(\x_1,\cdots,\x_n,$ $\ex_1\cdots,\ex_m)$,
where $\x'_i\in\f$, $\frac{\ex'_i}{\ex_i}\in\f$ and
$\x_1,\cdots,\x_n,\ex_1\cdots,\ex_m$ are algebraically independent
over $\f$. Let $u\in\pv$ and suppose that $u=\frac{P}{Q}$, where
$P,Q\in\f[\x_1,\cdots,\x_n,\ex_1\cdots,\ex_m]$ and $(P,Q)=1$. Then
for $i=1,2,\cdots,t$ and $j=1,2,\cdots s$ there are $\f-$ linear
forms $\Dl_i$ over the set $\{\x_i|1\leq i\leq n\}$ and power
products $\mpp_j$ over the set $\{\ex_i|1\leq i\leq m\}$ such that
$$\f\langle u\rangle=\f(\Dl_i,\mpp_j|1\leq i\leq t, 1\leq j\leq s
).$$ Moreover these forms can be explicitly computed from the
polynomials $P$ and $Q$. \et

\bpf

Let $u\neq 0$ and
$u=\frac{P}{Q}\in\f[\x_1,\cdots,\x_n,\ex_1\cdots,\ex_m]$,
$(P,Q)=1$. Rewrite $P$ and $Q$ as polynomials over the ring
$\f[\x_1,\cdots,\x_n][\ex_1\cdots,\ex_m]$. That is
$P=\sum^k_{i=0}a_{\pp_i}\pp_i$, $Q=\sum^l_{i=0}b_{\pq_i}\pq_i$,
where $a_{\pp_i},b_{\pq_i}\in\f[\x_1,\cdots,\x_n]$ and $a_{\pp_k}$
and $b_{\pq_l}$ are non zero. Now divide through $P$ and $Q$ by
$a_{\pp_k}$. Thus we obtain
\begin{equation}\label{eq main stru thm}
u=\frac{P}{Q}=\frac{\sum^k_{i=0}\frac{a_{\pp_i}}{a_{\pp_k}}\pp_i}
{\sum^l_{i=0}\frac{b_{\pq_i}}{a_{\pp_k}}\pq_i}
\end{equation}
and now the polynomials $P,Q$ becomes polynomials over the ring
$\mathbf{K}[\ex_1,$ $\cdots,\ex_m]$, where
$\mathbf{K}:=\f(\x_1,\cdots,\x_n)$. Hereafter we will call
$\sum^k_{i=0}\frac{a_{\pp_i}}{a_{\pp_k}}\pp_i$ as $P$ and
$\sum^l_{i=0}\frac{b_{\pq_i}}{a_{\pp_k}}\pq_i$ as $Q$. Note that
$P$ and $Q$ are relatively prime in the ring
$\mathbf{K}[\ex_1,\cdots,\ex_m]$.

We observe that $\pv\supseteq\f$ is a P-V extension and let $\G$
be the group of differential automorphisms of $\pv\supseteq\f$.
Thus there is a subgroup $\sg\leq\G$ such that $\f\langle
u\rangle$ is the fixed field of $\sg$. Let $\sigma\in\sg$. Then
$\sigma(u)=u$ and therefore we obtain \begin{equation}\label{usual
eqn} \sigma(P)Q=\sigma(Q)P.\end{equation} Since $(P,Q)=1$ in
$\mathbf{K}[\ex_1,\cdots,\ex_m]$, $P$ divides $\sigma(P)$ and $Q$
divides $\sigma(Q)$.

We observe that
$$\sigma(P)=\sum^k_{i=0}\frac{\sigma(a_{\pp_i})}
{\sigma(a_{\pp_k})}\pp_i(c_{1\sigma},\cdots,c_{m\sigma})\pp_i,$$
$$\sigma(Q)=\sum^l_{i=0}\frac{\sigma(b_{\pq_i})}
{\sigma(a_{\pp_k})}\pq_i(c_{1\sigma},\cdots,c_{m\sigma})\pq_i$$
and that $\frac{\sigma(a_{\pp_i})}
{\sigma(a_{\pp_k})}\in\mathbf{K}$ since $\mathbf{K}$ is a normal
extension of $\f$, and
$\pp_i(c_{1\sigma},\cdots,c_{m\sigma})\in\C$. Therefore
$\sigma(P),\sigma(Q)\in\mathbf{K}[\ex_1,\cdots,\ex_n]$.

Claim:
$\sigma(\frac{a_{\pp_i}}{a_{\pp_k}})=\frac{a_{\pp_i}}{a_{\pp_k}}$,
$\sigma(\frac{b_{\pq_i}}{a_{\pp_k}})=\frac{b_{\pq_i}}{a_{\pp_k}}$,
$\sigma(\frac{\pp_i}{\pp_k}) =\frac{\pp_i}{\pp_k}$ and
$\sigma(\frac{\pq_i}{\pp_k}) =\frac{\pq_i}{\pp_k}$.

From the facts that $P$ divides $\sigma(P)$,
\begin{equation*}P=\sum^k_{i=0}\frac{a_{\pp_i}}{a_{\pp_k}}\pp_i\ \
\text{and}\ \ \sigma(P)=\sum^k_{i=0}\frac{\sigma(a_{\pp_i})}
{\sigma(a_{\pp_k})}\pp_i(c_{1\sigma},\cdots,c_{m\sigma})\pp_i,\end{equation*}
we see  $\sigma(P)=\pp_k(c_{1\sigma},\cdots,c_{m\sigma})P$. Since
$\pp_i$ are linearly independent over $\mathbf{K}$, for each $i$,
we have
$$\frac{\sigma(a_{\pp_i})}{\sigma(a_{\pp_k})}=\frac{\pp_k(c_{1\sigma},
\cdots,c_{m\sigma})}{\pp_i(c_{1\sigma},\cdots,c_{m\sigma})}\frac{a_{\pp_i}}{a_{\pp_k}}.$$
Observe that $a_{\pp_i}\in\f[\x_1\cdots,\x_n]$ and now replace
$\frac{a_{\pp_i}}{a_{\pp_k}}$ by
$\frac{\alpha_{\pp_i}}{\beta_{\pp_i}}$, where
$\alpha_{\pp_i}:=\frac{a_{\pp_i}}{g_i}$,
$g_i:=(a_{\pp_i},a_{\pp_k})$ and
$\beta_{\pp_i}:=\frac{a_{\pp_k}}{g_i}$.  Thus we have
\begin{equation}\label{anot usual eqn}
\sigma(\alpha_{\pp_i})\beta_{\pp_i}=\frac{\pp_k(c_{1\sigma},
\cdots,c_{m\sigma})}{\pp_i(c_{1\sigma},\cdots,c_{m\sigma})}\alpha_{\pp_i}
\sigma(\beta_{\pp_i}).
\end{equation}
Clearly $(\alpha_{\pp_i},\beta_{\pp_i})=1$ and since
$\frac{\pp_k(c_{1\sigma},
\cdots,c_{m\sigma})}{\pp_i(c_{1\sigma},\cdots,c_{m\sigma})}\in\C$,
we have $\alpha_{\pp_i}$ divides $\sigma(\alpha_{\pp_i})$ and
$\beta_{\pp_i}$ divides $\sigma(\beta_{\pp_i})$. Apply proposition
\ref{trans poly prop} and obtain
$\sigma(\alpha_{\pp_i})=\alpha_{\pp_i}$,
$\sigma(\beta_{\pp_i})=\beta_{\pp_i}$ and thus from equation
\ref{anot usual eqn} we have
\begin{equation*}\frac{\pp_i(c_{1\sigma},
\cdots,c_{m\sigma})}{\pp_k(c_{1\sigma},\cdots,c_{m\sigma})}=1.
\end{equation*}
From this equation it is clear that
\begin{equation}\label{cruc expo eqn}\sigma(\frac{\pp_i}{\pp_k})
=\frac{\pp_i}{\pp_k}.\end{equation} Since
$\sigma(\frac{\alpha_{\pp_i}}{\beta_{\pp_i}})=\frac{\alpha_{\pp_i}}{\beta_{\pp_i}}$,
we have
$\sigma(\frac{a_{\pp_i}}{a_{\pp_k}})=\frac{a_{\pp_i}}{a_{\pp_k}}$
for each $i$. The claims $\sigma(\frac{b_{\pq_i}}{a_{\pp_k}})
=\frac{b_{\pq_i}}{a_{\pp_k}}$ and $\sigma(\frac{\pq_i}{\pp_k})
=\frac{\pq_i}{\pp_k}$ follows similarly.

We may apply theorem \ref{struct theorem for antideriv} for each
$\alpha_{\pp_i}$ and $\beta_{\pp_i}$ and obtain $\f-$ linear forms
over $\{\x_1,\cdots,\x_n\}$ so that the differential fields
$\f\langle \alpha_{\pp_i}\rangle$ and $\f\langle
\beta_{\pp_i}\rangle$ equals the field generated by their
corresponding linear forms. Thus we have linear forms
$\{D_{i1},\cdots,D_{it_i}\}$ such that
$$\f\langle\frac{\alpha_{\pp_i}}{\beta_{\pp_i}}\rangle=
\f\langle\alpha_{\pp_i},\beta_{\pp_i}\rangle=\f(D_{i1},\cdots,D_{it_i}).$$

Note that
$\frac{\alpha_{\pp_i}}{\beta_{\pp_i}}=\frac{a_{\pp_i}}{a_{\pp_k}}$
and therefore $$\f\langle \frac{a_{\pp_i}}{
a_{\pp_k}}\rangle=\f(D_{i1},\cdots,D_{it_i}).$$

Similarly we can obtain linear forms $\{E_{j1},\cdots,E_{js_j}\}$
so that
$$\f\langle \frac{b_{\pq_j}}{
a_{\pp_k}}\rangle=\f(E_{j1},\cdots,E_{js_j}).$$

Let $\{\Dl_i|1\leq i\leq t\}$
$=\{D_{i1},\cdots,D_{it_i}\}\cup\{E_{j1},\cdots,E_{js_j}\}$,
$\{\mpp_1,\cdots,\mpp_s\}=\{\pp_i|1\leq i\leq
k\}\cup$$\{\pq_j|1\leq j\leq l\}$ and $\mpp_1:=\pp_k$. Then
writing
$$u=\frac{\sum^k_{i=0}\frac{a_{\pp_i}}{a_{\pp_k}}\frac{\pp_i}{\pp_k}}
{\sum^l_{i=0}\frac{b_{\pq_i}}{a_{\pp_k}}\frac{\pq_i}{\pp_k}},$$ we
immediately see that
$$\sigma(u)=u\Leftrightarrow\sigma(\Dl_i)=\Dl_i, \sigma(\frac{\mpp_i}{\mpp_1})
=\frac{\mpp_i}{\mpp_1}.$$ Hence the theorem. \epf

\subsection{Tower of Extensions by Antiderivatives}
Let $\f$ be a differential field with an algebraically closed
field of Constants $\C$ and let $\f_\infty$ be a complete
Picard-Vessiot closure of $\f$ (every homogeneous linear
differential equation over $\f_\infty$ has a full set of solutions
in $\f_\infty$ and it has $\C$ as its field of constants and
$\f_\infty$ is minimal with respect to these properties). All the
differential fields under consideration are subfields of
$\f_\infty$.

A differential field extension $\pv$ of $\f$ is called a
\textsl{tower of extension by antiderivatives} if there are
differential fields $\pv_i$, $0\leq i\leq n$ such that
$$\pv:=\pv_n\supseteq\pv_{n-1}\supseteq\cdots,\supseteq\pv_1\supseteq\pv_0:=\f$$
and $\pv_i$ is an extension by antiderivatives of $\pv_{i-1}$ for
each $1\leq i\leq n$.

\bt\label{no alg extn thm}

Let $\M\supseteq\f$ be differential fields and let
$$\pv:=\pv_n\supset\pv_{n-1}\supset\cdots\supset\pv_1\supset\pv_0:=\f$$
be a tower of extensions by antiderivatives. Then $u\in\pv$ is
algebraic over $\M$ only if $u\in\M$.

\et

\bpf

We will use an induction on $n$ to prove this theorem. Consider
the tower
$$\M\cdot\pv:=\M\cdot\pv_n\supseteq\M\cdot\pv_{n-1}\supseteq\cdots\supseteq\M\cdot\pv_1\supseteq\M.$$
Clearly, the above tower is a tower of extension by
antiderivatives. Suppose that $u\in\pv$ is algebraic over $M$.

Observe that $u\in\M\cdot\pv$ and assume that if
$u\in\M\cdot\pv_{n-1}$ then $u\in\M$(this is our induction
hypothesis). Now, $\M\cdot\pv$ is a Picard-Vessiot (extension by
antiderivatives)extension of $\M\cdot\pv_{n-1}$ and the
differential Galois $\G(\M\cdot\pv|\M\cdot\pv_{n-1})$ is
isomorphic to $(\C,+)^m$ for some $m\in\N$. Note that $u$ is
algebraic over $\M\cdot\pv_{n-1}$ since
$\M\cdot\pv_{n-1}\supseteq\M$. Thus $[\M\cdot\pv_{n-1}\langle
u\rangle, \M\cdot\pv_{n-1}]<\infty.$

Then from the fundamental theorem we should have a finite
algebraic subgroup of
$\G(\M\cdot\pv|\M\cdot\pv_{n-1})\backsimeq(\C,+)^m$ fixing
$\M\cdot\pv_{n-1}\langle u\rangle$. Since the only finite
algebraic subgroup of $(\C,+)^m$ is the trivial group, we obtain
$\M\cdot\pv_{n-1}\langle u\rangle=\M\cdot\pv_{n-1}$ and thus
$u\in\M\cdot\pv_{n-1}$. Now we apply our induction hypothesis to
prove the theorem. \epf

Thus the above theorem shows that if
$\pv\supseteq\M\supsetneq\K\supseteq\f$ are differentials fields
and $\pv$ is a tower of extension by antiderivatives of $\f$ then
$\M$ is purely transcendental over $\K$.

\bt \label{simple version} Let $\pv\supseteq\f$ be a NNC
extension. If there is an $\x\in\pv\setminus\f$ such that
$\x'\in\f$ then for any $n\in\N$ and distinct
$\alpha_1,\cdots,\alpha_n\in\C$, the elements $\y_i\in\f_\infty$
such that $\y'_{\alpha_i}=\frac{1}{\x+\alpha_i}$ are algebraically
independent over $\f(\x)$. Moreover, the differential field
$\f(\y_\alpha,\x)$, where $\y'_\alpha=\frac{1}{\x+\alpha}$ and
$\alpha\in\C$ is not imbeddable in any Picard-Vessiot extension of
$\f$. \et

\bpf

Let there be an $\x\in\pv\setminus\f$ such that $\x'\in\f$.
Suppose that there are distinct constants
$\alpha_1,\cdots,\alpha_n\in\C$ such that the elements $\y_i\in\U$
are algebraically dependent over $\f(\x)$. Since $\y_i$ are
antiderivatives, by the Kolchin-Ostrowski Theorem, there are
constants $c_i\in\C$, not all zero, such that
$\sum^n_{i=1}c_i\y_i\in\f(\x)$. Let $P,Q\in\f[\x]$, $(P,Q)=1$ and
$Q$, a monic polynomial such that
\begin{equation}\sum^n_{i=1}c_i\y_i=\frac{P}{Q}.\end{equation}
Taking the derivative of the above equation, we obtain
\begin{equation*}Q^2(\sum^n_{i=1}\frac{c_i}{\x+\alpha_i})=P'Q-QP'.\end{equation*}
We rewrite the above equation as
\begin{equation}Q^2F=G(P'Q-QP'),\end{equation}
where $F=\sum^n_{i=1}\prod^n_{j=1, i\neq j}c_i(\x+\alpha_j)$ and
$G=\prod^n_{i=1}\x+\alpha_i$. Note that not all $c_i$ are zero and
assume that $c_1\neq 0$. Then $y:=\x+c_1$ divides $G$, $y^2$ does
not divide $G$ and $(F,G)=1$. Thus $y$ divides $Q^2$ and since $y$
is irreducible(and $\f[\x]$ is a UFD) $y$ divides $Q$. Now let
$l\in\N$ be the largest number such that $y^l$ divides $Q$. Then
$y^{2l}$ divides $Q^2$ and therefore $y^{l+1}$ divides
$G(P'Q-Q'P)$. Since $y^l$ divides $Q$, $y$ divides $G$ and
$(P,Q)=1$ we obtain $y^l$ divides $Q'$. Writing $Q=y^lH$ for some
$H\in\f[\x]$ and observing that $Q'=ly'y^{l-1}H+y^lH'$, we obtain
$y$ divides $H$, contradicting the maximality of $l$. Thus the
elements $\y_1,\cdots,\y_n$ are algebraically independent over
$\f(\x)$.

Suppose that there is an $\alpha\in\C$ such that
$\pv\supset\f(\y_{\alpha},\x)\supset\f$ for some Picard-Vessiot
extension $\pv$ of $\f$. Note that $\f(\x)$ is a Picard-Vessiot
sub-extension of $\pv\supset\f$ with differential Galois group
$\G(\f(\x)|\f)\backsimeq(\C,+)$, and every automorphism of
$\f(\x)$ fixing $\f$ lifts to an automorphism of $\pv$ over $\f$.
In particular, there is an automorphism $\sigma\in\G(\pv|\f)$ such
that $\sigma(\x)=\x+c$ for some $c\neq 0$. Observing that
\begin{align*}\y'_\alpha&=\frac{1}{\x+\alpha}\\
\implies\quad\sigma^i(\y_\alpha)'&=\frac{1}{\x+\alpha+ic}
\end{align*}
and that $\sigma\in\G(\pv|\f)$, we obtain
$\y_{\alpha+ic}:=\sigma^i(\y_\alpha)\in\pv$. Since $\alpha+ic$ are
distinct for $1=1,2,\cdots,m$,  the elements
$\y_{\alpha+c},\y_{\alpha+2c},\cdots\in\pv$ are algebraically
independent over $\f$. Thus we obtain a contradiction to the fact
that $\pv$, a Picard-Vessiot extension over $\f$, has a finite
transcendence degree over $\f$.\epf

\begin{remark}Thus if $\pv\supseteq\f$ are differential fields such that
$\x\in\pv\setminus\f$ and $\x'\in\f$ then the differential field
$\f(\y_\alpha,\x)$, $\y'_\alpha=\frac{1}{\x+\alpha}$ and
$\alpha\in\C$ is not imbeddable in any Picard-Vessiot extension of
$\f$ and thus $\y_\alpha\notin\f_1$. We may apply the above
theorem again for the element $\y_\alpha$ with $\f_1$ as the
ground field. Then for any $\z_\beta\in\f_\infty$ such that
$\z'_\beta=\frac{1}{\y_\alpha+\beta}$, $\beta\in\C$, we obtain
that the differential field $\f_1(\z_\beta, \y_\alpha)$ is not
imbeddable in any Picard-Vessiot extension of $\f_1$ and thus
$\z_\beta\notin\f_2$. A repeated application of the theorem proves
the following: If $\f$ is a differential field that has a proper
extension by antiderivatives then for given any $n$, $\f_n$ has
proper extensions by antiderivatives.\end{remark}

Let $\pv\supset\f$ be differential fields and let
$\x_1,\cdots,\x_l\in\pv$ be algebraically independent
antiderivatives of $\f$.

\begin{definition}
An antiderivative $\y$ of $\f(\x_1,\cdots,\x_l)$ is called an
Irreducible-explicit(I-E)antiderivative if $\y'=\frac{A}{CB}$,
where $A,B,C\in\f[\x_1,\cdots,\x_l]$, $(A,B)=$ $(A,B)=$ $(A,B)=1$
and $C$ is an irreducible polynomial.
\end{definition}

\begin{definition}
For each $i=1,2,\cdots,m$ let $\y_i\in\U$ be an antiderivative of
$\frac{A_i}{C_iB_i}$, where
$C_i,A_i,B_i\in\f[\x_1,\cdots,\x_l]$,$(A,B)=$$(A,B)=$$(A,B)=1$,
and satisfying the following conditions;
\begin{itemize}
    \item [C1:]$C_i$ is an irreducible polynomial, $C_i\nmid C_j$ if
    $i\neq j$ and $C_i\nmid B_j$ for any $1\leq i,j\leq m$.
    \item [C2:] for every $1\leq i\leq m$ there is an element $\x_{C_i}\in\{\x_1,\cdots,\x_l\}$
    such that the partial $\frac{\partial C_i}{\partial \x_{C_i}}\neq
    0$ and $\frac{\partial A_i}{\partial \x_{C_i}}=$$\frac{\partial B_i}
    {\partial \x_{C_i}}=0$.
\end{itemize}
We call $\y_1,\cdots,\y_m$ a J-I-E(Joint-Irreducible-Explicit)
antiderivatives of $\f(\x_1,\cdots,\x_l)$. We call the
differential field $\f(\y_1,\cdots,\y_m,\x_1,\cdots,\x_l)$, a
\textsl{2-tower} J-I-E \textit{extension} of $\f$.

\end{definition}

The following theorem shows any set of antiderivatives
$\y_1,\cdots,\y_m$ of $\f(\x_1,\cdots,\x_l)$,
$\y'_i=\frac{A_i}{C_iB_i}$, becomes algebraically independent over
$\f(\x_1,\cdots,\x_l)$ once it satisfies $\textsc{C}1$(see theorem
\ref{alg indp of antideriv}) and thus J-I-E antiderivatives of
$\f(\x_1,\cdots,\x_l)$ are algebraically independent
$\f(\x_1,\cdots,\x_l)$.

\bt\label{alg indp of antideriv} Let $\pv\supseteq\f$ be
differential fields, $\x_1,\cdots,\x_l\in\pv$ be antiderivatives
of $\f$ and assume that $\x_1,\cdots,\x_l$ are algebraically
independent over $\f$. For each $i=1,\cdots,m$ let $A_i,B_i,C_i
\in\f[\x_1,\cdots,\x_l]$, $(A_i,B_i)=$ $(A_i,C_i)=$ $(B_i,C_i)=1$
be polynomials satisfying the following condition

\textsc{C1:} $C_i$ is an irreducible polynomial, $C_i\nmid C_j$ if
    $i\neq j$ and $C_i\nmid B_j$ for any $1\leq i,j\leq m$.

Let $\y_1,\cdots,\y_m\in\U$ be antiderivatives of
$\f(\x_1,\cdots,\x_l)$ with $\y'_i=\frac{A_i}{C_iB_i}$. Then
$\y_1,\cdots,\y_m$ are algebraically independent over
$\f(\x_1,\cdots,\x_l)$. \et

\bpf

Suppose that $\y_1,\cdots,\y_m$ are algebraically dependent over
$\f(\x_1,\cdots,\x_l)$. Then the Kolchin-Ostrowski theorem
guarantees constants $\alpha_1,\cdots,c_m\in\C$, not all zero,
such that $\sum^m_{i=1}\alpha_i\y_i\in\f(\x_1,\cdots,\x_l)$.
Assume that $\alpha_1\neq 0$.

First we note that if $\sum^m_{i=1}\alpha_i\y_i\in\C$ then
$\sum^m_{i=1}\alpha_i\frac{A_i}{C_iB_i}=0$ and now writing
$\sum^m_{i=2}\alpha_i\frac{A_i}{C_iB_i}=\frac{F}{G}$,
$F,G\in\f[\x_1,\cdots,\x_l]$, we obtain
\begin{align*}&\alpha_1\frac{A_1}{C_1B_1}=-\frac{F}{G}\\
\implies&\alpha_1A_1G=-FC_1B_1.
\end{align*}
Since $A_1\neq 0$, we obtain $F\neq 0$ and thus we may assume $F$
and $G$ are relatively prime polynomials. Clearly, $C_1$ divides
$A_1G$ and since $A_1$ and $C_1$ are relatively prime, $C_1$
divides $G$. On the other hand
$\sum^m_{i=2}\alpha_i\frac{A_i}{C_iB_i}=\frac{F}{G}$ implies $G$
divides $\prod^m_{i=2}C_iB_i$, which implies $C_1$ divides
$\prod^m_{i=2}C_iB_i$ contradicting the condition \textsc{C1}.
Thus $\sum^m_{i=1}\alpha_i\y_i\in\f(\x_1,\cdots,\x_l)\setminus\C$.

Let $P,Q\in\f[\x_1,\cdots,\x_l]$ be relatively prime polynomials
such that \begin{equation}\label{fun
eqn}\sum^m_{i=1}\alpha_i\y_i=\frac{P}{Q}.\end{equation}

Let  $S,T\in\f[\x_1,\cdots,\x_l]$ be polynomials such that
$\frac{S}{T}=(\frac{P}{Q})'$
$=\sum^m_{i=1}\alpha_i\frac{A_i}{C_iB_i}$. We know that
$\sum^m_{i=1}\alpha_i\y_i\notin\C$ and therefore $S\neq 0$ and
thus we may assume $S$ and $T$ are relatively prime. Since
$\alpha_1\neq 0$, we see that $C_1$ divides $T$. And, $T$ divides
$\prod^m_{i=1}C_iB_i$ and $C_i, B_i$ satisfies condition
\textsc{C1} implies $C^2_1$ does not divide $T$. Thus $P,Q,S$ and
$T$ satisfies the hypothesis of theorem \ref{anti deriv in field
cri}. But, taking the derivative of equation \ref{fun eqn} we
obtain $(\frac{P}{Q})'=\frac{S}{T}$, which contradicts theorem
\ref{anti deriv in field cri}. \epf


\section{Differential Subfields of the J-I-E Tower}\label{the jie tower}

In the next section we will prove a structure theorem for the
differential subfields of a certain tower of extensions by
antiderivatives, namely J-I-E extensions.  These towers are made
by adjoining antiderivatives that appears in theorem \ref{alg indp
of antideriv}.

As usual, let $\C$ be an algebraically closed-characteristic zero
field, $\f$ be a differential field with field of constants $\C$
and let $\f_\infty$ be a complete Picard-Vessiot closure with $\C$
as its field of constants.

\subsection{Automorphisms of J-I-E towers}

Let $\y_{11},\cdots,\y_{1n_1}$ be algebraically independent
antiderivatives of $\f$ and for $i=1,2,\cdots,k,$ let
$\pv_i:=\pv_{i-1}(\y_{i1},\y_{i2},\cdots,\y_{in_i})$, where
$\pv_0:=\f$ and for $i\geq 2$ $\y_{i1},\y_{i2},\cdots,\y_{in_i}$
are I-E antiderivatives of $\pv_{i-1}$, that is,
$\y'_{ij}=\frac{A_{ij}}{C_{ij}B_{ij}}$ and for each $2\leq i\leq
k$ and for all $1\leq j\leq n_i$,
$A_{ij},B_{ij}$,$C_{ij}\in\pv_{i-2}[\y_{i-1 1},\cdots,\y_{i-1
n_{i-1}}]$ are polynomials such that $(A_{ij}, B_{ij})$ $=(B_{ij},
C_{ij})$ $=(A_{ij}, C_{ij})=1$ and satisfying conditions
\textsc{C1} and \textsc{C2}. Let $I_i:=\{\y_{ij}|1\leq j\leq
n_i\}$, $\Lambda_t:=Span_\C\cup^t_{i=1}I_i$, $\Lambda_0=\{0\}$ and
$\pv:=\pv_k$. We will also recall the conditions \textrm{C1} and
\textrm{C2} here

\begin{itemize}
    \item [C1:]$C_{ij}$ is an irreducible polynomial for each $i,j$. For every $i$, $C_{is}\nmid C_{it}$ (that is, they are non associates)if
    $s\neq t$ and $C_{is}\nmid B_{it}$ for any $1\leq s,t\leq n_i$.
    \item [C2:] for every $1\leq j\leq n_i$ there is an element $\y_{C_{ij}}\in\{\y_{i-1 1},\cdots,\y_{i-1
n_{i-1}}\}$
    such that the partial $\frac{\partial C_{ij}}{\partial \y_{C_{ij}}}\neq
    0$ and $\frac{\partial A_{ij}}{\partial \y_{C_{ij}}}=$$\frac{\partial B_{ij}}
    {\partial \y_{C_{ij}}}=0$.
\end{itemize}

\begin{definition}
We call
\begin{equation}\pv:=\pv_k\supset\pv_{k-1}\supset\cdots\supset\pv_2\supset\pv_1\supset\pv_0:=\f\end{equation}
a tower of extensions by J-I-E antiderivatives. Note that $\pv_1$
is an ordinary antiderivative extension of $\f$.
\end{definition}

Let $\G_\infty:=\G(\f_\infty|\f)$, the group of all differential
automorphisms of the complete Picard-Vessiot closure $\f_\infty$
of $\f$. We will show that the group of differential automorphisms
$\G(\pv|\f)$ is isomorphic to the additive group $(\C,+)^\delta$
for some $\delta\leq $ tr.d $\pv|\f$. Moreover, the action of
$\G(\pv|\f)$ on $\pv$ is given by
$\sigma(\y_{ij})=\y_{ij}+c_{ij\sigma}$, $c_{ij\sigma}\in\C$.

\begin{lemma}

For any $\sigma\in\G_\infty$ and $t\geq 2$ the elements of $I_t$,
namely, $\y_{t1},\y_{t2},\cdots,\y_{tn_t}$ are J-I-E
antiderivatives of the differential field
$\pv_{t-1}(\sigma(\pv_{t-1}))$, the compositum of differential
fields $\pv_{t-1}$ and $\sigma(\pv_{t-1})$.
\end{lemma}

\bpf  We observe that
$\pv_{t-1}(\sigma(\pv_s))=\pv_{t-1}(\cup^{s}_{i=1}\sigma(I_i))$
and since $\sigma(\y_{1j})=\y_{1j}+c_{j\sigma}$,
$\pv_{t-1}(\sigma(I_1))=\pv_{t-1}$.  For $2\leq s\leq t-1,$ let
$I^\sigma_s\subset \sigma(I_s)$ be a transcendence base of the
differential field $\pv_{t-1}(\sigma(\pv_s))$ over
$\pv_{t-1}(\sigma(\pv_{s-1}))$. Note that $\sigma(I_s)$ consists
of antiderivatives of $\pv_{t-1}(\sigma(\pv_{s-1}))$ and that
$\pv_{t-1}(\sigma(\pv_{s-1}))$ $(\sigma(I_s))$
$=\pv_{t-1}(\sigma(\pv_{s}))$. Thus $\pv_{t-1}(\sigma(\pv_s))$ is
an extension by antiderivatives of $\pv_{t-1}(\sigma(\pv_{s-1}))$
and therefore
$\pv_{t-1}(\sigma(\pv_{s-1}))(I^\sigma_s)=\pv_{t-1}(\sigma(\pv_{s}))$
for each $1\leq s\leq t-1$.

Thus
$\pv_{t-1}(\sigma(\pv_{t-1}))=\pv_{t-1}(\cup^{t-1}_{i=1}I^\sigma_i)$.
Since $\pv_{t-1}=\pv_{t-2}(\y_{t-1 1},\cdots, \y_{t-1 n_{t-1}})$,
$\y_{t-1 1},\cdots, \y_{t-1 n_{t-1}}$ are algebraically
independent over $\pv_{t-2}$(because they are J-I-E
antiderivatives) and the set $\cup^{t-1}_{i=1}I^\sigma_i$ is
algebraically independent over $\pv_{t-1}$, we obtain that
$\y_{t-1 1},\cdots, \y_{t-1 n_{t-1}}$ are algebraically
independent over $\pv_{t-2}(\cup^{t-1}_{i=1}I^\sigma_i)$. Also
note that
$$\pv_{t-2}(\cup^{t-1}_{i=1}I^\sigma_i)=\pv_{t-2}\sigma(\pv_{t-2})(I^\sigma_{t-1})$$
and that the elements of $I^\sigma_{t-1}$ are antiderivatives of
$\pv_{t-2}\sigma(\pv_{t-2})$. Thus
$\pv_{t-2}(\cup^{t-1}_{i=1}I^\sigma_i)$ is a differential field
which is also a fraction field of the polynomial ring
$\pv_{t-2}[\cup^{t-1}_{i=1}I^\sigma_i]$.

We will now show that $\y_{t 1}, \y_{t, 2},\cdots,$ $\y_{t, n_t}$
are J-I-E antiderivatives of the compositum $\pv_{t-1}$
$\sigma(\pv_{t-1})$.
Since $\y_{t 1}, \y_{t,2},\cdots,\y_{t, n_t}$ are J-I-E
antiderivatives of $\pv_{t-1}$, there are polynomials $A_{t-1 j},
B_{t-1 j}$, $ C_{t-1 j}\in$ $\pv_{t-2}[\y_{t-1 1}, \y_{t-1,
2},\cdots,$ $\y_{t-1, n_{t-1}}]$ such that $(A_{t-1 j},B_{t-1
j})=$ $(B_{t-1 j},C_{t-1 j})$ $=(A_{t-1 j},C_{t-1 j})=1$ and
satisfying conditions \textrm{C1} and \textrm{C2}. We observe that
all the above conditions on $A_{t-1 j}, B_{t-1 j}$ and $C_{t-1 j}$
holds in the polynomial ring $\pv_{t-2}[\cup^{t-1}_{i=1}$
$I^\sigma_i,\y_{t-1 1},$ $ \y_{t-1, 2},$ $\cdots,\y_{t-1,
n_{t-1}}]$ as well and therefore by ``Gauss' lemma" these
conditions hold in the ring
$$\pv_{t-2}(\cup^{t-1}_{i=1}I^\sigma_i)[\y_{t-1 1}, \y_{t-1,
2},\cdots,\y_{t-1, n_{t-1}}].$$ Thus $\y_{t 1}, \y_{t,
2},\cdots,\y_{t, n_t}$ become J-I-E antiderivatives of the field
$$\pv_{t-2}(\cup^{t-1}_{i=1}I^\sigma_i, \y_{t-1 1}, \y_{t-1,
2},\cdots,\y_{t-1, n_{t-1}})=\pv_{t-1}\sigma(\pv_{t-1}).$$ \epf

\bt \label{the big lemma} Let $\M$ be a differential subfield of
$\f_\infty$, $\x_1,\cdots,\x_l\in\f_\infty$ be algebraically
independent antiderivatives of $\M$ and for $i=1,2,\cdots,m$ let
$\y_i\in\f_\infty$ be J-I-E antiderivatives of
$\M(\x_1,\cdots,\x_l)$ ( that is, $\y'_i=\frac{A_i}{B_iC_i}$,
where $A_i,B_i$ and $C_i$ satisfies conditions
$(A_i,B_i)=$$(B_i,C_i)=$$(A_i,C_i)=1$, \textsc{C1} and
\textsc{C2}). Suppose that there is a subgroup $\sg$ of
$\G_\infty$ of differential automorphisms fixing $\M$ and an
element $\s:=\sum^e_{i=1}\alpha_i\y_i\in\M(\y_1,\cdots,\y_m,$
$\x_1,\cdots,\x_l)$, $\alpha_i\in\C$ such that for every
$\sigma\in\sg$, $\sigma(\s)\in$ $\M(\y_1,\cdots,\y_m,$
$\x_1,\cdots,\x_l)$. Then every $\sigma\in\sg$ fixes $A_i,B_i$ and
$\C_i$ whenever $\alpha_i\neq 0$, that is
$\sigma(\y_i)=\y_i+c_{i\sigma}$, for some $c_{i\sigma}\in\C$. In
particular, for every $\sigma\in\sg$ there is a
$c_\sigma:=\s(c_{1\sigma},\cdots,c_{l\sigma})\in\C$ such that
$\sigma(\s)=\s+c_\sigma$. \et

\bpf If $\sg$ is the trivial group then the proof is trivial.
Assume that $\sg$ is a nontrivial group. Since $\x'_i\in\M$,
$\M(\x_1,\cdots,\x_l)$ is an extension by antiderivatives of $\M$
and thus the differential field $\M(\x_1,\cdots,\x_l)$ is
preserved by $\sg$. In particular $\sigma(\frac{A_i}{C_iB_i})\in$
$\M(\x_1,\cdots,\x_l)$. Then $\M(\y_1,$
$\cdots,\y_m,\x_1,\cdots,\x_l)$ has $m+1$ antiderivatives
$\sum^e_{i=1}\alpha_i\sigma(\y_i)$, $\y_1,\cdots,\y_m$ of
$\M(\x_1,$ $\cdots,\x_l)$ and therefore the antiderivatives has to
be algebraically dependent over $\M(\x_1,\cdots,\x_l)$. Now from
the Kolchin-Ostrowski theorem we have constants $\gamma_i$, $1\leq
i\leq m+1$ not all zero such that
\begin{equation}\label{new antideriv}\sum^m_{i=1}\gamma_i\y_i+\gamma_{m+1}
\sum^e_{i=1}\alpha_i\sigma(\y_i)\in\M(\x_1,\cdots,\x_l).\end{equation}

Note that if $\gamma_{m+1}=0$ then $\y_i$'s become algebraically
dependent over $\M(\x_1,\cdots,\x_l)$, which is not true. so
$\gamma_{m+1}\neq 0$ and thus we may assume
$\gamma_{m+1}=1$(dividing through the equation \ref{new antideriv}
by $\gamma_{m+1}$).

First we will show that $\sigma(C_i)=C_i$ for all $\sigma\in\sg$
whenever $\alpha\neq 0$. Then we will use this to show that $\sg$
indeed fixes $A_i$ as well as $B_i$ whenever $\alpha_i\neq 0$.

Suppose that there is a $\rho\in\sg$ and an $i, 1\leq i\leq m$
such that $\alpha_i\neq 0$ and $\rho(C_i)\neq C_i$. For
convenience, let us assume that $i=1$. The automorphism $\rho$
acts on the ring $\M[\x_1,\cdots,\x_l]$ by sending
$\x_i\to\x_i+c_{i\rho}$ and if $\rho$ is nontrivial then clearly
$\rho$ has an infinite order. Thus we have $\rho(C_1)=$
$C_1(\x_1+c_{1\rho},\cdots,\x_l+c_{l\rho})$. From proposition
\ref{trans poly prop} we see that $C_1$ divides $\rho(C_1)$ only
if $C_1=\rho(C_1)$ and thus $\rho(C_1)$ and $C_1$ are not
associates (over $\M$). In fact, for any $i,j\in\N\cup\{0\}$,
$i\neq j$, the elements $\rho^i(C_1)$ and $\rho^j(C_1)$ are
non-associates. Since every polynomial in $\M[\x_1,\cdots,\x_l]$
has finitely many (non-associate) irreducibles and $\rho^i(C_1)$
is also an irreducible for each $i\in\N$, there is a $j\in\N$ such
that
$$\rho^j(C_1)\nmid B_1B_2\cdots
B_m.$$ We also note that $\rho^j(C_1)\nmid\rho^j(B_j)$ for any
$1\leq j\leq m$ and $\rho^j(C_1)\nmid\rho^j(C_i)$ for any $i\neq
1$; otherwise $C_1|B_j$, or $C_1|C_i$ for some $i\neq 1$ and in
either case, contradicts the condition \textsc{C1}. Thus
\begin{equation}\label{rho equation c_1}\rho^j(C_1)\ \ \text{does not
divide}\ \
B_1\prod^m_{i=2}B_i\rho^j(B_i)\rho^j(C_i).\end{equation}

The equation \ref{new antideriv} is true for all $\sigma\in\sg$
and thus there are polynomials $A,B\in\M[\x_1,\cdots,\x_l]$
$$\sum^m_{i=1}\gamma_i\y_i+
\sum^e_{i=1}\alpha_i\rho^j(\y_i)=\frac{A}{B}.$$

Let $S, T\in\M[\x_1,\cdots,\x_l]$ be relatively prime polynomials
such that
\begin{equation}\label{sf eqn}
\alpha_1\frac{\rho^j(A_1)}{\rho^j(C_1)\rho^j(B_1)}+\sum^m_{i=1}\gamma_i\frac{A_i}{C_iB_i}+
\sum^e_{i=2}\alpha_i\frac{\rho^j(A_i)}
{\rho^j(C_i)\rho^j(B_i)}=\frac{S}{T}\end{equation}

and let $F,G\in\M[\x_1,\cdots,\x_l]$ be relatively prime
polynomials such that
\begin{equation}\label{G factor eqn}\frac{F}{G}=-\sum^m_{i=1}\gamma_i\frac{A_i}{C_iB_i}+
\sum^e_{i=2}\alpha_i\frac{\rho^j(A_i)}
{\rho^j(C_i)\rho^j(B_i)}.\end{equation} Note that
\begin{equation}\label{G div eqn}G \ \ \text{divides}\ \ B_1\prod^m_{i=2}B_i\rho^j(B_i)\rho^j(C_i)\end{equation}

Suppose that $S=0$. Then
\begin{align}\label{imp rho eqn}&\alpha_1\frac{\rho^j(A_1)}{\rho^j(C_1)\rho^j(B_1)}=\frac{F}{G}\notag\\
\implies&\alpha_1\rho^j(A_1)G=\rho^j(C_1)\rho^j(B_1)F.\end{align}
 Since $A_1$ is a non zero
polynomial, so is $\rho^j(A_1)$ and thus $\alpha_1\neq 0$ implies
$F\neq 0$. From equation \ref{imp rho eqn} we obtain $\rho^j(C_1)$
divides $G$ and now equation \ref{G div eqn} contradicts equation
\ref{rho equation c_1}.

Thus $S\neq 0$. Substituting equation \ref{G factor eqn} in
equation \ref{sf eqn} we obtain
\begin{align}\label{rhoc1 div t}\alpha_1\frac{\rho^j(A_1)}{\rho^j(C_1)\rho^j(B_1)}-\frac{F}{G}&=\frac{S}{T}\notag\\
\left(\alpha_1\rho^j(A_1)G-\rho^j(C_1)\rho^j(B_1)F\right)T&=SG\rho^j(C_1)\rho^j(B_1).\end{align}

From the above equation \ref{rhoc1 div t} we obtain $\rho^j(C_1)$
divides $\alpha_1\rho^j(A_1)GT$. Again equations \ref{G div eqn}
and \ref{rho equation c_1} guarantees $\rho^j(C_1)$ does not
divide $G$ and clearly $\rho^j(C_1)$ does not divide
$\rho^j(A_1)$. Therefore $\rho^j(C_1)$ divides $T$, which implies
that $\rho^j(C_1)$ is an irreducible factor of $T$. Thus we have
produced polynomials $A,B,S,T\in\M[x_1,\cdots,\x_n]$ contradicting
theorem \ref{anti deriv in field cri}. Hence $\sigma(C_i)=C_i$ for
all $\sigma\in\sg$.

Now we will show that $\sg$ fixes $A_i$ and $B_i$ for every $i$.

Assume that $\alpha_1\neq 0$ and pick a $\sigma\in\sg$. Note that
$\sigma(C_1)=C_1$ and that $\sigma$ is an automorphism, therefore
$C_1\neq\sigma(C_j)$ for any $j\neq 1$. If
$P\in\M[\x_1,\cdots,\x_l]$ is a polynomial and $C_1$ divides
$\sigma(P)$ then $\sigma^{-1}(C_1)$ divides $P$. But
$\sigma(C_1)=C_1$ implies $\sigma^{-1}C_1=C_1$ and therefore $C_1$
divides $P$. Hence we note that
\begin{equation}\label{sigma equation c_1}C_1\ \ \text{does not
divide}\ \
B_1\prod^m_{i=2}B_i\sigma(B_i)\sigma(C_i).\end{equation}

Take the derivative of equation \ref{new antideriv} to obtain
\begin{equation}\label{cru rho eqn2}B^2\left(\sum^m_{i=2}\gamma_i\frac{A_i}{C_iB_i}+
\sum^e_{i=2}\alpha_i\frac{\sigma(A_i)}
{\sigma(C_i)\sigma(B_i)}+\gamma_1\frac{A_1}
{C_1B_1}+\alpha_1\frac{\sigma(A_1)}
{C_1\sigma(B_1)}\right)=BA'-AB'.\end{equation}

Let $F,G\in\M[\x_1,\cdots,\x_l]$ be relatively prime polynomials
such that
\begin{equation}\label{G factors in b c}\sum^m_{i=2}\gamma_i\frac{A_i}{C_iB_i}+
\sum^e_{i=2}\alpha_i\frac{\sigma(A_i)}
{\sigma(C_i)\sigma(B_i)}=\frac{F}{G},\end{equation} and let
$S,T\in\M[\x_1,\cdots,\x_l]$ be relatively prime polynomials such
that \begin{equation}\label{crude s t}\gamma_1\frac{A_1}
{C_1B_1}+\alpha_1\frac{\sigma(A_1)}
{C_1\sigma(B_1)}+\frac{F}{G}=\frac{S}{T}\end{equation}

Note that $(\frac{A}{B})'=\frac{S}{T}$ and that
\begin{equation}\label{G div eqn sig fix}G \ \text{divides}\
\prod^m_{i=2}B_i\sigma(B_i)\sigma(C_i).\end{equation} We rewrite
equation \ref{crude s t} as
\begin{equation}\label{std c t b eqn}TG\left(\gamma_1A_1\sigma(B_1) +\alpha_1\sigma(A_1)B_1\right)
+TFC_1\sigma(B_1)B_1=SGC_1\sigma(B_1)B_1\end{equation}

Again, we will split our into two cases; $S\neq 0$ and $S=0$. In
both the cases, we will show that $C_1$ divides
$\gamma_1A_1\sigma(B_1) +\alpha_1\sigma(A_1)B_1$. Assume for a
moment that we proved $C_1$ divides $\gamma_1A_1\sigma(B_1)
+\alpha_1\sigma(A_1)B_1$. Then from \textsc{C2} we have
$\x_{C_1}\in\{\x_1,\cdots,\x_l\}$ such that $\frac{\partial
C_1}{\partial \x_{C_1}}\neq
    0$ and $\frac{\partial A_1}{\partial \x_{C_1}}=$$\frac{\partial B_1}
    {\partial \x_{C_1}}=0$. Since $\sigma(\x_i)=\x_i+c_{i\sigma}$ for some $c_{i\sigma}\in\C$, $\sigma$
    is an automorphism of the ring
    $\M[\{\x_1,\cdots,\x_l\}\setminus\{\x_{C_1}\}]$ and therefore $\gamma_1A_1\sigma(B_1)$
    $+\alpha_1\sigma(A_1)B_1\in\M[\{\x_1,\cdots,\x_l\}\setminus\{\x_{C_1}\}]$.
    Thus $C_1$
divides $\gamma_1A_1\sigma(B_1) +\alpha_1\sigma(A_1)B_1$ implies
$\gamma_1A_1\sigma(B_1) +\alpha_1\sigma(A_1)B_1=0$, that is,
$\sigma\left(\frac{A_1}{B_1}\right)=-\frac{\gamma_1}{\alpha_1}\frac{A_1}{B_1}$.
Then $A_1$ divides $\sigma(A_1)$ and $B_1$ divides $\sigma(B_1)$
and therefore from proposition \ref{trans poly prop} we obtain
$\sigma(A_1)=A_1$ and $\sigma(B_1)=B_1$.

Let us show that $C_1$ divides $\gamma_1A_1\sigma(B_1)
+\alpha_1\sigma(A_1)B_1$.

Case $S\neq 0$:

From equation \ref{std c t b eqn} we observe that $C_1$ divides
$TG\left(\gamma_1A_1\sigma(B_1) +\alpha_1\sigma(A_1)B_1\right)$
and from equations \ref{G div eqn sig fix} and \ref{sigma equation
c_1} that $C_1$ does not divide $G$ and therefore  $C_1$ has to
divide $T(\gamma_1A_1\sigma(B_1)$ $ +\alpha_1\sigma(A_1)B_1)$. If
$C_1$ divides $T$ then the polynomials
$A,B,S,T\in\M[\x_1,\cdots,\x_l]$ contradicts theorem \ref{anti
deriv in field cri}. Thus $C_1$ divides $\gamma_1A_1\sigma(B_1)
+\alpha_1\sigma(A_1)B_1$.

Case $S=0$: From equation \ref{std c t b eqn} we have
$$G\left(\gamma_1A_1\sigma(B_1) +\alpha_1\sigma(A_1)B_1\right)
=-FC_1\sigma(B_1)B_1.$$ As noted earlier, $C_1$ does not divide
$G$ and thus $C_1$ divides $\gamma_1A_1\sigma(B_1)
+\alpha_1\sigma(A_1)B_1$.

Thus we see that for every $\sigma\in\sg$, $\sigma(A_i)=A_i$,
$\sigma(B_i)=B_1$ and $\sigma(C_i)=C_i$ and therefore
$$\left(\sigma(\y_i)\right)'=\sigma\left(\frac{A_i}{C_iB_i}\right)=\frac{A_i}{C_iB_i}.$$
Since $\y'_i=\frac{A_i}{C_iB_i}$, we obtain
$\sigma(\y_i)=\y_i+c_{i\sigma}$ for some $c_{i\sigma}\in\C$.
Clearly, for every $\sigma\in\sg$, $\sigma(\s)=\s+c_\sigma$ where
$c_\sigma:=\s(c_{1\sigma},\cdots,c_{l\sigma})\in\C$.\epf

Before we classify the differential subfields of a general J-I-E
tower we will first work with a two step tower.

\bt \label{gen anti deriv thm} Let $\f(\x_1,\cdots,\x_l)\supset\f$
be an extension by algebraically independent antiderivatives
$\x_1,\cdots,\x_l$ of $\f$. Let $\y_1,\cdots,\y_m$ be J-I-E
antiderivatives of $\f(\x_1,\cdots,\x_l)$. Then every differential
subfield of $\f(\y_1,\cdots,\y_m,$ $\x_1,\cdots,\x_m)$ is of the
form $\f(\rs,\rt)$, where $\rs$ and $\rt$ are finite subsets of
$span_\C\{\y_1,$ $\cdots,\y_m,$ $ \x_1,\cdots,\x_l\}$ and
$span_{\C}\{\x_1,\cdots,\x_m\}$ respectively. \et

\bpf

Let $\pv:=\f(\y_1,\cdots,\y_m,\x_1,\cdots,\x_l)$,
$\IL:=\f(\x_1,\cdots,\x_l)$ and $\pv\supseteq\K\supseteq\f$ be an
intermediate differential field. Note that $\IL$ is an extension
by antiderivatives of $\f$ and $\IL\supseteq\K\cap\IL\supseteq\f$
is an intermediate subfield. Thus there is a finite set
$\rt\subset span_\C\{\x_1,\cdots,\x_l\}$, algebraically
independent over $\f$ such that $\K\cap\IL=\f(\rt)$. Let
$\overline{\rt}\subset\{\x_1,\cdots,\x_l\}$ be a transcendence
base of $\IL$ over $\f(\rt)$. We observe that
$\f(\overline{\rt},\rt)=\IL$, $|T|+|\overline{T}|=l$, and
$\overline{\rt}$ is algebraically independent over $\K$;
otherwise, $\overline{\rt}$ becomes algebraically dependent over
$\K\cap\IL=\f(\rt)$ which contradicts the choice of
$\overline{\rt}$.

Thus $\K(\x_1,\cdots,\x_l)=\K(\overline{\rt})$. We observe that
$\pv\supseteq\K(\overline{\rt})\supset\IL$ and that $\pv$ is a
(Picard-Vessiot) extension by antiderivatives of $\IL$. Thus there
is a finite set $\rs^\sharp \subset$
$span_{\C}\{\y_1,\cdots,\y_m\}$ such that
$\K(\overline{\rt})=\IL(\rs^\sharp)$. We may also assume that
$\rs^\sharp$ is algebraically independent over $\IL$. Since
$\K(\overline{\rt})$ is a (Picard-Vessiot) extension by
antiderivatives of $\K$, for every $\s\in \rs^\sharp$ and
$\rho\in\G:=\G(\K(\overline{\rt})|\K)$, the element
$\rho(\s)\in\K(\overline{\rt})$. Thus $\rho(\s)\in\pv$ for every
$\rho\in\G$ and for every $\s\in \rs^\sharp$.

We have
$$\begin{CD}
\IL @> >> \K(\overline{\rt})\\
@V\  VV @VV  V\\
\f @> >> \K
\end{CD}$$
where are arrows are inclusions. Thus there is a natural injective
map $\phi:\G(\K(\overline{\rt})|K)$$\to$$\G(\IL|\f)$ of algebraic
groups such that $\rho(\x_i)=\phi(\rho)(\x_i)$ for all
$\rho\in\G(\K(\overline{\rt})|K)$, and there is an algebraic
subgroup $\sg$ of $\G(\IL|\f)$ such that the image
$\phi(\G(\K(\overline{\rt})|K))=\sg$. Note that the action of
$\rho$ on $\x_i$ completely determines $\rho$ for all
$\rho\in\G(\K(\overline{\rt})|K)$.

Thus $\sigma(\s)\in\pv$ for every $\sigma\in\sg$ and for every
$\s\in \rs^\sharp$. Now from theorem \ref{the big lemma} we obtain
$\sigma(\s)=\s+c_\sigma$ for all $\sigma\in\sg$, $c_\sigma\in\C$.
Thus $\s'\in\IL^\sg$ and in particular $\sigma(\s')=\s'$ for all
$\sigma\in\sg$. Since $\rho(\x_i)=\phi(\rho)(\x_i)$ for all
$\rho\in\G$ and $\phi$ is surjective, $\rho(\s')=\s'$ for every
$\rho\in\G$ and therefore $\s'\in\K^{\G}=\K$. Then
$\s\in\K(\overline{\rt})$ is an antiderivative of $\f$ and
therefore the set $\overline{\rt}\cup\{\s\}$ has to be
algebraically dependent over $\K$. From The Kolchin-Ostrowski
theorem, there is an element $\ft_\s\in span_\C\overline{\rt}$
such that $\s+\ft_\s\in\K$. We also observe that $\s'\in\K$ and
$\s'\in\IL$ and therefore $\s'\in\f(\rt)$. Now we let
$\rs:=\{\s+\ft_\s|\s\in \rs^\sharp\}$ and observe that
$\K\supset\f(\rs,\rt)\supset\f(\rt)\supset\f$. Let
$\overline{\rs}\subset\{\y_1,\cdots,\y_m\}$ be a transcendence
base of $\pv$ over $\K(\overline{\rt})=\IL(\rs^\sharp)$. Then
$|\overline{\rs}|+|\rs^\sharp|=m$ and in particular
$\IL(\overline{\rs}\cup \rs^\sharp)=\pv$.

We know that \begin{equation}\label{transcendence deg
arg}\textrm{tr.d} \ \pv|\f=\textrm{tr.d}\ \pv|\K+\textrm{tr.d}\
\K|\f(S, T)+\textrm{tr.d}\ \f(S,T)|\f\end{equation}
tr.d$\pv|\K=|\overline{\rs}|+|\overline{\rt}|$ and tr.d
$\f(\rs,\rt)|\f=|\rs|+|\rt|$. Note that $|\rs|=|\rs^\sharp|$ and
that $|\overline{\rs}|+|\overline{\rt}|+|\rs^\sharp|+|\rt|=$ tr.d
$\pv|\f=l+m$. Thus tr.d $\pv|\f=|\overline{\rs}|+|\overline{\rt}|$
$+|\rs|+|\rt|=$ tr.d $\pv|\K+$tr.d $\f(S,T)|\f$ and therefore from
equation \ref{transcendence deg arg} we obtain tr.d $\K|\f(S,
T)=0$. Thus $\K$ is algebraic over $\f(\rs,\rt)$. Now letting
$\M:=\f(\rs,\rt)$ and applying theorem \ref{gen anti deriv thm},
we obtain $\K=\f(\rs,\rt)$. \epf

\bt\label{group  action on ie towers} If there is an
$\s=\sum_{j=1}^{n_t}\alpha_{tj}\y_{tj}+\sum^{t-1}_{i=1}\sum^{n_i}_{j=1}\alpha_{ij}\y_{ij}\in\Lambda_t\setminus\Lambda_{t-1}$
for some $1\leq t\leq k$ and a subgroup $\bh$ of
$\G(\f_\infty|\f)$ such that for every $\sigma\in\bh$,
$\sigma(\s)\in\pv_k=:\pv$ then
$\sigma(\y_{ij})=\y_{ij}+c_{ij\sigma}$ for every $\sigma\in\sg$
provided the coefficient $\alpha_{ij}$ of $\y_{ij}$ in $\s$ is
nonzero.

\et

\bpf

We will use an induction on $t$ to prove this theorem.

$t=1$: Then $\s$ is a linear combination of antiderivatives
$\y_{11},\cdots,\y_{1n_1}$ of $\f$. Therefore for every
$\sigma\in\G(\f_\infty|\f)$ we have
$$\y'_{1j}=\sigma(\y'_{1j})=\sigma(\y_{1j})'.$$
Since $\f_\infty$ and $\f$ has the same field of constants, there
is a $c_{1j\sigma}\in\C$ such that
$\sigma(\y_{1j})=\y_{1j}+c_{1j\sigma}$.

Assume that our theorem is true for $t-1$.

$t\geq 2$: For
$$\s=\sum_{j=1}^{n_t}c_{tj}\y_{tj}+\sum^{t-1}_{i=1}\sum^{n_i}_{j=1}c_{ij}\y_{ij},$$
where $\alpha_{tj}\neq 0$ for some $j$, suppose that
$\sigma(\s)\in\pv$. Then
\begin{align}
\sigma(\s)=&\sum_{j=1}^{n_t}c_{tj}\sigma(\y_{tj})+\sum^{t-1}_{i=1}\sum^{n_i}_{j=1}c_{ij}\sigma(\y_{ij})\in\pv\\
\implies&\sum_{j=1}^{n_t}c_{tj}\sigma(\y_{tj})\in\pv(\sigma(\pv_{t-1}));
\ \text{since}\
\sum^{t-1}_{i=1}\sum^{n_i}_{j=1}c_{ij}\sigma(\y_{ij})\in\sigma(\pv_{t-1})
\end{align}

Suppose that for $i\geq t+1$,
$\sigma(\s)\in\pv_i(\sigma(\pv_{t-1}))$. Then note that
$\pv_i(\sigma(\pv_{t-1}))$ is an extension by algebraically
independent antiderivatives $\y_{i1},\cdots,\y_{in_i}$ of
$\pv_{i-1}(\sigma(\pv_{t-1}))$. Also note that $\sigma(\s)$ is an
antiderivative of $\sigma(\pv_{t-1})$ and therefore an
antiderivative of $\pv_{i}(\sigma(\pv_{t-1}))$. Thus there are
constants $\alpha_{i0},\alpha_{ij}\in\C$, $1\leq j\leq n_i$ not
all zero such that
$$\alpha_{i0}\sigma(\s)+\sum^{n_i}_{j=1}\alpha_{ij}\y_{ij}\in\pv_{i-1}(\sigma(\pv_{t-1})).$$
But, if $\alpha_{ij}\neq 0$ for some $1\leq j\leq n_i$ then from
the above equation and from the facts that
$\sigma(\s)\in\sigma(\pv_t)$ and
$\sigma(\pv_{t-1})\subset\sigma(\pv_t)$ we have
\begin{align*}\sum^{n_i}_{j=1}\alpha_{ij}\y_{ij}\in\pv_{i-1}(\sigma(\pv_t))\end{align*}
and since $t\leq i-1$, $\pv_t\subset\pv_{i-1}$, which implies
$$\sum^{n_i}_{j=1}\alpha_{ij}\y_{ij}\in\pv_{i-1}(\sigma(\pv_{i-1})),$$ a
contradiction to theorem...

Thus $\sigma(\s)\in\pv(\sigma(\pv_{t-1}))$ implies
$\sigma(\s)\in\pv_t(\sigma(\pv_{t-1}))$. Let
$\M:=\pv_{t-2}(\sigma(\pv_{t-2}))$. We know that
$I_{t-1}=\{\y_{t-1, 1},\cdots,\y_{t-1, n_{t-1}}\}$ is
algebraically independent over $\M$. Now let
$I^\sigma_{t-1}\subset \sigma(I_{t-1})$ be a transcendence base of
$\pv_{t-1}(\sigma(\pv_{t-1}))$ over $\M(I_{t-1})$. Then
$\M(I_{t-1},I^\sigma_{t-1})=$ $\pv_{t-1}(\sigma(\pv_{t-1}))$ and
$\pv_t(\sigma(\pv_{t-1}))=\M(I_t, I_{t-1},I^\sigma_{t-1})$. Thus
we have the following tower of antiderivatives $$\M(I_t,
I_{t-1},I^\sigma_{t-1})\supset\M(I_{t-1},I^\sigma_{t-1})\supset\M.$$
We also know that $I_t$ consists of J-I-E antiderivatives of
$\M(I_{t-1},I^\sigma_{t-1})$. Now applying lemma \ref{the big
lemma} we obtain that $\sigma(\y_{tj})=\y_{tj}+c_{tj\sigma}$ for
every $\sigma\in\sg$. Also note that
$\sigma(\y_{tj})=\y_{tj}+c_{tj\sigma}$ implies
$$\sigma(\sum^{n_t}_{j=1}\alpha_{tj}\y_{tj})=\sum^{n_t}_{j=1}\alpha_{tj}\y_{tj}
+\sum^{n_t}_{j=1}\alpha_{tj}c_{tj\sigma}.$$ Thus
\begin{align*}&\sigma(\s)\in\pv\\
\implies&\sigma(\s)-\sum^{n_t}_{j=1}\alpha_{tj}\y_{tj}\in\pv\\
\implies&\sum^{n_t}_{j=1}\alpha_{tj}c_{tj\sigma}+\sigma(\sum^{t-1}_{i=1}\sum^{n_i}_{j=1}\alpha_{ij}\y_{ij})\in\pv\\
\implies&\sigma(\sum^{t-1}_{i=1}\sum^{n_i}_{j=1}\alpha_{ij}\y_{ij})\in\pv.
\end{align*}

Now we apply our induction hypothesis to the sum
$\sum^{t-1}_{i=1}\sum^{n_i}_{j=1}\alpha_{ij}\y_{ij}$ to prove our
theorem.

\begin{corollary}
The group of differential automorphisms of $\pv$ over $\f$ is a
subgroup of $(\C,+)^n$, where $n=tr.d(\pv|\f)$.
\end{corollary}

From theorem we observe that if $\sigma(\y_{ij})\in\pv$ then
$\sigma(\y_{ij})=\y_{ij}+c_{ij\sigma}$ for some
$c_{ij\sigma}\in\C$. Thus $\G(\pv|\f)$ is a subgroup of
$(\C,+)^n$. \epf

Now we will prove a generalization of the Ostrowski theorem for a
tower of extensions by J-I-E antiderivatives.

\bt \label{kol-ost gen}{\textsf Generalized Ostrowski Theorem}

Let $\pv_k\supset\K\supset\f$ be an intermediate differential
field  and let $T_i\subseteq I_i$ be subsets such that $T_i$ is a
set of antiderivatives of $\K(\cup^{i-1}_{j=1}T_j)$ for each
$1\leq i\leq k$. If $\cup^k_{j=1}T_j$ is algebraically dependent
over $\K$ then there is a nonzero $\s\in\K\cap\Lambda_k$.

\et

\bpf

Suppose that $\cup^k_{j=1}T_i$ is algebraically dependent over
$\K$. Then there is a $t$ such that $T_t$ is algebraically
dependent over $\K(\cup^{t-1}_{j=1}T_j)$. Then by The
Kolchin-Ostrowski theorem, there is a non zero
$\ft_t\in\K(\cup^{t-1}_{j=1}T_j)\cap\Lambda_t$. Let $\sg_{t-1}$ be
the group of all differential automorphisms of
$\K(\cup^{t-2}_{j=1}T_j)(T_{t-1})$ over $\K(\cup^{t-2}_{j=1}T_j)$.
Note that for every $\sigma\in\sg_{t-1}$,
$\sigma(y)\in\K(\cup^{t-2}_{j=1}T_j)(T_{t-1})\subseteq\pv_k$ for
every $y\in\K(\cup^{t-2}_{j=1}T_j)(T_{t-1})$ and that $\sg_{t-1}$
can be realized as a subgroup of $\G$. Thus we may apply theorem
\ref{group  action on ie towers} and obtain that
$\sigma(\ft_t)=\ft_t+\alpha_{\ft_t \sigma}$ for some
$\alpha_{\ft_t \sigma}\in\C$. This shows us that
$\ft_t\in\K(\cup^{t-2}_{j=1}T_j)(T_{t-1})$ is an antiderivative of
$\K(\cup^{t-2}_{j=1}T_j)$ and therefore the set $\{\ft_t\}\cup
T_{t-1}$ is algebraically dependent over
$\K(\cup^{t-2}_{j=1}T_j)$; observe that $\ft_t\notin
\cup^{t-1}_{j=1}T_j$. Again by the Kolchin-ostrowski theorem there
is a $\ft_{t-1}\in\Lambda_{t-1}$ and a constant $c_{\ft_t,t-1}$,
where $\ft_{t-1}$ or $c_{\ft_t,t-1}$ is nonzero such that
$$c_{\ft_t,t-1}\ft_t+\ft_{t-1}\in\K(\cup^{t-2}_{j=1}T_j)\cap\Lambda_{t-1}.$$

Now a repeated application of thereom \ref{group  action on ie
towers} and the Kolchin-Ostrowski theorem will prove the existence
of a nonzero $\s\in\K\cap\Lambda_t$.\epf

\bt \label{class jie thm} For every differential subfield $\K$ of
$\pv:=\pv_k$, the field generated by $\f$ and
$S_k:=\K\cap\Lambda_k$ equals the differential field $\K$. That is
$$\K=\f(S_k).$$
Moreover $\K$ itself is a tower of extensions by antiderivatives,
namely
$$\K=\f(S_k)\supset\f(S_{k-1})\supset\f(S_{k-2})\supset\cdots\supset\f(S_1)\supset\f,$$
where $S_i:=S_k\cap\Lambda_i$.

\et

\bpf

We will use an induction on $k$ to prove this theorem. $k=1$: Here
$\pv:=\pv_1$ is an extension by antiderivatives of $\f$ and
therefore from theorem \ref{struct theorem for antideriv} our
desired result follows immediately

$k\geq 2$: Assume that for any differential subfield of
$\pv_{k-1}$ our theorem is true. Let $S_i:=\K\cap\Lambda_i$ and
note that $S_i\supset S_{i-1}$ and the following containments
\begin{equation}\pv\supseteq\K\supseteq\f(S_k)\supseteq\f.
\end{equation}
We will first show that $\f(S)$ is a differential field. Applying
our induction hypothesis to the differential field $\langle
\f(S'_i)\rangle\subseteq\pv_{k-1}$, where $S'_i=\{\y'|\y\in S_i\}$
we obtain that $\langle \f(S'_i)\rangle=\f(T)$, where $T=\langle
\f(S'_i)\rangle\cap\Lambda_{i-1}$. Also note that $\langle
\f(S'_i)\rangle\subseteq\K$ and therefore
$$T=\langle \f(S'_i)\rangle\cap\Lambda_{i-1}\subseteq K\cap\Lambda_{i-1}\subseteq K\cap\Lambda_i=S_i.$$
Thus $\f(S_i)\supseteq\f(T)$ and since $S'_i\subset\f(T)$ and
$\f(T)$ is a differential field, $\f(S_i)$ is also a differential
field. Hence $\f(S_k)$ is a differential field and

$$\f(S_k)\supset\f(S_{k-1})\supset\f(S_{k-2})\supset\cdots\supset\f(S_1)\supset\f,$$
is a tower of extension by antiderivatives.

Let $\bar{S}_i\subset I_i$ be a transcendence base of $\pv_i$ over
the differential field $\pv_{i-1}(S_i)$.  Since $\pv_i$ is purely
transcendental over $\pv_{i-1}$ it is also purely transcendental
over $\pv_{i-1}(\bar{S}_i)$ too and therefore
$\pv_{i-1}(S_i,\bar{S}_i)=\pv_i$. We note that
$\f(S_1,\bar{S}_1)=\pv_1$,
$\f(S_2,\bar{S}_1,\bar{S}_2)=\pv_1(S_2,\bar{S}_2)=\pv_2$ and in
general we have $\f(S_t)(\cup^t_{i=1}\bar{S}_i)=\pv_t$. Since
$\K\supseteq\f(S_k)$ we have
\begin{equation}\label{full tower}\pv=\K(\cup^k_{i=1}\bar{S}_i)\supseteq\K(\cup^{k-1}_{i=1}\bar{S}_i)\supseteq\cdots\supseteq\K(\bar{S}_2,\bar{S}_1)
\supseteq\K(\bar{S}_1)\supseteq\K\supseteq\f(S_k)\supseteq\f\end{equation}

We know that $\cup^t_{i=1}\bar{S}_i$ is algebraically independent
over $\f(S_t)$. Since $S_t=\K\cap\Lambda_k$ we obtain from theorem
\ref{kol-ost gen} that $\cup^t_{i=1}\bar{S}_i$ is algebraically
independent over $\K$. Now from equation \ref{full tower} we
obtain
\begin{equation}\label{trd eqn 1}tr.d
(\pv|\f)=\sum^k_{i=1}|\bar{S}_i|+ tr.d (\K|\f(S_k))+ tr.d
(\f(S_k)|\f).\end{equation}

On the other hand we have
$$\pv=\f(S_k)(\cup^k_{i=1}\bar{S}_i)\supseteq\f(S_k)\supseteq\f
$$and thus \begin{equation}\label{trd eqn 2}tr.d
(\pv|\f)=\sum^k_{i=1}|\bar{S}_i|+ tr.d (\f(S_k)|\f)\end{equation}
From equation \ref{trd eqn 1} and \ref{trd eqn 2} we obtain tr.d
$(\K|\f(S_k))=0$, that is, $\K$ is algebraic over $\f(S_k)$. Now
from theorem \ref{no alg extn thm} we obtain $\K=\f(S_k)$. \epf

\subsection{Example}
Let $\C:=\mathbb{C}$ denote the complex numbers, $\C_\infty$ the
complete Picard-Vessiot closure of $\C$, $x\in\C_\infty$ be an
element whose derivative is $1$, $\tan^{-1}x\in\C_\infty$ be an
element such that
$$(\tan^{-1}x)'=\frac{1}{1+x^2}$$ and let
$\tan^{-1}(\tan^{-1}x)\in\C_\infty$ be an element such that
$$\left(\tan^{-1}(\tan^{-1}x)\right)'=\frac{1}{(1+(\tan^{-1}x)^2)(1+x^2)}.$$

We will use theorem \ref{class jie thm} to compute the
differential field $\C\langle \tan^{-1}(\tan^{-1}(x))\rangle$.
First we observe that $(\tan^{-1}(x))'=\frac{1}{1+x^2}$
$=\frac{1}{(x+i)(x-i)}$ and thus $\tan^{-1}x$ is an I-E(J-I-E)
antiderivative of $\C(x)$. We also observe that
$\tan^{-1}(\tan^{-1}(x))$ is an I-E(J-I-E) antiderivative of
$\C(x,\tan^{-1}(x))$(note that
$\left(\tan^{-1}(\tan^{-1}(x))\right)'=\frac{1}{(1+(\tan^{-1}x)^2)(1+x^2)})$.
Thus $x, \tan^{-1}(x), \tan^{-1}(\tan^{-1}(x))$ are algebraically
independent over $\C$. Also from theorem \ref{class jie thm} we
see that there should be a linear combination of the form
$c_1\tan^{-1}x+c_2x$, where $c_1$ is non zero (since the
$\frac{1}{(1+(\tan^{-1}x)^2)(1+x^2)}$ $\in\C\langle
\tan^{-1}(\tan^{-1}x)\rangle$ ). Thus by differentiating
$c_1\tan^{-1}x+ c_2x$, we see that $x\in\C\langle
\tan^{-1}(\tan^{-1}x)\rangle$ and therefore
$\tan^{-1}(x)\in\C\langle $ $\tan^{-1}(\tan^{-1}x)\rangle$ since
$c_1\tan^{-1}x+c_2x\in\C\langle \tan^{-1}(\tan^{-1}x)\rangle$.
Hence $$\C\langle \tan^{-1}(\tan^{-1}x)\rangle=\C(
\tan^{-1}(\tan^{-1}x),\tan^{-1}x, x).$$

We observe that
$$\left(\frac{1}{2i}\ln(x-i)-\frac{1}{2i}\ln(x+i)\right)'=\frac{1}{x^2+1}$$
and since $(\tan^{-1}x)'=\frac{1}{x^2+1}$ there is a $c\in\C$ such
that
$$\tan^{-1}x=\frac{1}{2i}\ln(x-i)-\frac{1}{2i}\ln(x+i)+c.$$
Also note that
\begin{align*}&\frac{1}{2i}\big(\ln(\frac{1}{2i}\ln(x-i)-\frac{1}{2i}\ln(x+i)+c-i)
-(\ln(\frac{1}{2i}\ln(x-i)-\frac{1}{2i}\ln(x+i)+c+i)\big)'\\
&=\frac{1}{x^2+1}
\left(\frac{1}{\frac{1}{2i}\ln(x-i)-\frac{1}{2i}\ln(x+i)+c-i}-\frac{1}{\frac{1}{2i}\ln(x-i)-\frac{1}{2i}\ln(x+i)+c+i}\right)\\
&=\frac{1}{(1+(\tan^{-1}x)^2)(x^2+1)}
\end{align*}
and since
$(\tan^{-1}(\tan^{-1}x))'=\dfrac{1}{(1+(\tan^{-1}x)^2)(x^2+1)},$
there is a constant $d\in\C$ such that
\begin{align*}\tan^{-1}(\tan^{-1}(x))&=\frac{1}{2i}\big(\ln(\frac{1}{2i}\ln(x-i)-\frac{1}{2i}\ln(x+i)+c-i)\\
&-(\ln(\frac{1}{2i}\ln(x-i)-\frac{1}{2i}\ln(x+i)+c+i)\big)+d.\end{align*}
Hence $$\tan^{-1}(\tan^{-1}x)\in\C(x,y_1,y_2,z),$$ where
\begin{align*}z:=&\frac{1}{2i}\Big(\ln(\frac{1}{2i}\ln(x-i)-\frac{1}{2i}\ln(x+i)+c-i)\\
-&\ln(\frac{1}{2i}\ln(x-i)-\frac{1}{2i}\ln(x+i)+c+i)\Big)\\
y_1:=&\ln(x-i)\\
y_2:=&\ln(x+i).
\end{align*} Clearly $$\C\langle \tan^{-1}(\tan^{-1}x)\rangle=\C(
z,\frac{1}{2i}(y_1-y_2), x).$$

\begin{remark}
The J-I-E extensions may have non-elementary functions. For
example; if $a_i\in\C$ are distinct constants for $i=1,\cdots,n$
then the elements $\y_i:=\int\frac{\ln(x)}{x-a_i}$ are J-I-E
antiderivatives of the differential field $\C(x,\ln(x))$ with
$\y'_i:=\frac{A_i}{C_iB_i}$ where $A_i:=\ln(x)$, $B_i:=1$ and
$C_i:=x-a_i$. These $\y_i$'s are non-elementary functions, see
\cite{elen}. From theorem \ref{alg indp of antideriv} we see that
these $\y_i$'s are algebraically independent over $\C(x,\ln(x))$
and from therorem\ref{class jie thm} we see that any differential
field $\K$, $\C(x,\ln(x),\y_i|1\leq i\leq
n)\supseteq\K\supseteq\C$ is of the form $\C(S)$, where
$S\subset$span$_\C\{x,\ln(x),\y_i|1\leq i\leq n\}$ is a finite
set. Moreover $\C(S)$ itself is a tower of (Picard-Vessiot)
extensions by antiderivatives.

\end{remark}


\section{Extensions by iterated logarithms}\label{the iter logs}

In this section we will provide an example of a J-I-E tower
namely, the extensions by iterated logarithms. Though many of the
results for iterated logarithms setting can be deduced from the
J-I-E tower setting from section \ref{the jie tower}, we  will
still prove those results here separately and this will help us in
writing an algorithm for computing the finitely differentially
generated subfields of the extensions by iterated logarithms.

{\bf{\small Iterated Logarithms}}

Let $\C$ be an algebraically closed-characteristic zero
differential field with a trivial derivation and let$\C_\infty$ be
the complete Picard-Vessiot closure of $\C$. Let
$\li[0,0]\in\C_\infty$ be an element such that $\li'[0,0]=1$. We
will often denote $\li[0,0]$ by $x$. Given
$\vec{c}=(c_1,\cdots,c_n)\in\C^n$ let $\li[\vec{c},n]\in\C_\infty$
be an element such that
\begin{equation}\label{defn of iter log}\li'[\vec{c},n]=\frac{\li'[\pi(\vec{c}),n-1]}{\li[\pi(\vec{c}),n-1]+\psi_n(\vec{c})},\end{equation}
where $\psi_n:\C^n\to\C$ is the map $\psi_n(c_1,\cdots,c_n)=c_n$
and
$\pi:\C^n\to\C^{n-1}$ is the map $$\left\{%
\begin{array}{ll}
   \pi(c_1,\cdots,c_n)=(c_1,\cdots,c_{n-1}) , & \hbox{when $n>1$;} \\
    \pi(c)=0, & \hbox{when $n=1$.} \\
\end{array}%
\right.$$

Whenever we write $\li[\vec{c},n]$, it is understood that
$\vec{c}\in\C^n$. We observe that for $\vec{c}=(c)\in\C$
\begin{align*}\li'[\vec{c},1]&=\frac{\li'[\pi(\vec{c}),0]}{\li[\pi(\vec{c}),0]+\psi_1(\vec{c})}\\
&=\frac{\li'[0,0]}{\li[0,0]+c}\\
&=\frac{1}{\li[0,0]+c}.
\end{align*}
Thus for $c\in\C$, the element $\li[\vec{c},1]$ can be seen as the
element $\ln(x+c)$. Similarly for
$\vec{c}=(c_1,\cdots,c_n)\in\C^n$, the element $\li[\vec{c},n]$
can be seen as the element
$\ln(\ln(\cdots(\ln(x+c_1)+c_2)\cdots+c_{n-1})+c_n)$.

For $1\leq k\leq n-1$, let $\pi^k:\C^n\to\C^{n-k}$ be the map
$\pi^k(c_1,\cdots,c_n)=(c_1,\cdots,c_{n-k})$ and let
$\pi^n:\C^n\to\C^0:=\{0\}$ be the zero map. For $1 \leq k\leq n$
let $\psi_k:\C^n\to\C$ be the map $\psi_k(c_1,\cdots,c_n)=c_k$.
Under these notations, we can rewrite equation \ref{defn of iter
log} as
\begin{equation}\label{iter log written using proj}
 \li'[\vec{c},n]=\Big(\prod^{n-1}_{i=1}\frac{1}{\li[\pi^{i+1}(\vec{c}),n-(i+1)]+\psi_{n-i}(\pi^i(\vec{c}))}\Big)
\frac{1}{\li[\pi(\vec{c}),n-1]+\psi_n(\vec{c})}.
\end{equation}
This above equation is obtained simply by clearing the derivative
that appears in the numerator of the RHS of the equation \ref{defn
of iter log}. Note that
\begin{equation}\label{diff iter
log}\li'[\pi(\vec{c}),n-1]=\prod^{n-1}_{i=1}\frac{1}{\li[\pi^{i+1}(\vec{c}),n-(i+1)]+\psi_{n-i}(\pi^i(\vec{c}))}.\end{equation}

\begin{definition}
When $n\in\N$ we will call $\li[\vec{c},n]$ an \textsl{$n^{th}$
level iterated logarithm} or simply an \textsl{iterated
logarithm}, without specifying its level.

We note that $\li[0,0]$, whose derivative equals 1, is not an
iterated logarithm under our definition. Hereafter we will call
$\li[0,0]$ as $x$.
\end{definition}
{\bf{\small More notations}}

Let $\Lambda_0:=\{x\}$,
$\Lambda_n:=\{\li[\vec{c},n]|\vec{c}\in\C^n\}$ and
$\Lambda_\infty=\cup^\infty_{i=0}\Lambda_i$ and let
$\Cd_0=\C(\Lambda_0)$, $\Cd_n:=\C(\cup^n_{i=0}\Lambda_i)$ and
$\Cd_\infty=\C(\Lambda_\infty)$. Note that $\Cd_0$, $\Cd_n$ and
$\Cd_\infty$ are differential fields(follows from equation
\ref{iter log written using proj}).

Let $\vec{c}\in\C^n$. We define
$\pi^k(\li[\vec{c},n]:=\li[\pi^k(\vec{c}),n-k]$ whenever $k\leq
n$. Note that $\pi^n(\li[\vec{c},n])=\li[0,0]=x$. When $k>n$ we
define $\pi^k(\li[\vec{c},n]):=x$ and $\pi^k(x):=x$ for any
$k\in\N$. Now we may also define $\pi^k(S)$ for a non empty set
$S\subset\Lambda_\infty$ as $\pi^k(S)=\{\pi^k(y)|y\in S\}$. Thus
\begin{equation} \label{pi property}\text{if}\ y\in\Lambda_n,\
\text{then}\ \pi(y)\in\Lambda_{n-1},\ \pi^2(y)\in\Lambda_{n-2},\
\cdots,\ \pi^n(y)=x\in\Lambda_0.\end{equation}
We also see that if $E\subset\Lambda_\infty$ is a finite set, then
there is an $n\in\N$ such that $\pi^n(E)=\{x\}$. Given a nonempty
set $E\subset\Lambda_\infty$ it is not necessary that $\C(E)$ is a
differential field. For example $\C(\li[\vec{0},1])$, that is the
field $\C(\ln(x))$ is not a differential field. whereas,
$\C(\ln(x),x)=\C(\li[\vec{0},1],\pi(\li[\vec{0},1])=x)$ is a
differential field. \big(note that $x\notin\C(\ln(x))$; in fact
$x$ and $\ln(x)$ are algebraically independent over $\C$. We will
later show that any collection of iterated logarithms is
algebraically independent over $\C(x)$.\big) More in general we
have the following propositions.

\bprop\label{diff field gen by a iter log} Let
$\li[\vec{c},n]\in\Lambda_\infty$ be an iterated logarithm. Then
$$\C(\li[\vec{c},n],\li[\pi(\vec{c}),n-1],\li[\pi^2(\vec{c}),n-2],\cdots,\li[\pi^n(\vec{c}),0]=x),$$
is a differential field \eprop

\bpf

We will use an induction on $n$ to prove our proposition.
\newline$n=1$. Note that $\li'[c,1]=\frac{1}{x+c}$ and $x'=1$. Therefore $\C(\li[c,1],x)$ is also a differential field.
We recall that if $\vec{v}\in\C^{n}$ then $\pi^n(\vec{v})=0$ and
therefore $\li[\pi^n(\vec{v}),n-n]=\li[0,0]=x$. Let us assume for
any $\vec{v}\in\C^{n}$ that
$\C(\li[\vec{v},n],\li[\pi(\vec{v}),n-1],$ $\cdots,x)$ is a
differential field and let $\vec{c}\in\C^{n+1}$. From our
induction hypothesis, we know that
$\f:=\C(\li[\pi(\vec{c}),n-1],\li[\pi^2(\vec{c}),n-2],\cdots,$
$x]$ is a differential field since $\pi(\vec{c})\in\C^{n-1}$. Thus
$$\li'[\vec{c},n]=\frac{\li'[\pi(\vec{c}),n-1]}{\li[\pi(\vec{c}),n-1]+\psi_n(\vec{c})}\in\f.$$
Hence
$\f(\li[\vec{c},n])=\C(\li[\vec{c},n],\li[\pi(\vec{c}),n-1],\cdots,
x)$ is a differential field. \epf

\bprop Let $E\subset\Lambda_\infty$ be a finite set of iterated
logarithms. Then $$\C(E,\pi(E),\pi^2(E),\cdots,x)$$ is a
differential field \eprop

\bpf If $E=\emptyset$ then $\C(E,\pi(E),\pi^2(E),\cdots,x)=\C(x)$
which is a differential field and we are done.  Let $E=\{y_j|1\leq
j\leq s\}$. We know from proposition \ref{diff field gen by a iter
log} that $\mathbf{K}_j:=\C(y_j,\pi(y_j),\cdots,\pi^{n_j}(y_j)=x)$
is a differential field and since $\C(E,\pi(E),\pi^2(E),\cdots,x)$
is a compositum of differential fields $\mathbf{K}_j$, we see that
$\C(E,\pi(E),\pi^2(E),\cdots,x)$ is also a differential field.\epf

\begin{definition}
For $E\subset\Lambda_\infty$, we will call the field $\C(E)$ an
\textsl{extension by iterated logarithms} if $E$ contains at least
one iterated logarithm, that is, if $E$ has an element from
$\Lambda_\infty$ other than $x$. And, we will call the
differential field $\C(E,$ $\pi(E),\pi^2(E),\cdots,x)$  as the
\textsl{Container Differential Field}[{\em CDF}] for the set $E$.
\end{definition}

\subsection{The Two Towers and a Structure Theorem for $\Cd_n$}

Let $E\subset\Lambda_\infty$ be a finite non empty set. Then there
is a minimal $n\in\N$ such that $\pi^n(E)=\{x\}$. Once this
minimal $n$ is chosen, it is clear that $E$ contains at least one
element from $\Lambda_n$ and no elements from $\Lambda_i$ for any
$i>n$. Hereafter we will use the symbol $\Ge$ to denote
$\cup^n_{i=0}\pi^i(E)$, where $n$ satisfies the above minimality
condition. Thus $\C(E,\pi(E),\pi^2(E),\cdots,x)$, the container
differential field of $E$ is the field $\C(\Ge)$. Note that
$\pi(\Ge)\subset\Ge$ and let $T_i:=\Lambda_i\cap \Ge$ for all
$1\leq i\leq n$. Then the $T_i$'s are disjoint and partitions
$\Ge$ in such a way that each $T_i$ contains iterated logarithms
only from level $i$. Clearly $E\subseteq\Ge$, and $\Ge$ may
contain more elements than $E$, but those elements that are in
$\Ge$ but not in $E$ has to come from $\cup^{n-1}_{i=0}\Lambda_i$.
Thus $T_n:=\Lambda_n\cap\Ge=\Lambda_n\cap E$. Also we observe that
$\C(\Ge)$ is a differential field and it contains $\C(E)$.
\begin{definition}
 We will call this partition $T_0, T_1,\cdots, T_n$ of $\Ge$ as the \textsl{levelled partition} of
$\Ge$.\end{definition} We observe that
\begin{align*}\pi(T_i)&=\pi(\Lambda_i\cap \Ge)\\
&\subseteq\pi(\Lambda_i)\cap\pi(\Ge)\\
&\subseteq\pi(\Lambda_i)\cap\Ge\\
&\subseteq T_{i-1}.
\end{align*}
Thus $\pi(T_i)\subseteq T_{i-1}$ for all $1\leq i\leq n$. We also
note that $T_0=\{x\}$ since $E$ is non empty. We will use this
partition of $\Ge$ to prove that the iterated logarithms are
algebraically independent over $\C(x)$ and this will be done in
subsection \ref{subsect Alg indp}.

Now we will construct a tower of Picard-Vessiot extensions by
antiderivatives(iterated logarithms) to reach $\C(\Ge)$ from $\C$
using this \textsl{levelled partition} of $\Ge$. \big(Note that
this tower is not imbeddable in the Picard-Vessiot closure of
$\C$.


The construction of this tower is obvious. Let ${\bf
K}_0:=\C(T_0)=\C(x)$ and let ${\bf K}_i:={\bf K}_{i-1}(T_i)$ for
all $i\in\N$. That is ${\bf K}_i=\C(\cup^i_{j=0}T_j)$ for $0\leq
i\leq n$. Clearly ${\bf K}_0$ is an extension by antiderivatives
of $\C$. Also, for $y\in T_i$, $\pi^j(y)\in\cup^{i-1}_{k=0}T_k$
for all $i,j\in\N$ and in fact, $\pi^j(y)=x$ for all $j\geq i$.
Now from equation \ref{iter log written using proj} we see that
$y'\in{\bf K}_{i-1}$ and thus ${\bf K}_i$ is also an extension by
antiderivatives of ${\bf K}_{i-1}$. Therefore we have a tower of
P-V extensions by antiderivatives namely
\begin{equation} \label{the level tower}\C(\Ge)={\bf K}_n
\supset{\bf K}_{n-1}\supset\cdots \supset{\bf K}_1\supset{\bf
K}_0\supset\C. \end{equation} We will call this the
\textsl{levelled partition tower} of $\C(\Ge)$.

There is another useful way of dividing the set
$\Ge=\cup^n_{i=0}\pi^i(E)$. Let
$$\pr=E\setminus\cup^n_{i=1}\pi^i(E).$$ We claim that $\cup^n_{i=0}
\pi^i(\pr)=\Ge$. Before we prove this claim, we note that
$\pr\cap\pi^i(\pr)=\emptyset$ for all $i$, $1\leq i\leq n$ and
this statement immediately follows from the definition of $\pr$.
We also note that $\pr\subseteq E$.

Now we will use an induction argument to show that
$\cup^n_{i=0}\pi^i(\pr)=\Ge$. First we observe that
$\Ge=\cup^n_{i=0}(\Lambda_i\cap\Ge)$. From the choice of $n$ it is
clear that $\Lambda_n\cap E\neq \emptyset$. From equation \ref{pi
property} we see that for every $y\in\Lambda_n\cap E$,
$y\notin\cup^{n}_{i=1}\pi^i(E)$. Thus $\Lambda_n\cap E\subseteq
\pr$ and therefore $\Lambda_n\cap E\subseteq
\pr\subseteq\cup^n_{i=0}\pi^i(\pr)$. Assume that there is a $k\leq
n$ such that for any $i$, $k\leq i\leq n$,
$\Lambda_{i}\cap\Ge\subseteq\cup^n_{i=0}\pi^i(\pr)$. We will show
that $\Lambda_{k-1}\cap\Ge\subseteq\cup^n_{i=0}\pi^i(\pr)$. Let
$y\in \Lambda_{k-1}\cap\Ge$. If $y\in\pr$, we are done. So, we
suppose that $y\notin\pr$. Then $y\in\cup^n_{i=1}\pi^i(E)$ and
therefore there is a $z\in E$ and a $j\in \N$ such that
$\pi^j(z)=y$. Clearly such a $z\in\cup^n_{i=k}\Lambda_i\cap E$.
That is, $z$ has to be a higher level iterated logarithm than $y$
is (see equation \ref{pi property}). Now from our induction
hypothesis we obtain $z\in\cup^n_{i=0}\pi^i(\pr)$ and since
$\cup^n_{i=0}\pi^i(\pr)$ is invariant under $\pi$, we obtain
$y\in\cup^n_{i=0}\pi^i(\pr)$. Thus $\cup^n_{i=0}\pi^i(\pr)=\Ge$.

\begin{definition}
We will call the set $\pr\subset E$ as the $\pi-$\textsl{base} of
$\Ge$.\end{definition}

We may also construct a tower of Picard-Vessiot extension by
antiderivatives(by iterated logarithms) to reach $\C(\Ge)$ by
defining $\pt_i:=\pt_{i-1}(\pi^{n-i}(\pr))$ for $1\leq i\leq n$,
where $\pt_0:=\C(x)$. Then $\pt_i=\C(\cup^i_{j=0}\pi^{n-j}(\pr))$
for $0\leq i\leq n$ and clearly, $\pt_i$ is a differential field.
Thus we see that
\begin{equation}\label{the pi tower} \C(\Ge)={\bf P}_n\supset{\bf
P}_{n-1}\supset\cdots \supset{\bf P}_1\supset{\bf P}_0\supset\C.
\end{equation}
We will call the above tower as the $\pi-$\textsl{tower} of
$\C(\Ge)$.

We observe that $\pr\subset\cup^n_{i=0}\Lambda_i$ and therefore
$\pi(\pr)\subset\cup^{n-1}_{i=0}\Lambda_i$,
$\pi^2(\pr)\subset\cup^{n-2}_{i=0}\Lambda_i$ and in general
$\pi^j(\pr)\subset\cup^{n-j}_{i=0}\Lambda_i$. Thus
$\pi^{n-j}(\pr)\subset\cup^j_{i=0}\Lambda_i$ and from this fact we
also obtain
$\cup^m_{j=0}\pi^{n-j}(\pr)\subset\cup^m_{i=0}\Lambda_i$ for any
$m, 0\leq m\leq n$. Since $\cup^n_{i=0}\pi^i(\pr)=\Ge$,
\begin{align*}\cup^m_{i=0}\pi^{n-j}(\pr)
&\subseteq\cup^m_{i=0}\Lambda_i\cap\Ge\\
&=\cup^m_{i=0}T_i
\end{align*}
and thus $\cup^m_{i=0}\pi^{n-j}(\pr)\subseteq\cup^m_{i=0}T_i$.
This shows that $\pt_m\subseteq{\bf K}_m$ for every $0\leq m\leq
n$. Nonetheless the inequality could be strict and we will now
provide an example for the same.

Let $\C:=\mathbb{C}$ and let $E=\{\ln(\ln(x+e)+5), \ln(\ln(x)),
\ln(x), \ln(x+1)\}$. In our notation, the set
$E=\{\li[\vec{v}_1,2],\li[\vec{v}_2,2], \li[\vec{v}_3,1],
\li[\vec{v}_4,1]\}$, where $\vec{v}_1=(\exp,5)$,
$\vec{v}_2=(0,0)$, $\vec{v}_3=(0)$ and $\vec{v}_4=(1)$. Then we
immediately see that $\pi(\ln(\ln(x+\exp)+5))=\ln(x+e)$,
$\pi^2(\ln(\ln(x+\exp)+5))=\pi(\ln(x+e))$$=x$,
$\pi(\ln(\ln(x)))=\ln(x)$, $\pi^2(\ln(\ln(x)))=x$,
$\pi(\ln(\ln(x+1)))=\ln(x+1)$, $\pi(\ln(x+1))=x$ and $\pi(x)=x$.
Thus the set $\Ge=\{\ln(\ln(x+e)+5), \ln(\ln(x)), \ln(x),
\ln(x+1), \ln(x+e),x\}$.

Let us obtain the \textsl{levelled partition} of $\Ge$. The set
$T_0=\Ge\cap\Lambda_0=\{x\}$, $T_1=\Lambda_1\cap\Ge$
$=\{\ln(x),\ln(x+1),\ln(x+e)\}$ and the set $T_2=\Ge\cap\Lambda_2$
$=\{\ln(\ln(x)),\ln(\ln(x+e)+5)\}$. Therefore the \textsl{levelled
partition tower} would be
$$\C(\Ge)\supset\C(\ln(x),\ln(x+1),\ln(x+e),x)\supset
\C(x)\supset\C.$$

Note that the $\pi-$base $\pr$ of $\Ge$ is given by
$\pr=E\setminus\cup^2_{i=1}\pi^i(E)$. Since $\cup^2_{i=1}\pi^i(E)$
$=\{\ln(x+e),\ln(x),x\}$ we see that $\pr=\{\ln(\ln(x+e)+5),
\ln(\ln(x)), \ln(x+1)\}$. Thus the $\pi-$\textsl{partition tower}
of $\C(\Ge)$ is
$$\C(\Ge)\supset\C(\ln(x),\ln(x+e),x)\supset
\C(x)\supset\C.$$

Therefore, if we assume that the iterated logarithms are
algebraically independent over $\C(x)$ then
$\ln(x+1)\notin\C(\ln(x),\ln(x+e),x)$ and thus the two towers are
distinct.

 {\bf Structure theorem for $\Cd_n$:} Here we will assume that
the iterated logarithms are algebraically independent over
$\C(x)$. That is, the set $\Lambda_\infty$ is algebraically
independent over $\C$. A proof for this fact is provided in
subsection \ref{subsect Alg indp}, theorem \ref{alg indp of trans
log over C}. Thus $\Cd_\infty$ is the field of fractions of the
polynomial ring $\C[\Lambda_\infty]$. For $y\in\Lambda_\infty$ let
$\frac{\partial }{\partial y}$ denote the standard partial
derivation on the polynomial ring $\C[\Lambda_\infty]$.

Let $u\in\Cd_n\setminus\Cd_{n-1}$. Then there is a finite non
empty set $S\subset\cup^n_{i=0}\Lambda_i$ such that
$u=\frac{P}{Q}$, $P, Q\in\C[S]$ and $(P,Q)=1$(that is the G.C.D of
$P$ and $Q$ in the polynomial ring $\C[S]$ is 1 ). It is
conceivable that some of the elements of $S$ may not be necessary
to express $u$. So, we define a set $E_u$ as
\begin{equation}E_u:=\Big\{y\in S\Big|\dfrac{\partial P}{\partial
y}\neq 0\ or\ \dfrac{\partial Q}{\partial y}\neq
0\Big\}.\end{equation}

\begin{definition}
The set $E_u$ is called the set of all \textsl{essential elements}
of $u$
\end{definition}
We observe that $u\in\C(E_u)$ and that if $u\in\C(S)$ for some set
$S\subset\Lambda_\infty$ then $E_u\subset S$. Sometimes we drop
the suffix $u$ and simply write $E$ instead of $E_u$. Since
$\C[\Lambda_\infty]$ is a polynomial ring (over a field), the set
$E_u$ is unique. The following theorem proves the uniqueness of
$E_u$.

\bt({\small Uniqueness of $E_u$})\label{unique of E} Let
$u\in\Cd_\infty$ and let $E_u$ be a set of \textsl{essential
elements} of $u$. Then $u\in\C(S)$ for some
$S\subset\Lambda_\infty$ only if $E_u\subseteq S$ and thus the set
$E_u$ is unique for a given $u$.\et

\bpf Let $S\subset\Lambda_\infty$ and let $u\in\C(S)$. Then
\begin{equation}u=\frac{P}{Q}
=\frac{A}{B},\end{equation} for some $A,B\in\C[S]$ and
$P,Q\in\C[E_u]$, where $(P,Q)=1$. Since $(P,Q)=1$, from the above
equation it is clear that $P$ divides $A$ and $Q$ divides $B$ in
the polynomial ring $\C[S\cup E_u]$. Thus there are
$R,T\in\C[S\cup E_u]$ such that $PR=A$ and $QT=B$. Note that if
$y\in E_u$ then $\frac{\partial P}{\partial y}\neq 0\ or\
\frac{\partial Q}{\partial y}\neq 0$. Suppose that there is a
$y\in E_u$ such that $\frac{\partial P}{\partial y}\neq 0$.
Consider the equation $PR=A$. Then $\deg_{y}(P)\geq 1$ and note
that $PR=A$ implies $\deg_{y}(P)+\deg_{y}(R)=\deg_{y}(A)$. Thus
$\deg_{y}(A)\geq 1$. Hence $y\in S$. Similarly if $\frac{\partial
Q}{\partial y}\neq 0$ and $\frac{\partial P}{\partial y}=0$, we
may use the equation $QT=B$ to show that $y\in S$ and thus
$E_u\subset S$. \epf

The following corollary is a direct consequence of the above
theorem. \bco \label{cor alg indp of trans log over C} Let
$S\subset\Lambda_\infty$ be any nonempty set and for $1\leq j\leq
s$ let $y_j\in\Lambda_\infty$ be distinct. Then for any constants
$a_j\in\C^*$ such that $\sum^s_{j=1}a_jy_j\in\C(S)$, the element
$y_j\in S$ for each $j$. \eco

\bpf Suppose that there are $a_j\in\C^*$ and such that
$\sum^s_{j=1}a_jy_j\in\C(S)$. Since $a_j\in\C^*$, the essential
elements of $\sum^s_{j=1}a_jy_j$ is the set $E:=\{y_j|1\leq j\leq
s\}$. Now from theorem \ref{unique of E} we obtain $E\subset S$.
\epf

Now we will state the structure theorem for singly generated
differential subfields of $\Cd_n$.

\bt Let $u\in\Cd_n\setminus\Cd_{n-1}$, $E$ the \textsl{essential
elements} of $u$ and $\C(\Ge)$ the container differential field of
$E$. Let $\pr\subseteq E$ be the $\pi-$ base of $\Ge$. Then the
differential field
$$\C\langle u\rangle=\C(\ls,\pi(\pr),\pi^2(\pr),\cdots,x),$$ where
$\ls$ is a finite nonempty subset of $span_\C \pr$. Moreover for
every $y\in\pr$, $\ls$ contains at least one linear combination in
which $y$ appears nontrivially. 
\et The above structure theorem is proved in subsection
\ref{subsection structure theorem}. There we will also generalize
this theorem to finitely generated differential subfields of
$\Cd_n$ and  give an algorithm to find the set $\ls$ and $\pr$
that appears in the above structure theorem.

\begin{remark}
Given a $u \in\Lambda_\infty$, there is a set finite
$E\subset\cup^n_{i=0}\Lambda_i$ and we may also choose a minimal
$n$ such that the above inclusion holds. Then $\C\langle u\rangle$
becomes a subfield of the container differential field $\C(\Ge)$
of $E$. The field $\C(\Ge)$ is an \textsl{elementary extension} of
$\C$. The above stated theorem(and its generalized version) shows
that every differential subfield of $\C(\Ge)$, more in general, a
finitely differentially generated subfields of $\Cd_\infty$ has to
be a \textsl{generalized elementary extension} of a special form.
For a definition of elementary and generalized elementary
extension and results related to our theorem in a more general
context, one may refer to the following papers \cite{M.Sing},
\cite{M.Sing 2} and \cite{R.Risch}.

\end{remark}

\subsection{Algebraic Independence of Iterated logarithms}\label{subsect Alg indp}

Here we will show that the set $\Lambda_\infty$ is algebraically
independent over $\C$. For $i=1,2,\cdots,n$ let $c_i\in\C$ be
distinct constants. By choosing $C_i:=x+c_i$, $A_i=B_i=1$, we see
that $\C(x,\li[\vec{c}_1,1],\cdots,\li[\vec{c}_n,1])$, where
$\vec{c}_i:=(c_i)$ is an extension by J-I-E antiderivatives of
$\C(x)$ and thus $\li[\vec{c}_1,1],\cdots,\li[\vec{c}_n,1]$ are
algebraically independent over $\C(x)$. Assume that every finite
subset of $\Lambda_{t-1}$, $t\geq 2$ consists of J-I-E
antiderivatives of $\C(\cup^{t-2}_{j=0}\Lambda_j)$. For
$i=1,2,\cdots,n$ let $\vec{c}_i:=(c_{1i},c_{2i},\cdots,c_{ti})$
$\in\C^t\setminus\{0\}$ be distinct vectors. Note that
\begin{equation}\label{iter log written using proj}
 \li'[\vec{c}_i,t]=\Big(\prod^{n-1}_{j=1}\frac{1}{\li[\pi^{j+1}(\vec{c}_i),n-(j+1)]+\psi_{n-j}(\pi^j(\vec{c}_i))}\Big)
\frac{1}{\li[\pi(\vec{c}_i),t-1]+\psi_t(\vec{c}_i)}
\end{equation}
and therefore choosing $A_i=1$,
$B_j:=\li[\pi^{j+1}(\vec{c}_i),t-(j+1)]+\psi_{t-j}(\pi^j(\vec{c}_i))$
and $C_j:=\li[\pi(\vec{c}_i),t-1]+\psi_t(\vec{c}_i)$ we see that
$\C(\cup^{t-1}_{j=0}\Lambda_j,$
$\li[\vec{c}_1,t],\cdots,\li[\vec{c}_n,t] )$ is an extension by
J-I-E antiderivatives of $\C(\cup^{t-1}_{j=0}\Lambda_j)$ and thus
$\Lambda_t$ is algebraically independent over
$\C(\cup^{t-1}_{j=0}\Lambda_j)$. Now we will give a proof for the
algebraic independence of the iterated logarithms without
appealing to results from section\ref{the jie tower}.

\begin{lemma}\label{ gen trans log alg dep}

Let $S_{n-1}\subset\Lambda_{n-1}$ be a finite set of
antiderivatives of  a differential field $\f$ and let
$S_n\subset\Lambda_n$ be such that $\pi(S_n)\subseteq S_{n-1}$.
Suppose that $S_{n-1}$ is algebraically independent over $\f$.
Then $S_n$ is algebraically independent over $\f(S_{n-1})$.

\end{lemma}

\bpf

Note that $\f(S_{n-1})$ is a differential field and since
$\pi(S_n)\subseteq S_{n-1}$, from equations \ref{defn of iter log}
and \ref{diff iter log} it is clear that $\f(S_{n-1})(S_n)$ is
also a differential field. Let $S_n=\{\li[\vec{c}_i,n]|1\leq i\leq
s\}$, $\vec{c}_i=(c_{1i},c_{2i},\cdots,c_{ni})$ and Suppose that
$S_n$ is algebraically dependent over $\f(S_{n-1})$. Then by
theorem \ref{Kol-Ost} there are constants $\alpha(\vec{c}_i)\in\C$
not all zero such that
$\sum^s_{i=1}\alpha(\vec{c}_i)\li[\vec{c}_i,n]\in\f(S_{n-1})$. We
may assume that $\alpha(\vec{c}_1)\neq 0$ and rewrite the sum as
$X+\sum^t_{j=1}\alpha(\vec{b}_j)\li[\vec{b}_j,n]$ where
$\{\vec{b}_j\}\subseteq\{\vec{c}_i\}$ is the set of all vectors
such that $\pi(\vec{b}_j)=\pi(\vec{c}_1)$
$=(c_{11},c_{21},\cdots,c_{n-1 1})$ and
$X=\sum^s_{i=1}\alpha(\vec{c}_i)\li[\vec{c}_i,n]-\sum^t_{j=1}\alpha(\vec{b}_j)\li[\vec{b}_j,n]$.
We may order the set $\{\vec{b}_j\}$ so that
$\vec{b}_1=\vec{c}_1$. Let
$\mathbf{K}:=\f(S_{n-1}\setminus\{\li[\pi(\vec{c_1}),n-1]\})$ and
let $X+\sum^t_{j=1}\alpha(\vec{b}_j)\li[\vec{b}_j,n]=\frac{P}{Q},$
where $P,Q \in\mathbf{K}[\li[\pi(\vec{c}_1),n-1]]$, $(P,Q)=1$ and
$Q$ a monic polynomial. Then
$$X'+\sum^t_{j=1}\frac{\alpha(\vec{b}_j)\li'[\pi(\vec{c}_1),n-1]}{(\li[\pi(\vec{c}_1),n-1]+c_{j
n})}=\frac{QP'-PQ'}{Q^2}.$$ Let $f:=\li'[\pi(\vec{c}_1),n-1]$ and
let
$\frac{F}{G}=\sum^t_{j=1}\frac{\alpha(\vec{b}_j)}{\li[\pi(\vec{c}_1),n]+c_{j
n}}$, where $F$ and $G$ are obtained by clearing the denominator
of the sum
$\sum^t_{j=1}\frac{\alpha(\vec{b}_j)}{\li[\pi(\vec{c}_1),n]+c_{j
n}}$. Note that $(F,G)=1$. Now we have
\begin{equation}\label{lo poly}Q^2(GX'+fF)=G(QP'-PQ').\end{equation}
From the definition of $X$, it is clear that
$X=\sum^t_{j=1}\alpha(\vec{a}_j)\li[\vec{a}_j,n]$ where
$\{\vec{a}_j\}\subset\{\vec{c}_j\}$ is the set of all vectors such
that $\pi(\vec{a}_j)\neq\pi(\vec{c}_1)$. Therefore
$X'\in\mathbf{K}$. Thus equation \ref{lo poly} is a polynomial in
$\li[\pi(\vec{c_1}),n-1]$ over the field $\mathbf{K}$.  Let
$y:=\li[\pi(\vec{c}_1),n-1]+c_{1n}$. Since $y$ divides $G$ and
$(F,G)=1$, $y$ does not divide $F$. Thus $y$ does not divide
$GX'+fF$ and therefore from \ref{lo poly} $y$ divides $Q^2$. Hence
$y$ divides $Q$. Let $l\in\N$ be the greatest positive integer
such that $y^l$ divides $Q$. Then $y^{2l}$ divides $Q^2$ and
therefore $y^{l+1}$ divides $Q^2$, which implies $y^{l+1}$ divides
$G(QP'-PQ')$. Since $y$ divides $G$ and $y^2$ does not divide $G$,
$y^l$ divides $QP'-PQ'$. But $y^l$ divides $Q$ and therefore $y^l$
divides $PQ'$. Since $(P,Q)=1$, we see that $y^l$ divides $Q'$.
Write $Q=y^lH$ and consider $Q'=ly^{l-1}y'H+y^lH'$. Note that
$y^l$ divides $Q'$ implies $y^l$ divides $ly^{l-1}y'H$ and since
$y'\in\mathbf{K}$, $y$ divides $H$. Thus $y^{l+1}$ divides $Q$,
contradicting the maximality of $l$.\epf

\bt \label{alg indp of trans log over C} Let
$E\subset\Lambda_\infty$ be a nonempty finite set. Then $E$ is
algebraically independent over $\C$.\et

\bpf As usual, let $\Ge:=\cup^n_{i=0}\pi^i(E)$ where $n$ is the
least positive integer such that $E\subset\cup^n_{i=0}\Lambda_i$
and let $\{T_i|0\leq i\leq n\}$ be the \textsl{levelled partition}
of $\Ge$. As we noted earlier $\pi(T_i)\subseteq T_{i-1}$,
$T_n\neq \emptyset$ and $\pi^n(T_n)=\{x\}=T_0$. Clearly, $T_0$ is
algebraically independent over $\C$ (see theorem \ref{anti elem})
and since $\pi(T_1)\subset T_0$, from lemma \ref{ gen trans log
alg dep} we get $T_1$ is algebraically independent over $\C(T_0)$.
Since $\pi(T_i)\subset T_{i-1}$, a repeated application of lemma
\ref{ gen trans log alg dep} will show us that
$\Ge=\cup^n_{j=0}T_j$ is algebraically independent over $\C$.
Since $E\subset \Ge$, $E$ is also algebraically independent over
$\C$. \epf


\subsection{Normality of $\Cd_n$ and Some Consequences:}
Let $\C_\infty$ be the complete Picard-Vessiot closure of $\C$ and
let $\Phi\in\G(\C_\infty|\C)$. Let $(v_i)_{i\in\N}$ be a sequence
in $\C$ and let $\vec{v}_n:=(v_1,\cdots,v_n)$ for all $n\in\N$
(the vector $\vec{v}_1=(v_1)$). Thus in our notation
$\pi(\vec{v}_n)=\vec{v}_{n-1}$. We observe that
$\Phi(x)=x+\alpha_{\Phi}$ for some $\alpha_{\Phi}\in\C$. Since
$\li'[\vec{v}_1,1]=\frac{1}{x+v_1}$ we see that
$\Phi\big(\li[\vec{v}_1,1]\big)'=\frac{1}{\Phi(x)+v_1}=\frac{1}{x+\alpha_{\Phi}+v_1}$
$=\li'[\Phi(\vec{v}_1),1]$, where
$\Phi(\vec{v}_1):=(v_1)+(\alpha_{\Phi})$. Since any two
antiderivatives differ by a constant,
$\Phi\big(\li[v_1,1]\big)=\li[\Phi(\vec{v}_1),1]+\alpha_{\Phi(\vec{v}_1)}$,
for some $\alpha_{\Phi(\vec{v}_1)}\in\C$. Assume that
$\Phi\big(\li[\vec{v}_{n-1},n-1]\big)=\li[\Phi(\vec{v}_{n-1})
,n-1]+\alpha_{\Phi(\vec{v}_{n-1})}$ where
$\Phi(\vec{v}_{n-1})=(v_1+\alpha_\Phi,v_2+\alpha_{\Phi(\vec{v}_1)},\cdots,v_{n-1}+\alpha_{\Phi(\vec{v}_{n-2})})$
and $\alpha_{\Phi(\vec{v}_{n-1})}\in\C$. Since
$$\li'[\vec{v}_n,n]=\frac{\li'[\vec{v}_{n-1},n-1]}{\li[\vec{v}_{n-1},n-1]+v_n},$$
we see that
\begin{align*}
\Phi\big(\li[\vec{v}_n,n]\big)'&=\frac{\li'[\Phi(\vec{v}_{n-1}),n-1]}
{\li[\Phi(\vec{v}_{n-1}),n-1]+v_n+\alpha_{\Phi(\vec{v}_{n-1})}}\\
&=\li'[\Phi(\vec{v}_n),n]
\end{align*}
where
$\Phi(\vec{v}_n)=(v_1+\alpha_\Phi,v_2+\alpha_{\Phi(\vec{v}_1)},\cdots,v_n+\alpha_{\Phi(\vec{v}_{n-1})})$.
Since any two antiderivatives differ by a constant, we obtain
\begin{equation}\label{action on iter log set}\Phi\big(\li[\vec{v}_n,n]\big)=\li[\Phi(\vec{v}_n),n]+\alpha_{\Phi(\vec{v}_n)}\end{equation}
for some $\alpha_{\Phi(\vec{v}_n)}\in\C$.

From equation\ref{action on iter log set}, we see that for every
$\Phi\in\G(\C_\infty|\C)$,
\begin{equation}\label{aut trans on lambda}\Phi(\Lambda_i)\subseteq\Lambda_i+\C\end{equation} for all $i\in\N$. Thus
$\Cd_n$ is a normal differential subfield of $\C_\infty$.

\begin{remark} \label{the normality rem}
Let $\Phi\in\G(\Cd_\infty|\C)$ and for $n\in\N\cup\{0\}$ let
$$\Phi\big(\li[\vec{v}_n,n]\big)=\li[\vec{v}_n,n]+\alpha_{\Phi(\vec{v}_n)},$$
with $\alpha_{\Phi(\vec{v}_n)}\in\C^*$. Then from the above
discussion, we see that for any $m<n$
$$\Phi\big(\li[\vec{v}_m,m]\big)=\li[\vec{v}_m,m].$$

For any $m>n$ and $k\in\N$
$$\Phi^k\big(\li[\vec{v}_m,m]\big)=\li[\Phi^k(\vec{v}_m),m]+\alpha_{\Phi^k(\vec{v}_m)},$$
where
$\Phi^k(\vec{v}_m)=(v_1,\cdots,v_n,v_{n+1}+k\alpha^{\Phi}_{\vec{v}_n},\cdots,v_{m}+k\alpha^{\Phi}_{\vec{v}_{m-1}})$.
Since $\alpha^{\Phi}_{\vec{v}_n}\neq 0$, $\Phi^i(\vec{v}_m)\neq
\Phi^j(\vec{v}_m)$ when $i\neq j$. Thus
$\li[\Phi^i(\vec{v}_m),m]\neq\li[\Phi^j(\vec{v}_m),m]$ for any
$i\neq j$ and for any $m>n$. Hence the set
$\{\li[\vec{v}_m,m],\li[\Phi^j(\vec{v}_m),m]|i\in\N\}$ is
algebraically independent over $\C$ for any $m>n$(follows from
theorem \ref{alg indp of trans log over C}).
\end{remark}

Now we will prove a theorem which will help us to prove the
structure theorem for the differential subfields of $\Cd_n$.

\bt \label{trans log in PV}Let $\f$ be a differential field
finitely generated over its constants $\C$, $\pv$ be a
Picard-Vessiot extension of $\f$, and let
$\f\subset\pv\subset\mathfrak{L}_\infty$. If
$\sum^s_{j=1}a_jy_j\in\pv$ for some $a_j\in\C\setminus\{0\}$,
$y_j\in \cup^\infty_{i=0}\Lambda_i$ and $s\in\N$ then
$\pi^i(y_j)\in\f$ for all $i\in\N$ and thus $y'_j\in\f$.

\et

\bpf Let there be $y_j\in\cup^\infty_{i=0}\Lambda_i$ and
$a_j\in\C^*$ such that $\sum^s_{j=1}a_jy_j\in\pv$. Note that $\pv$
is finitely generated over $\f$ and $\f$ is finitely generated
over $\C$ and thus $\pv$ is finitely generated over $\C$. Let
$u_1,\cdots,u_t\in\pv$ such that $\C(u_1,\cdots,u_t)=\pv$,
$E_{u_i}$ be the set of \textsl{essential elements} of $u_i$, and
let $S:=\cup^t_{i=1}E_{u_i}\cup\{y_j|1\leq j\leq s\}$. From the
definition of $S$ it is quite clear that we have the following
containments
\begin{equation}\label{cruc containments}\C(S)\supseteq\pv(y_1,\cdots,y_s)\supseteq\pv\supseteq\f\supseteq\C.\end{equation}

Since $\Cd_n$ and $\pv$ are normal differential subfields of the
complete Picard-Vessiot closure $\f_\infty$ of $\f$, every
automorphism $\phi\in\G(\pv|\f)$ extends to an automorphism
$\Phi\in\G(\mathfrak{L}_n|\f)$ and every automorphism
$\Phi\in\G(\Cd_n|\f)$ restricts to an automorphism
$\phi\in\G(\pv|\f)$.

Let $\Phi\in\G(\Cd_n|\f)$. Since $\pv$ is a normal differential
subfield of $\Cd_n|\f$, $\Phi(\pv)\subseteq\pv$ and therefore
\begin{equation}\sum^s_{j=1}a_j\Phi^k(y_j)\in\mathcal{\pv}.\end{equation}

Let $y_j=\li[\vec{v}_{jm_j},m_j]$, where
$\vec{v}_{jm_j}=(v_{j1},\cdots,v_{jm_j})$. Then
$\Phi^k(y_j)=\li[\Phi^k(\vec{v}_{jm_j}),m_j]+\alpha_{\Phi^k(\vec{v}_{jm_j})}$,
where $\alpha_{\Phi^k(\vec{v}_{jm_j})}\in\C$. Therefore
$$\sum^s_{j=1}a_j\li[\Phi^k(\vec{v}_{jm_j}),m_j]+\sum^s_{j=1}a_j\alpha_{\Phi^k(\vec{v}_{jm_j})}\in\mathcal{\pv}$$
and thus
$$\sum^s_{j=1}a_j\li[\Phi^k(\vec{v}_{jm_j}),m_j]\in\mathcal{\pv}\subseteq\C(S).$$
Now from corollary \ref{cor alg indp of trans log over C} we see
that
$$\li[\Phi^k(\vec{v}_{jm_j}),m_j]\in S$$ for every $j,k\in\N$. For
a fixed $j$, consider the set $T:=\{\li[\vec{v}_{jm_j},m_j],$
$\li[\Phi^k(\vec{v}_{jm_j}),m_j]|k\in\N\}$. From the action of
$\Phi$ on $\vec{v}_{jm_j}$, it is clear that if
$\Phi(\vec{v}_{jm_j})\neq\vec{v}_{jm_j}$ then $T$ is infinite. But
$T$ cannot be infinite because it sits inside the finite set $S$.
Hence $\Phi(\vec{v}_{jm_j})=\vec{v}_{jm_j}$ and therefore
\begin{align*}\Phi(\li[\vec{v}_{jm_j},m_j])&=\li[\Phi(\vec{v}_{jm_j}),m_j]+\alpha_{\Phi(\vec{v}_{jm_j})}\\
&=\li[\vec{v}_{jm_j},m_j]+\alpha_{\vec{v}_{jm_j}}.
\end{align*}
Now from the remark \ref{the normality rem} it follows that
$\Phi(\pi^i(y_j))=\pi^i(y_j)$ for all $i\in\N$. This shows that
$\pi^i(y_j)\in\Cd^{\G(\Cd_n|\f)}_n=\f$.\epf

\subsection{Differential Subfields of $\Lambda_\infty$}\label{subsection structure theorem}

In this section we will classify the finitely generated
differential subfields of $\Cd_n$. First we will point out an
interesting property that every differential subfield $\f\neq \C$
of $\Lambda_n$ possesses, which is that $x\in\f$ and this result
is a consequence of the structure theorem.

\bprop\label{trans logs extn has no alg extn}

Let $u\in \Cd_n\setminus\Cd_{n-1}, n\in\N$, $E$ be the set of
essential elements of $u$, $\Ge:=\cup^n_{j=0}\pi^i(E)$ and let
$\{T_i|0\leq i\leq n\}$ be the \textsl{levelled partition} of
$\Ge$. Then $u$ is not algebraic over $\C(\cup^i_{j=0}T_j)$ for
any $0\leq i\leq n-1$. \eprop

\bpf

Let $u=\frac{P}{Q}$, $P,Q\in\C(\Ge\setminus\{y\})[y]$, where $y\in
T_n$. The \textsl{levelled partition} of $\Ge$ is constructed in
such a way that $T_n\neq\emptyset$ and $T_n\subseteq E$. Since $E$
consists of essential elements of $u$ and $y\in E$,
$u\notin\f:=\C(\Ge\setminus\{y\})$. Let
$\mathbf{K}_i:=\C(\cup^i_{j=0}S_j)$ for each $i, 1\leq i\leq n-1$.
Then $\mathbf{K}_i\subset\f$. Since $y'\in\f$ and $y\notin\f$,
$\pv=\f(y)$ is a Picard-Vessiot extension of $\f$ with a
differential Galois group $\G:=(\C,+)$. Note that $\G$ has no non
trivial algebraic subgroups(in particular no nontrivial finite
subgroups). Since $\f\langle u\rangle\varsupsetneq\f$, $\f\langle
u\rangle=\pv$, which implies $u$ is not algebraic over $\f$. Thus
$u$ is not algebraic over $\mathbf{K}_i$ for any $0\leq i\leq
n-1$. \epf

Thus we have just shown that if $$\C(\Ge)={\bf K}_n\supset{\bf
K}_{n-1}\supset\cdots \supset{\bf K}_1\supset{\bf K}_0\supset\C.$$
is the \textsl{levelled partition} tower of $\Ge$, where
$\Ge:=\cup^n_{j=0}\pi^i(E)$ and $E$ is the set of essential
elements of an element $u\in\Cd_n\setminus\Cd_{n-1}$ then $u$ is
not algebraic over ${\bf K}_i$ for any $0\leq i\leq n-1$.

Note that if $u\in\C(x)$ then $\C\langle u\rangle=\C$ or $\C(x)$
depending whether $u$ is a constant or not. Thus if $\f$ is a
differential subfield(need not be finitely generated) of $\C(x)$
then $\f=\C(x)$ or $\C$ depending whether $\f$ contains a
nonconstant or not. Thus it is enough to state the structure
theorem only for elements in $u\in\Cd_n\setminus\Cd_{n-1}$.

\bt\label{singly gen diff subfields of log closure}  Let
$u\in\Cd_n\setminus\Cd_{n-1}$, $E$ the \textsl{essential elements}
of $u$ and $\C(\Ge)$ the container differential field of $E$. Let
$\pr\subseteq E$ be the $\pi-$ base of $\Ge$. Then the
differential field
$$\C\langle u\rangle=\C(\ls,\pi(\pr),\pi^2(\pr),\cdots,x),$$ where
$\ls$ is a finite nonempty subset of $span_\C \pr$. Moreover, for
every $y\in\pr$, $\ls$ contains at least one linear combination in
which $y$ appears nontrivially.\et

\bpf For $i\geq 1$ let $\pt_{n-i}$ denote the differential field
$\C(\pi^i(\pr),\pi^{i+1}(\pr),$ $\cdots,x)$ and let
$\pt_{n-i}\langle u\rangle$ be the differential field generated by
$\pt_{n-i}$ and $u$. Note that
$\C(\Ge)=\pt_n=\C(\pr,\pi(\pr),\pi^2(\pr)\cdots,x)$ is a
Picard-Vessiot extension of
$\pt_{n-1}=\C(\pi(\pr),\pi^2(\pr),\cdots$ $,x)$ with Galois group
$\G:=(\C,+)^m$. Note that the transcendence degree of
$\C(\Ge)=\pt_n$ over $\pt_{n-1}$ is $|\pr|$ since
$\pr\cap\pi^j(\pr)=\emptyset$ for any $1\leq j\leq n$ and
therefore $m=|\pr|$. Clearly $\pt_{n-1}\langle u\rangle$ is an
intermediate differential field. Since
$u\in\Cd_n\setminus\Cd_{n-1}$, we see that $\pt_{n-1}\langle
u\rangle\neq\pt_{n-1}$. Let $\sg\leq\G$ be the group of all
automorphisms that fixes $\pt_{n-1}$ and let
$\{L_i(x_1,\cdots,x_m)|1\leq i\leq t\}$ be the system of
polynomials  for which $\sg$ is the set of solutions. Then it is
easy to see that
\begin{equation}\label{eq struct of diff field}
\pt_{n-1}\langle
u\rangle=\pt_{n-1}(L_i(y_1,\cdots,y_m)),\end{equation} where
$y_j\in \pr.$ Note that $L_i(y_1,\cdots,y_m)'\in\pt_{n-1}$ and
thus $\pt_{n-1}(L_i(y_1,$ $\cdots,y_m))$ is a differential field.

Let $D_i$ be the set of \textsl{essential elements} of
$L_i(y_1,\cdots,y_m)$. Then from equation \ref{eq struct of diff
field} $u\in\C(U)$, where
$U=(\cup^t_{i=1}D_i)\cup(\cup^n_{i=1}\pi^i(E))$. Since $E$ is the
\textsl{essential elements} of $u$, we obtain $\pr\subset E\subset
U$. Now, $\pr\cap(\cup^n_{i=1}\pi^i(E))=\emptyset$ will imply
$\pr\subset\cup^t_{i=1}D_i$. Hence for every $y_j\in \pr$ there is
an $L_i(y_1,\cdots,y_m)$ such that the coefficient of $y_j$ is
nonzero. Let us denote the set $\{L_i(y_1,\cdots,y_m)|1\leq i\leq
t\}$ by $\ls$.

Since $\pt_{n-1}$ is a Picard-Vessiot extension of $\pt_{n-2}$, we
see that $\pt_{n-1}\langle u\rangle$ is a Picard-Vessiot extension
of $\pt_{n-2}\langle u\rangle$. Also,
$L_i(y_1,\cdots,y_m)\in\pt_{n-1}$ for each $i$. Thus from theorem
\ref{trans log in PV} we see that for each $y_j\in \pr$,
$\pi(y_j)\in\pt_{n-2}\langle u\rangle$ and thus
$\pi(\pr)\subset\pt_{n-2}\langle u\rangle$. This shows that
$\pt_{n-1}\langle u\rangle=\pt_{n-2}\langle u\rangle$. Since
$\pt_{n-2}\langle u\rangle$ is a Picard-Vessiot extension of $
\pt_{n-3}\langle u\rangle$, again applying theorem \ref{trans log
in PV} we see that $\pi^2(\pr)\subset\pt_{n-3}\langle u\rangle$
and therefore $\pt_{n-2}\langle u\rangle=\pt_{n-3}\langle
u\rangle$. Thus $\pt_{n-1}\langle u\rangle= \pt_{n-2}\langle
u\rangle=\pt_{n-3}\langle u\rangle$. Assume that
$\pt_{n-(i-1)}\langle u\rangle=\pt_{n-i}\langle u\rangle$. Then
$\pi^{i-1}(\pr)\subset\pt_{n-i}\langle u\rangle$ and therefore
applying theorem \ref{trans log in PV} to the Picard-Vessiot
extension $\pt_{n-i}\langle u\rangle|$ $\pt_{n-(i+1)}\langle
u\rangle$, we see that $\pi^i(\pr)\subset\pt_{n-(i+1)}\langle
u\rangle$. This shows us that $\pt_{n-i}\langle
u\rangle=\pt_{n-(i+1)}\langle u\rangle$. Thus the above induction
argument shows \begin{equation}\label{the fat u}\pt_{n-1}\langle
u\rangle=\C\langle u\rangle\end{equation} and therefore from
equation \ref{eq struct of diff field} we obtain
$$\C\langle u\rangle=\C(\ls,
\pi(\pr),\pi^2(\pr),\cdots,x),$$ where
$\ls=\{L_i(y_1,\cdots,y_m)|1\leq i\leq t\}\subset span_\C \pr$.
\epf

As we noted earlier, $\Ge=\cup^n_{i=0}\pi^i(\pr)$ and therefore
$\pr\subseteq E$ implies $\pi(E)\subset\cup^n_{i=1}\pi^i(\pr)$.
Thus $\cup^n_{i=1}\pi^i(\pr)=\cup^n_{i=1}\pi^i(E)$ and hence we
also have
$$\C\langle u\rangle=\C(\ls,
\pi(E),\pi^2(E),\cdots,x).$$

\begin{remark}
From theorem \ref{singly gen diff subfields of log closure} we
also see that, if $u\in\Cd_n\setminus\Cd_{n-1}$ and $E$ the set of
\textsl{essential elements} of $u$ then
\begin{equation}  \C(\Ge)\supseteq\C\langle u\rangle\supset\pt_{n-1}
\supset\cdots\pt_1\supset\pt_0\supset\C.\end{equation} In
particular, if $u\in\Cd_\infty\setminus\C$ then $x\in\C\langle
u\rangle$.
\end{remark}

Now we will generalize theorem \ref{singly gen diff subfields of
log closure} to any finitely generated differential subfield of
$\Cd_n$.

\bt \label{struct of fin gen dif fields of log closure} Let
$\mathbf{K}:=\C\langle u_1,\cdots,u_m\rangle$ be a finitely
differentially generated subfield of $\Cd_n\setminus\Cd_{n-1}$ and
let $E:=\cup^m_{i=1}E_i$, where $E_i$ is the set of
\textsl{essential elements} of $u_i$. For each $i$, let $n_i\in\N$
be minimal such that $E_i\subset\cup^{n_i}_{j=0}\Lambda_j$ and let
$\pr_i\subset E_i$ be the $\pi-$base of
$\Ge_i:=\cup^{n_i}_{j=0}\pi^j(E_i)$. Then there are finite sets
$\ls_i\subset span_\C \pr_i$ such that
$$\mathbf{K}=\C(\ls,\pi(\pr),\pi^2(\pr),\cdots,x),$$ where $\ls=
\cup^m_{i=1}\ls_i$ and $\pr=\cup^m_{i=1}\pr_i$. Moreover, for
every $y\in\pr$, $\ls$ contains at least one linear combination in
which $y$ appears nontrivially. \et

\bpf Since $\mathbf{K}$ is a compositum of singly generated
differential fields, the proof follows from theorem \ref{singly
gen diff subfields of log closure}. \epf

\bt \label{fin gen imp single gen} Every finitely generated
differential subfield of $\Cd_\infty$ is singly generated.  \et

\bpf Let $\mathbf{K}$ be a finitely generated differential
subfield of $\Cd_n\setminus\Cd_{n-1}$. Then from theorem
\ref{struct of fin gen dif fields of log closure} there are sets
$\ls$ and $\pr$ such that
$$\mathbf{K}=\C(\ls,\pi(\pr),\pi^2(\pr),\cdots,x).$$  Let
$\ls=\{L_i|1\leq i\leq m\}$, $u=\sum^n_{i=1}x^iL_i$,
$\pv:=\C(\pr,\pi(\pr),\pi^2(\pr),\cdots,x)$ and let
$\f:=\C(\pi(\pr),\pi^2(\pr),\cdots,x)$. We see that $\pv|\f$ is a
Picard-Vessiot extension(antiderivative extension), and since
$L_i\in span_\C\pr$ we obtain $L'_i\in\f$ and thus $\mathbf{K}$ is
an intermediate Picard-Vessiot sub-extension of $\pv|\f$. Consider
the Picard-Vessiot extension $\mathbf{K}|\f$. Since
$\f(\ls)$=$\mathbf{K}$ is an antiderivative extension of $\f$ and
$u\in\mathbf{K}$, we see that for any $\Phi\in\G(\mathbf{K}|\f)$
\begin{align*}
\Phi(u)&=\sum^n_{i=1}x^i\Phi(L_i)\\
&=\sum^n_{i=1}x^i(L_i+c_i)\\
&=\sum^n_{i=1}x^iL_i+\sum^n_{i=1}c_ix^i\\
&=u+\sum^n_{i=1}c_ix^i,
\end{align*}
where $c_i\in\C$. Thus if $\Phi$ fixes $u$, we obtain
$\sum^n_{i=1}c_ix^i=0$ and
 therefore $\Phi$ has to be the identity. Thus $\f\langle
 u\rangle=\mathbf{K}$. Consider $$\frac{\partial u}{\partial y}=
 \sum^n_{i=1}x^i\frac{\partial L_i}{\partial y}.$$ We observe from
 theorem \ref{struct of fin gen dif fields of log closure} that
 for $y\in\pr$ there is an $i$ such that
 $\frac{\partial L_i}{\partial y}\neq 0$, and we also recall that $\pr\cup\{x\}$
 is algebraically independent over $\C$. Thus $\frac{\partial u}{\partial y}\neq
 0$ for any $y\in\pr$
and we also obtain that $E:=\pr\cup\{x\}$ is the set of
\textsl{essential elements} of $u$. It can be easily seen that the
$\pi-$base of $\Ge:=\cup^n_{i=0}\pi^i(E)$ is again $\pr$ and
therefore applying theorem \ref{singly gen diff subfields of log
closure}, we see that
$\pi(\pr),\pi^2(\pr),\cdots,x\subset\C\langle u\rangle$. Thus
$\f\subset\C\langle u\rangle$ and therefore ${\bf K}=\f\langle
u\rangle=\C\langle u\rangle$ and we are done. \epf

{\bf An Algorithm to Compute the Differential field $\C\langle
u\rangle$}

\bt \label{the alg theo for iter log} Let
$u\in\Cd_n\setminus\Cd_{n-1}$ and let $P,Q\in\C[E]$, where $E$ is
the set of \textsl{essential elements} of $u$, $(P,Q)=1$ and
$u=\frac{P}{Q}$. Then the set $\ls$ and $\pr$ from theorem
\ref{singly gen diff subfields of log closure} can be computed
from $P$ and $Q$. \et

\bpf

Since $\pr=E\setminus\cup^n_{i=1}\pi^i(E)$, we see that the set
$\pr$ can be computed once the set $E$ of essential elements is
known. From equation \ref{the fat u} we see that
$\pi(\pr),\pi^2(\pr)\cdots$$\pi^n(\pr)=\{x\}$$\subset\C\langle
u\rangle$. That is $\pt_{n-1}\subset\C\langle u\rangle$ and thus
$\C\langle u\rangle$ is an intermediate differential field of the
Picard-Vessiot extension $\C(\Ge)|\pt_{n-1}$. That is
\begin{equation}\label{trap u in a PV}\C(\Ge)\supseteq\C\langle u\rangle\supset\pt_{n-1}.\end{equation}
Also note that $\C(\Ge)$ is an extension by antiderivatives of
$\pt_{n-1}$ and that $\C(\Ge)=\pt_{n-1}(\pr)$ and
$\pr\cap\pt_{n-1}=\emptyset$ since $\Ge$ is algebraically
independent over $\C$. Thus $\C(\Ge)|\pt_{n-1}$ is a pure
transcendental extension of transcendence degree $|\pr|$. Now we
may apply theorem \ref{struct theorem for antideriv} to obtain the
set $\ls$. Thus from equation \ref{the fat u}, we see that
$\C\langle u\rangle=\pt_{n-1}(\ls)$.\epf

{\bf Algorithm:} Write out two polynomial expressions, say $A, B$,
over $\C$ with elements from $\Lambda_\infty$ as indeterminates.
The following steps will find the differential field $\C\langle
u\rangle$, where $u=\frac{A}{B}$, in the form of a finitely
generated field expressed in theorem \ref{singly gen diff
subfields of log closure}.

\begin{description}
  \item[Step 0]
First we form a finite set $S$ by picking elements from
$\Lambda_\infty$ that appear in the expression of $A$ or $B$. Then
compute the set $E$ of essential elements of $u$. That is, find
the set
$$E:=\Big\{y\in S\Big|\dfrac{\partial P}{\partial y}\neq 0\ or\
\dfrac{\partial Q}{\partial y}\neq 0\Big\}.$$ Also find the set
$\pr=E\setminus\cup^n_{i=1}\pi^i(E)$, where $n$ is the least
positive integer such that $\pi^n(E)=\{x\}$ and let $\Ge:=
\cup^n_{i=0}\pi^i(E)$.\\

\item[Step 1] From equation \ref{trap u in a PV}, we obtain
$\pt_{n-1}\subset\C\langle u\rangle$. In particular
$\pi(\pr),\pi^2(\pr)$ $\cdots\pi^n(\pr)=\{x\}$ $\subset\C\langle
u\rangle$. Since $\C(\Ge)$ is an antiderivative
    extension of $\pt_{n-1}$, we obtain that $\C\langle u\rangle$ is
    an intermediate differential subfield of the Picard-Vessiot
    extension $\C(\Ge)$ of $\pt_{n-1}$. \\

 \item[Step 2]
We replace $A, B$ by some $P,Q\in\C[E]$ such that $(P,Q)=1$. This
can be done in two ways. We may use {\textsl MATHEMATICA 5.2} and
compute the $GCD$ of $A,B$ and divide $A,B$ by the $GCD$ to get
$P,Q$ such that $\frac{A}{B}=\frac{P}{Q}$ and $GCD$ of $P,Q$ is 1.
In case, when \textsl{MATHEMATICA} 5.2 fails to compute the $GCD$,
we way  compute the Gr\"{o}bner basis \cite{Ad-Lou} for the Ideal
$<A,B>$ generated over $\C[S]$ and use {\textsf Gaydar's formula}
\cite{Duz-Chm} to compute the $GCD$ and then use the multivariable
division algorithm \cite{Ad-Lou} to find out $P,Q$ such that
$\frac{A}{B}=\frac{P}{Q}$ and $GCD$ of $P,Q$ is 1.

Thus we note that finding a relatively prime polynomials for a
given pair of polynomial from $\C[\Lambda_\infty]$ is a finite
process.

Now we have $u=\frac{P}{Q}$, $P,Q\in\C[E]$ and $(P,Q)=1$.\\

     \item[Step 3] Write $P$ and $Q$ as polynomials over $\R:=$ $\C[\pi(\pr),$ $\pi^2(\pr)$,$\cdots$,$x]$ with
    elements of $\pr$ as variables. Then $\pt_{n-1}$ becomes the
    fraction field of $\R$. Note that $\C(\Ge)|\pt_{n-1}$ is a Picard-Vessiot extension( by
    antiderivatives) of transcendence degree $p:=|\pr|$
    and thus if  $\sigma\in\G(\Ge|\mathbf{K})$ then
    $\sigma(P)=P(y_1+c_{1\sigma},\cdots,y_p+c_{p\sigma})$ and $\sigma(Q)=Q(y_1+d_{1\sigma},\cdots,y_p+d_{p\sigma})$
where $c_{i\sigma}$, $d_{j\sigma}\in\C$ and $y_i\in\pr$. Also from
theorem \ref{struct theorem for antideriv},
we see that $\sigma(u)=u$ if and only if $\sigma(P)=P$ and $\sigma(Q)=Q$.\\
    \item[Step 4] From proposition \ref{trans poly prop} we obtain that if $\sigma$ fixes $P$ and $Q$ then it
fixes each of the homogeneous components of $P$ and $Q$ and from
this fact (following the proof of proposition \ref{trans poly
prop}) we obtain linear forms over $\R$ such that the field
generated by $\pt_{n-1}$ and the linear forms equals the field
$\C\langle u\rangle$. Thus, we compute a system of linear forms
$\{D_j\}$ over $\R$ such that $\sigma(P)=P$ and $\sigma(Q)=Q$ if
and only if
    $D_j(c_{1\sigma},\cdots,c_{p\sigma})=0$.  \\

    \item[Step 5] Since $\R$ is a polynomial ring, using proposition
     \ref{reduction of poly to c poly}, we
    could compute a system of linear forms $\{L_j\}$ over $\C$
from the system $\{D_j\}$ such that the set of solutions of $L_j$ and $D_j$ over $\C^p$ are the same.\\
    \item[Step 6] Finally, from theorem \ref{singly gen diff subfields of log closure}
we see that the field $$\C\langle
u\rangle=\C(\ls,\pi(\pr),\pi^2(\pr),\cdots,x),$$ where
$\ls=\{L_j(y_1,\cdots,y_p)|y_i\in\pr\}$.
\end{description}

\section{Examples}\label{examples}
In this section we will apply our algorithm to compute the
differential fields generated by an element of $\Cd_\infty$ and
$\C$. Also we assume $\C:=\mathbb{C}$, the field of complex
numbers.

{\bf Example 1} Consider the field $\Cd_1$ and  Let
$$u=\frac{5x^3\ln(x+1)+\ln(x+e)+27x^3\ln(x+\sqrt{2})}{\ln(x)+x\big(\ln(x+2)-17\ln(x+3)\big)^2}\in\Cd_1.$$

\begin{description}
    \item[Step 0] Let $A:=5x^3\ln(x+1)+\ln(x+e)+27x^3\ln(x+\sqrt{2})$
    and $B:=\ln(x)+x\big(\ln(x+2)-17\ln(x+3)\big)^2$.
     We observe that $u\in\C(S)$, where $S=\{x, \ln(x), \ln(x+1),
    \ln(x+2), \ln(x+3), \ln(x+e), \ln(x+\sqrt{2})\}$. We easily see that the essential elements $E$
    equals the set $S$. The set $\Ge=\cup^1_{i=0}\pi(E)$ and in
    this case, we see that $\Ge=E$. The $\pi-$base of $\pr$ of
    $\Ge$ is the set $\pr=\{\ln(x), \ln(x+1),
    \ln(x+2), log(x+3), \ln(x+e), \ln(x+\sqrt{2})\}$.\\
    \item[Step 1] Since $u\in\Cd_1$, we have $n=1$ and thus $\C(\Ge)\supseteq\C\langle
    u\rangle\supset\pt_0=\C(x)$. The differential field $\C(\Ge)$ is an antiderivative
    extension of $\C(x)$ and therefore $\C\langle u\rangle$ is
    an intermediate differential subfield of the Picard-Vessiot
    extension $\C(\Ge)$ of $\C(x)$.
    \item[Step 2]We note that $A$ and $B$ are relatively prime
    and thus we may choose $P:=A$ and $Q:=B$.
    \item[Step 3]We rewrite $P$ and $Q$ as polynomials over $\R:=\C[x]$. Then
    $P=x^3\big(5\ln(x+1)+27\ln(x+\sqrt{2})\big)+\ln(x+e)$ and
    $Q=\ln(x)+x\big(\ln(x+2)-17\ln(x+3)\big)^2$. Let
    $y_1:=\ln(x+1)$, $y_2:=\ln(x+\sqrt{2})$,
    $y_3:=\ln(x+e)$, $y_4:=\ln(x)$, $y_5:=\ln(x+2)$ and $y_6:=\ln(x+3)$.
     We observe that if $\sigma\in
    \G(\C(\Ge)|\pt_0)$, then $\sigma(y_i)=y_i+c_{i\sigma}$ for
    each $y_i\in\pr$ and we also observe that for any $\sigma\in\G(\C(\Ge)
    |\C(x))$, $\sigma(u)=u$ if and only if $\sigma(P)=P$ and
    $\sigma(Q)=Q$.
    \item[Step 4]Note that $P$ is a homogeneous polynomial of
total degree 1 over $\C[x]$. If $\sigma$ fixes $P$ then
\begin{align*}&\sigma(P)=P\\
\iff&\sum^3_{i=1}c_{i\sigma}\frac{\partial P}{\partial y_i}=0\\
\iff&x^3\big(5c_{1\sigma}+27c_{2\sigma}\big)+c_{3\sigma}=0.
\end{align*}
Let $D_1:=x^3\big(5y_1+27y_2\big)+y_3$. Then we see that for any
$\sigma\in\G(\C(\Ge)|\C(x))$, $\sigma(D_1)=D_1$ if and only if
$x^3\big(5c_{1\sigma}+27c_{2\sigma}\big)+c_{3\sigma}=0$.

If $\sigma$ fixes $Q$ then $\sigma$ fixes the homogeneous
components of $Q$ and thus $\sigma$ fixes $y_4:=\ln(x)$ and
$x\big(y_5-17y_6\big)^2$. Now
\begin{align*}&\sigma(x(y_5-17y_6)^2)=x(y_5-17y_6)^2\\
\iff&\sum^6_{i=5}c_{i\sigma}\frac{\partial Q}{\partial y_i}=0\\
\iff&x(c_{5\sigma}-17c_{6\sigma})(y_{5}-17y_{6})=0\\
\iff&c_{5\sigma}-17c_{6\sigma}=0.
\end{align*}
Let $D_2:=y_5-17y_6$. Then  for any $\sigma\in\G(\C(\Ge)|\C(x))$,
$\sigma(D_2)=D_2$ if and only if $c_{5\sigma}-17c_{6\sigma}=0$.

\item[Step 5]Note that
$x^3\big(5c_{1\sigma}+27c_{2\sigma}\big)+c_{3\sigma}=0$ if and
only if $c_{3\sigma}=0$ and $5c_{1\sigma}+27c_{2\sigma}=0$. That
is, $\sigma$ fixes $P$ if and only if it fixes $y_3$ and
$5y_1+27y_2$. We also observe that the linear form $D_2$ is
already over $\C$.

Thus we have proved that for any $\sigma\in\G(\C(\Ge)|\C(x))$,
$\sigma$ fixes $u$ if and only if $\sigma$ fixes $x, y_3, y_4,
5y_1+27y_2$ and $y_5-17y_6$.
    \item[Step 6]
$$\C\langle u\rangle=\C(x, \ln(x+e), \ln(x),
5\ln(x+1)+27\ln(x+\sqrt{2}), \ln(x+2)-17\ln(x+3))$$
\end{description} $\square$

{\bf Example 2}

Let $y_1:=\ln(\ln(\ln(x-i)+2)+3)$, $y_2:=\ln(\ln(x+i)+\sqrt{3})$,
$y_3:=\ln(x+\frac{5}{6})$,
$y_4:=\ln(\ln(x+\frac{1}{2})+\frac{1}{2})$,
$y_5:=\ln(x+\sqrt{5})$, $y_6:=\ln(x+5+i)$,
$y_7:=\ln(\ln(\ln(x)+i))$ and let
$$u=\frac{\ln(x+i)^2\ln(x-i)(y_1-y_3)^5+x^3\ln(x)(y_2-y_5)^2}
{\ln(\ln(x)+i)^2(y_5-y_7)^7+x\ln(x-i)^3\ln(\ln(x-i)+2)^2(y_6-y_4)^{12}}\in\Cd_3.$$
We will apply the algorithm to compute the differential field
generated by $\C$ and $u$.

\begin{description}
    \item[Step 0] Let
    $A:=\ln(x+i)^2\ln(x-i)(y_1-y_3)^5+x^3\ln(x)(y_2-y_5)^2$,
    $B:=\ln(\ln(x)+i)^2(y_5-y_7)^7+x\ln(x-i)^3\ln(\ln(x-i)+2)^2(y_6-y_4)^{12}$
    and $S:=\{y_1, y_2, y_3, y_4, y_5, y_6, y_7, \ln(x-i),
    \ln(x+i),\ln(\ln(x)+i), \ln(x), x,$ $\ln(\ln(x-i)+2)\}$. We observe that
    the set of essential elements $E$ of $u$ equals the set $S$.
    Since $\pi(E)=\{\ln(x+\frac{1}{2}), \ln(\ln(x)+i), \ln(\ln(x-i)+2),
    \ln(x+i)\}$, $\pi^2(E):=\{\ln(x-i),x,\ln(x)\}$ and
    $\pi^3(E)=\{x\}$, we see that
    $\Ge=\cup^3_{i=0}\pi^i(E)$ $=E\cup\{\ln(x+\frac{1}{2})\}$.
    Then the $\pi-$base $\pr$ of $E$ is the set
    $E\setminus\cup^3_{i=1}\pi^i(E)$ $=\{y_1,y_2,\cdots,y_7\}$.
    \item[Step 1] We know that $\cup^3_{i=1}\pi^i(\pr)$ $=\{\ln(x-i), \ln(\ln(x-i)+2), \ln(x+i), \ln(x+\frac{1}{2}),
    \ln(\ln(x)+i), \ln(x), x\}$ and that $\pt_2=\C(\cup^3_{i=1}$ $\pi^i(\pr))\subset\C\langle
    u\rangle$. Thus $\C\langle u\rangle$ is an intermediate
    subfield of the Picard-Vessiot extension(antiderivative extension) $\pt_3:=\C(\Ge)$ of
    $\pt_2$. Also note that $\pt_3=\pt_2(y_1,y_2,\cdots,y_7)$.

    \item[Step 2] One can easily see that $A$ and $B$ are
    relatively prime and thus choose $P:=A$ and $Q:=B$.
    \item[Step 3] The polynomials $P$ and $Q$ are already
    presented as polynomials over the field
    $\C(\cup^3_{i=1}\pi^i(\pr))$ with $y_1,y_2,\cdots,y_7$ as
    variables. We note that if $\sigma\in
    \G(\C(\Ge)|\pt_{2})$, then $\sigma(y_i)=y_i+c_{i\sigma}$ for
    each $y_i\in\pr$ and we also observe that for any $\sigma\in\G(\C(\Ge)
    |\pt_2)$ such that $\sigma(u)=u$ then $P$ divides $\sigma(P)$ and $Q$ divides $\sigma(Q)$. Then from
    proposition \ref{trans poly prop} we have $\sigma(u)=u$ if and only if $\sigma(P)=P$ and
    $\sigma(Q)=Q$.
    \item[Step 4] Let $\sigma=(c_{1\sigma},\cdots,c_{7\sigma})\in\G(\C(\Ge)
    |\pt_2)$ be an automorphism such that $\sigma(u)=u$. Then
    $\sigma(P)=P$ and $\sigma(Q)=Q$ and now we shall use proposition
    \ref{trans poly prop} to compute the linear forms. Note that
    $\sigma$ fixes $u$ if and only if it fixes $H_8:=\ln(x+i)^2\ln(x-i)(y_1-y_3)^5$, $H_6=x^3\ln(x)(y_2-y_5)^2$
    $H_{18}=x\ln(x-i)^3\ln(\ln(x-i)+2)^2(y_6-y_4)^{12}$ and
    $H_9=\ln(\ln(x))^2(y_5-y_7)^7$. Thus $\sum^7_{i=1}c_{i\sigma}\frac{\partial
    H_j}{\partial
    y_i}=0$ for $j=6, 8, 9$ and $18$, which gives us the following
    equations
    \begin{align*}&\ln(x+i)^2\ln(x-i)(c_{1\sigma}-c_{3\sigma})=0,\\
    &x^3\ln(x)(c_{2\sigma}-c_{5\sigma})=0,\\ &\ln(\ln(x)+i)^2(c_{5\sigma}-c_{7\sigma})=0,\\
    &x\ln(x-i)^3\ln(\ln(x-i)+2)^2(c_{6\sigma}-c_{4\sigma})=0.\end{align*}
We also observe that the $\pt_2-$linear forms of the field
$\C\langle u\rangle$ are $H_j$, $j=6,8,9$ and $16$. That is
$\C\langle u\rangle=\pt_2(H_6, H_8, H_9, H_{18})$.

    \item[Step 5] From the above displayed equations, it is clear
    that $\sigma(u)=u$ if and only if $c_{1\sigma}-c_{3\sigma}=0$,
    $c_{2\sigma}-c_{5\sigma}=0$, $c_{5\sigma}-c_{7\sigma}=0$ and
    $c_{6\sigma}-c_{4\sigma}=0$.
\item[Step 6] \begin{align*}\C\langle u\rangle=&\C(\ln(x-i),
    \ln(x+i),\ln(\ln(x)+i), \ln(x+\frac{1}{2}), \ln(x), x,\\ &\ln(\ln(x-i)+2), y_1-y_3, y_2-y_5, y_6-y_4,y_5-y_7).\end{align*}
\end{description}


\end{document}